\documentclass[a4paper,12pt]{article}     

\usepackage{graphicx}		  % to deal with graphicals
\usepackage{amsmath,amssymb}	% to deal with mathematics

\newtheorem{lemma}{Lemma}[section]
\newtheorem{theorem}[lemma]{Theorem}
\newtheorem{proposition}[lemma]{Proposition}
\newtheorem{corollary}[lemma]{Corollary}
\newtheorem{assumption}{Assumption}
\newtheorem{remark}[lemma]{Remark}

\def\authorfont{\footnotesize}

\def\ccode#1{\par
\vspace*{8pt}
{\authorfont{\leftskip18pt\rightskip\leftskip
\noindent #1\par}}\par}

\newenvironment{Proof}{
\hspace*{-9mm}
{ \it Proof.}}
{\hfill {$\square$}\vspace{1.5em}}

\begin{document}

\begin{center}{
{\Large 
 Properties of minimal charts and
 their applications IX: 
 charts of type $(4,3)$}
\vspace{10pt}
\\ 
Teruo NAGASE and Akiko SHIMA\footnote{The second author is supported by JSPS KAKENHI Grant Number 21K03255.}
}
\end{center}

%%%%%%%%%%%%%%%%%%%%%%%%%%%%%%%
%%%%   abstract   %%%%%%%%%%%%%
%%%%%%%%%%%%%%%%%%%%%%%%%%%%%%%
\begin{abstract}
Charts are oriented labeled graphs in a disk.
Any simple surface braid (2-dimensonal braid) can be described by using a chart.
Also, a chart represents an oriented closed surface
embedded in 4-space.
In this paper, we investigate embedded surfaces in 4-space
by using charts.
Let $\Gamma$ be a chart,
and we denote by $\Gamma_m$
the union of all the edges of label $m$.
A chart $\Gamma$ is of type $(4,3)$
if there exists a label $m$
such that 
$w(\Gamma)=7$,
$w(\Gamma_m\cap\Gamma_{m+1})=4$,
$w(\Gamma_{m+1}\cap\Gamma_{m+2})=3$
where 
$w(G)$ is the number of white vertices in $G$.
In this paper, we prove that there is 
no minimal chart of 
type $(4,3)$.
\end{abstract}

%
%
%%%%   Make Content   %%%%%%%%%%
%%%%%%%%%%%%%%%%%%%%%%%%%%%%%%%%
%\tableofcontents
%
%
%%%%%%%%%%%%%%%%%%%%%%%%%%%%
%%%%%%%%%%%%%%%%%%%%%%%%%%%%

\ccode{2020 Mathematics Subject Classification. Primary 57K45,05C10; Secondary 57M15.}
\ccode{ {\it Key Words and Phrases}. surface link, chart, C-move, white vertex. }

%05C10 Planar graphs; geometric and topological aspects of graph theory
%57M15 Relations of low-dimensional topology with graph theory
%57K45 Higher-dimensional knots and links.

\setcounter{section}{0}
\section{Introduction}

%\large
%\baselineskip 20pt

Charts are oriented labeled graphs in a disk (see  \cite{KnottedSurfaces},\cite{BraidBook}, and see Section~\ref{s:Prel}  for the precise definition of charts).
Let $D_1^2, D_2^2$ be 2-dimensional disks.
Any simple surface braid (2-dimensoinal braid) can be described 
by using a chart,
here a simple surface braid is a properly embedded surface
$S$ in the 4-dimensional disk $D_1^2\times D_2^2$ such that
a natural map $\pi:S\subset D_1^2\times D_2^2\to D_2^2$ is 
a simple branched covering map of $D_2^2$ and
the boundary $\partial S$ is a trivial closed braid in
the solid torus $D_1^2\times \partial D_2^2$
(see \cite{BraidThree}, \cite[Chapter 14 and Chapter 18]{BraidBook}).
Also, from a chart, 
we can construct a simple closed surface braid in 4-space ${\Bbb R}^4$. This surface is an oriented closed surface 
embedded in ${\Bbb R}^4$.
On the other hand, any oriented embedded closed surface 
 in ${\Bbb R}^4$ is ambient isotopic to a simple
closed surface braid
 (see \cite{BraidThree},\cite[Chapter 23]{BraidBook}). 
A C-move 
is a local modification between two charts
in a disk (see Section~\ref{s:Prel} for C-moves).
A C-move between two charts induces 
an ambient isotopy between oriented closed surfaces 
corresponding to the two charts.
In this paper, we investigate oriented closed surfaces in 4-space
by using charts.

We will work in the PL category or smooth category. All submanifolds are assumed to be locally flat.
In \cite{ONS},
we showed that there is no minimal chart with exactly five vertices
 (see Section~\ref{s:Prel} for the precise definition of minimal charts). 
Hasegawa proved that there exists a minimal chart with exactly
six white vertices \cite{H1}. 
This chart represents a 2-twist spun trefoil.
In \cite{INS} and \cite{NST},
we investigated minimal charts with exactly four white vertices.
In this paper, 
we investigate properties of minimal charts 
which support a conjecture that
there is no minimal chart with exactly seven white vertices
(see \cite{ChartApp1},\cite{ChartAppII},
\cite{ChartAppIII},\cite{ChartAppIV},
\cite{ChartAppV}, \cite{ChartAppVI},
\cite{ChartAppVII},\cite{ChartAppVIII},\cite{ChartAppX}).
This paper is the second paper from the last in this series.

Let $\Gamma$ be a chart.
For each label $m$, we denote by $\Gamma_m$
the union of all the edges of label $m$.

Now we define a type of a chart:
Let $\Gamma$ be a chart with at least one white vertex, 
and $n_1,n_2,\dots,n_k$ integers.
The chart $\Gamma$ is of {\it type $(n_1,n_2,\dots,n_k)$} if there exists a label $m$ of $\Gamma$ satisfying the following three conditions:
\begin{enumerate}
\item[(i)] For each $i=1,2,\dots, k$, 
the chart $\Gamma$ contains exactly $n_{i}$ white vertices in $\Gamma_{m+i-1}\cap \Gamma_{m+i}$.
\item[(ii)] If $i<0$ or $i>k$, then $\Gamma_{m+i}$ does not contain any white vertices.
\item[(iii)] Both of the two subgraphs $\Gamma_m$ and $\Gamma_{m+k}$ contain at least one white vertex.
\end{enumerate}
If we want to emphasize the label $m$,
then we say that $\Gamma$ is of {\it type $(m;n_1,n_2,\dots,n_k)$}. 
Note that $n_1\ge1$ and $n_k\ge1$ by Condition~(iii).

We proved in \cite[Theorem 1.1]{ChartAppII} that
if there exists a minimal $n$-chart $\Gamma$ with exactly seven white vertices,
then $\Gamma$ is a chart of 
type $(7),(5,2),(4,3),(3,2,2)$ or $(2,3,2)$ 
(if necessary we change the label
$i$ by $n-i$ for all label $i$).
In \cite{ChartAppV},
we showed that
there is no minimal chart of type $(3,2,2)$.
In \cite{ChartAppVI} and \cite{ChartAppVII},
there is no minimal chart of type $(2,3,2)$.
In \cite{ChartAppVIII},
there is no minimal chart of type $(7)$.
In this paper we shall show the following:

\begin{theorem}
\label{MainTheorem} 
There is 
no minimal chart of 
type $(4,3)$.
\end{theorem}

In the future paper \cite{ChartAppX},
we shall show there is no minimal chart of type
$(5,2)$.
Therefore we shall show that
there is no minimal chart with exactly seven white vertices.

The paper is organized as follows.
In Section~\ref{s:Prel},
we define charts and minimal charts.
Let $\Gamma$ be a minimal chart, and $m$ a label of $\Gamma$. 
In Section~\ref{s:DiskTwo},
we investigate
a disk $D$ with exactly two white vertices 
of $\Gamma_m$
such that
$\Gamma_m\cap\partial D$
consists of at most two points.
In Section~\ref{s:3AngledDiskOneWhiteVertex},
 we investigate a disk called 
a 3-angled disk of $\Gamma_m$ with
at most one white vertex in its interior, 
where for a positive integer $k$, 
a $k$-angled disk is a disk 
whose boundary contains exactly $k$ white vertices
and consists of edges of label $m$.
In Section~\ref{s:2-angledDisks},
we investigate a  2-angled disk of $\Gamma_m$.
In Section~\ref{s:3-angledDisksTwoWhiteNoFeelers},
 we investigate a 3-angled disk $D$ of $\Gamma_m$ 
without feelers,
where a feeler of $D$ is an edge $e$ of label $m$ in $D$
with  $e\cap\partial D\not=\emptyset$ and
 $e\not\subset \partial D$.
In Section~\ref{s:3-angledDisksOneFeeler},
 we investigate a 3-angled disk with exactly one feeler 
such that the feeler contains a black vertex.
In Section~\ref{s:IOC},
we review IO-Calculation(a property of numbers of 
inward arcs of label $m$ 
and outward arcs of label $m$ in a closed domain $F$
with $\partial F\subset\Gamma_{m-1}\cup\Gamma_m\cup\Gamma_{m+1}$
for some label $m$).
In Section~\ref{s:DiskLemma}, we review New Disk Lemma.
We shall give a condition of non-minimal charts.
Let $\Gamma$ be a minimal chart, and $k$ a label of $\Gamma$. In Section~\ref{s:ProperEdge},
we investigate a $4$-angled disk of $\Gamma_k$
with a proper edge of label $k+\delta$
for some $\delta\in\{+1,-1\}$.
In Section~\ref{s:ProperEdge5AngledDisk}, 
we investigate a $5$-angled disk of $\Gamma_k$ 
with a proper edge of label $k+\delta$
for some $\delta\in\{+1,-1\}$.
In Section~\ref{s:ProperEdge7AngledDisk}, 
we investigate a $7$-angled disk of $\Gamma_k$ 
with a proper edge of label $k+\delta$
for some $\delta\in\{+1,-1\}$.
From Section~\ref{s:Chart43} to Section~\ref{s:D1OneD2Three},
we shall prove Theorem~\ref{MainTheorem}.

%%%%%%%%%%%%%%%%%%%%%%%%%%%%%%%%%%%%%%
%%%%%%%%%%%%%%%%%%%%%%%%%%%%%%%%%%%%%%
%%%%%%%%%%%%%%%%%%%%%%%%%%%%%%%%%%%%%%
%%%%%%%%%%%%%%%%%%%%%%%%%%%%%%%%%%%%%%
%%%%%%%%%%%%%%%%%%%%%%%%%%%%%%%%%%%%%%
%%%%%%%%%%%%%%%%%%%%%%%%%%%%%%%%%%%%%%

%\newpage
\section{Preliminaries}
\label{s:Prel}

In this section, 
we introduce 
the definition of charts and its related words.

\begin{figure}[htb]
\begin{center}
\includegraphics{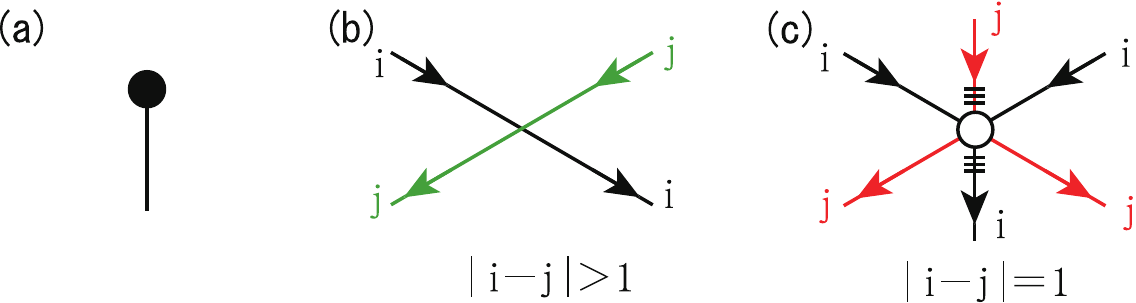}
\end{center}
\caption{ \label{fig01} (a) A black vertex. (b) A crossing. (c) A white vertex. 
Each arc with three transversal short arcs is a middle arc at the white vertex. }
\end{figure}

Let $n$ be a positive integer.
An $n$-{\it chart}  
(a braid chart of degree $n$ \cite{KnottedSurfaces}
or a surface braid chart of degree $n$ \cite{BraidBook}) 
is 
an oriented labeled graph in the interior of a disk,
which may be empty 
or
have closed edges without vertices
satisfying the following four conditions
(see Fig.~\ref{fig01}):
\begin{enumerate}
\item[(i)] 
Every vertex has degree $1$, $4$, or $6$.
\item[(ii)] 
The labels of edges are 
in $\{1,2,\dots,n-1\}$.
\item[(iii)]
In a small neighborhood of
each vertex of degree $6$,
there are six short arcs,
three consecutive arcs are
oriented inward 
and
the other three are outward,
and
these six are labeled $i$ and $i+1$
alternately for some $i$,
where the orientation and label of
each arc are inherited from
the edge containing the arc.
\item[(iv)]
For each vertex of degree $4$,
diagonal edges have the same label
and
are oriented coherently,
and the labels $i$ and $j$ of
the diagonals satisfy $|i-j|>1$.
\end{enumerate}
We call a vertex of degree $1$ a {\it black vertex},
a vertex of degree $4$ a {\it crossing}, and 
a vertex of degree $6$ a {\it white vertex}
respectively.

Among six short arcs
in a small neighborhood of
a white vertex,
a central arc of each three consecutive arcs
oriented inward (resp. outward) 
is called a   
{\it middle arc} at the white vertex
(see Fig.~\ref{fig01}(c)).
For each white vertex $v$, 
there are two middle arcs at $v$ 
in a small neighborhood of $v$.
An edge is said to be {\it middle at} a white vertex $v$ if it contains a middle arc at $v$.

Let $e$ be an edge connecting $v_1$ and $v_2$.
If $e$ is oriented from $v_1$ to $v_2$,
then we say that 
$e$ is oriented {\it outward at $v_1$}
and {\it inward at $v_2$}.
%%%%%%%%%%%%%%%%%%

%%%%%%%%%%%%%%%%%%

\begin{figure}
\begin{center}
\includegraphics{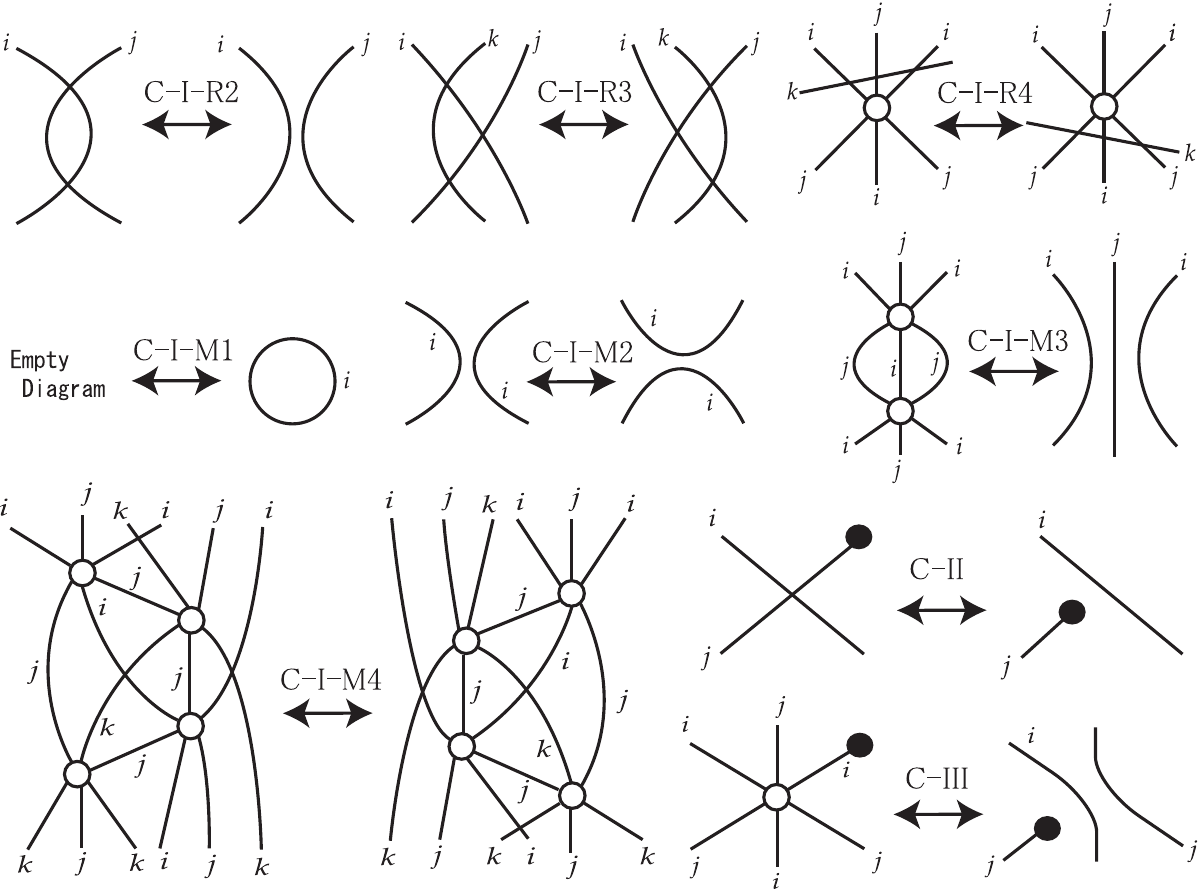}
\end{center}
\caption{ \label{fig02} For the C-III move, 
the edge with the black vertex is not middle at
a white vertex in the left figure. }
\end{figure}

Now {\it C-moves} are local modifications 
of charts as shown in Fig.~\ref{fig02}
(cf. \cite{KnottedSurfaces}, 
\cite{BraidBook} and \cite{Tanaka}).
Two charts are said to be {\it C-move equivalent}  if there exists
a finite sequence of C-moves 
which modifies one of the two charts 
to the other.

%%%%%%%%%%%%%%%%%%%
%%%%%%%%%%%%%%%%%%

An edge in a chart is called 
a {\it free edge}
if it has
two black vertices.

For each chart $\Gamma$,
let $w(\Gamma)$ and $f(\Gamma)$ be the number of white vertices, and the number of free edges respectively.
The pair $(w(\Gamma), -f(\Gamma))$ is called a {\it complexity} of the chart (see \cite{BraidThree}).
A chart $\Gamma$ is called a {\it minimal chart} if its complexity is minimal among the charts C-move equivalent to the chart $\Gamma$ with respect to the lexicographic order of pairs of integers.

We showed the difference of a chart in a disk and in a 2-sphere (see \cite[Lemma 2.1]{ChartApp1}).
This lemma follows from that there exists a natural one-to-one correspondence between $\{$charts in $S^2\}/$C-moves and $\{$charts in $D^2\}/$C-moves, conjugations
(\cite[Chapter 23 and Chapter 25]{BraidBook}).
To make the argument simple, we assume that 
the charts lie on the 2-sphere instead of the disk.
\begin{assumption}
In this paper,
all charts are contained in the $2$-sphere $S^2$.
\end{assumption}
We have the special point in the 2-sphere $S^2$, called the point at infinity,
 denoted by $\infty$.
In this paper, all charts are contained in a disk such that the disk 
does not contain the point at infinity $\infty$.

An edge in a chart is called 
a {\it terminal edge}
if it has
a white vertex and a black vertex.

Let $\Gamma$ be a chart,
and $m$ a label of $\Gamma$. 
A {\it hoop} is a closed edge of $\Gamma$ without vertices 
(hence without crossings, neither).
A {\it ring} is a simple closed curve in $\Gamma_m$ containing at least one crossing but not containing any white vertices.
A hoop is said to be {\it simple} 
if one of the two complementary domains
of the hoop
does not contain any white vertices.

We can assume that
all minimal charts $\Gamma$
satisfy the following four conditions 
(see \cite{ChartApp1},\cite{ChartAppII},\cite{ChartAppIII},
\cite{StI}):

%%%%%%%%%%%%%%%%%%%%%%%%%%%%%%%%%%%%%%%%%
\begin{assumption}
\label{AssumeTerminal}
If an edge of $\Gamma$
contains a black vertex,
then the edge is a free edge 
or a terminal edge.
Moreover 
any terminal edge contains a middle arc.
\end{assumption}

%%%%%%%%%%%%%%%%%%%%%%%%%%%%%%%%%%%%%%%
\begin{assumption}
\label{NoSimpleHoop}
All free edges and simple hoops in $\Gamma$ 
are moved into a small neighborhood $U_\infty$ 
of the point at infinity $\infty$. 
Hence
we assume that 
$\Gamma$ does not contain free edges
nor simple hoops, 
otherwise mentioned. 
\end{assumption}

%%%%%%%%%%%%%%%%%%%%%%%%%%%%%%%%%%%%%%%%	
\begin{assumption}
\label{Ring}
Each complementary domain of
any ring and hoop must contain 
at least one white vertex. 
\end{assumption}

\begin{assumption}
\label{Infinity}
The point at infinity $\infty$ is moved in any complementary domain of $\Gamma$.
\end{assumption}

In this paper
for a subset $X$ in a space
we denote 
the interior of $X$,
the boundary of $X$ and
the closure of $X$
by Int$X$, $\partial X$
and $Cl(X)$
respectively.

Let $\alpha$ be a simple arc or an edge,
and $p,q$ the endpoints of $\alpha$.
We denote 
$\partial \alpha=\{p,q\}$
and ${\rm Int}\alpha=\alpha-\{p,q\}$.

%%%%%%%%%%%%%%%%%%%%%%

%\newpage

\section{A disk with exactly two white vertices}
\label{s:DiskTwo}

In this section
for a minimal chart $\Gamma$,
we investigate
a disk $D$ with exactly two white vertices 
of $\Gamma_m$
such that
$\Gamma_m\cap\partial D$
consists of at most two points.

Let $\Gamma$ be a chart,
and $m$ a label of $\Gamma$. 
A {\it loop} is a simple closed curve in $\Gamma_m$ with exactly one white vertex
(possibly with crossings).

Let $\Gamma$ be a chart, and $\ell$ a loop of label $m$.
Let $e$ be the edge of label $m$ 
such that $e\cap\ell$ is a white vertex.
Then the loop $\ell$ bounds two disks on the 2-sphere.
One of the two disks does not contain the edge $e$.
The disk is called the {\it associated disk} of the loop $\ell$.

Let $X$ be a set in a chart $\Gamma$.
Let
 $$w(X)=\text{the number of white vertices in $X$.}$$

\begin{lemma}
{\rm (\cite[Lemma 4.2]{ChartAppII})}
\label{LemmaInsideLoopOutsideLoop}
Let $\Gamma$ be a minimal chart, and
$m$ a label of $\Gamma$.
Let $\ell$ be a loop of label $m$, and
$D$ the associated disk of the loop $\ell$.
Then $w(\Gamma\cap {\rm Int}D)\ge2$
and $w(\Gamma\cap (S^2-D))\ge2$.
\end{lemma}

\begin{lemma}$(${\rm \cite[Theorem 1.1]{ChartAppIV}}$)$
\label{LemmaNoLoop}
There is no loop in any minimal chart with exactly seven white vertices.
\end{lemma}

Let $\Gamma$ be a minimal chart.
A disk $D$ is said to be {\it admissible} 
provided that 
the boundary $\partial D$ does not contain any vertex of 
$\Gamma$, and
any edge $e$ with $e\cap\partial D\not=\emptyset$
intersects $\partial D$ transversely.

By Lemma~\ref{LemmaInsideLoopOutsideLoop},
we have the following lemma:

\begin{lemma}
\label{LoopInAdmissibleDisk}
Let $\Gamma$ be a minimal chart.
Let $D$ be an admissible disk.
If $w(\Gamma\cap D)\le2$,
then $D$ does not contain any loop.
\end{lemma}

Let $\Gamma$ be a chart, 
and $m$ a label of $\Gamma$.
Let $L$ be the closure of a connected component 
of the set obtained by taking out 
all the white vertices from $\Gamma_m$.
If $L$ contains at least one white vertex
but does not contain any black vertex,
then $L$ is called an {\it internal edge of label $m$}.
Note that an internal edge may contain a crossing of $\Gamma$.

\begin{lemma}
\label{LemmaAtMostOnePoint}
Let $\Gamma$ be a minimal chart, 
and
$m$ a label of $\Gamma$.
Let $D$ be an admissible disk.
Suppose that
$\Gamma_m\cap \partial D$ is at most one point.
If $w(\Gamma_{m}\cap D)\ge1$,
then $w(\Gamma\cap D)\ge2$.
\end{lemma}

\begin{Proof}
Suppose that $w(\Gamma\cap D)=1$.
Then $w(\Gamma_m\cap D)=1$.
Let $G$ be a connected component of $\Gamma_m\cap D$
with $w(G)=1$.
Since $\Gamma_m\cap \partial D$ is at most one point,
the set $G\cap \partial D$ is at most one point.

Let $w$ be the white vertex in $G$,
and $e_1,e_2$ internal edges (possibly terminal edges)
of label $m$ at $w$
not middle at $w$.
Since $G\cap \partial D$ is at most one point,
one of $e_1,e_2$ does not intersect $\partial D$,
say $e$.
Since $D$ does not contain any loop 
by Lemma~\ref{LoopInAdmissibleDisk},
the edge $e$ is not a loop.
Thus the condition $w(G)=1$ implies that 
the edge $e$ is a terminal edge.
However the edge $e$ is not middle at $w$.
This contradicts Assumption~\ref{AssumeTerminal}.
Hence we have $w(\Gamma\cap D)\ge2$.
\end{Proof}

Let $\Gamma$ be a chart. 
Suppose that an object consists of 
some edges of $\Gamma$, arcs in edges of 
$\Gamma$ and arcs around white vertices.
Then the object is called a {\it pseudo chart}.

%%%%%%%%%%%%%

In \cite[Lemma~5.3(1)]{ChartAppV},
the chart $\Gamma$ has the condition without loops.
However by Lemma~\ref{LemmaInsideLoopOutsideLoop},
we can show that we do not need this condition. 
Thus we have the following lemma.

\begin{figure}[htb]
\centerline{\includegraphics{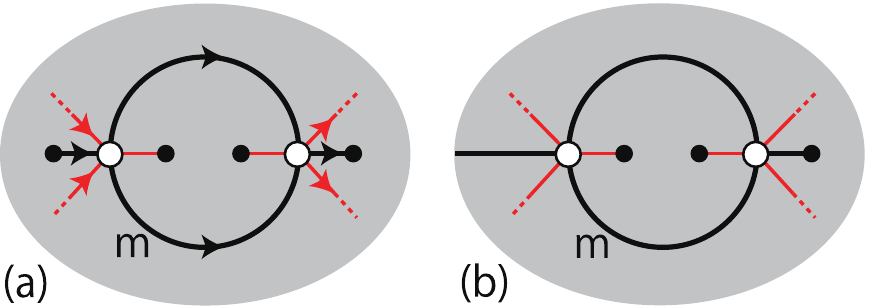}}
\caption{\label{fig03}
The gray region is the disk $D$,
and $m$ is a label.}
\end{figure}

\begin{lemma}
{\rm (cf. \cite[Lemma 5.3(1)]{ChartAppV})}
\label{LemmaTwoWhiteVertices}
Let $\Gamma$ be a minimal chart, 
and
$m$ a label of $\Gamma$.
Let $D$ be an admissible disk.
Suppose that
$\Gamma_m\cap \partial D$ is at most one point.
If $w(\Gamma\cap D)=w(\Gamma_{m}\cap D)=2$,
then the disk $D$ contains one of the 
two pseudo charts
as shown in 
Fig.~\ref{fig03}.
\end{lemma}

%%%%%%%%%%%%%
%%%%%%%%%%%%%%
%%%%%%%%%%%%

Let $\Gamma$ be a chart, $m$ a label of $\Gamma$, 
$D$ a disk with $\partial D\subset \Gamma_m$, 
and $k$ a positive integer.
If $\partial D$ contains exactly
$k$ white vertices, 
then $D$ is called 
{\it a $k$-angled disk of $\Gamma_m$}. 
Note that 
the boundary $\partial D$ may contain crossings.

Let $\Gamma$ be a chart, and
$m$ a label of $\Gamma$.
An edge of label $m$ is called a {\it feeler} of a $k$-angled disk $D$ of $\Gamma_m$
if the edge intersects $N-\partial D$
where $N$ is a regular neighborhood of $\partial D$ in $D$.

\begin{figure}
\centerline{\includegraphics{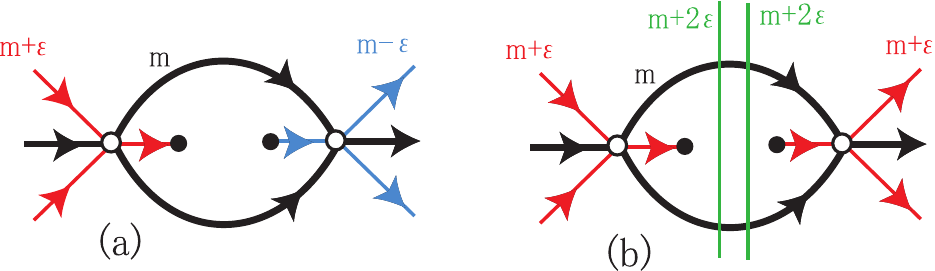}}
\caption{\label{fig04}
$m$ is a label,
and $\varepsilon\in\{+1,-1\}$.}
\end{figure}

\begin{lemma}
\label{Theorem2AngledDisk}
{\rm (\cite[Corollary 5.8]{ChartAppII})}
Let $\Gamma$ be a minimal chart.
Let $D$ be a $2$-angled disk of $\Gamma_m$ with at most one feeler.
If $w(\Gamma\cap{\rm Int}D)=0$,
then a regular neighborhood of $D$ contains one of two pseudo charts as shown in Fig.~\ref{fig04}.
\end{lemma}

%%%%%%%%%%%%%%%%%%
%%%%%%%%%%%%%%%%%% Figure
%%%%%%%%%%%%%%%%%%

\begin{figure}[htb]
\centerline{\includegraphics{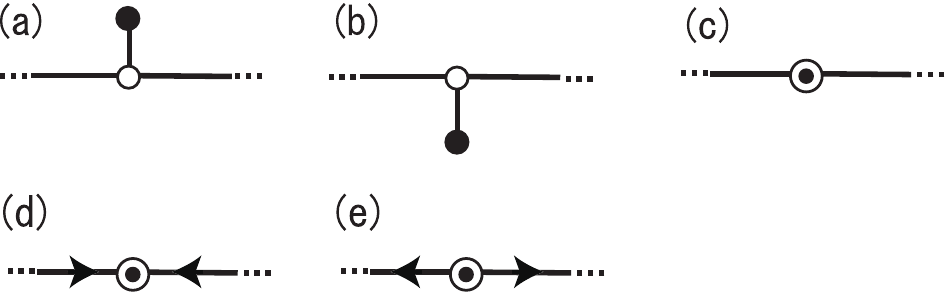}}
\caption{\label{fig05}
(a),(b) White vertices in terminal edges.
(c),(d),(e) BW-vertices.}
\end{figure}

 In our argument  we often construct a chart $\Gamma$. 
On the construction of a chart $\Gamma$, for a white vertex $w\in\Gamma_m$ for some label $m$,  
among the three edges of $\Gamma_m$ 
containing $w$, 
if one of the three edges is a terminal edge 
(see Fig.~\ref{fig05}(a) and (b)), 
then we remove the terminal edge and
put a black dot at the center of the white vertex  as shown in Fig.~\ref{fig05}(c).
Namely
Fig.~\ref{fig05}(c) means 
Fig.~\ref{fig05}(a) or 
Fig.~\ref{fig05}(b).
We call the vertex in Fig.~\ref{fig05}(c) 
a {\it BW-vertex} with respect to $\Gamma_m$.

%%%%%%%%%%%%%%%%%%
%%%%%%%%%%%%%%%%%% Figure
%%%%%%%%%%%%%%%%%%

\begin{lemma}
\label{OriBWvertex}
{\rm (\cite[Lemma 3.1]{ChartAppV})}
In a minimal chart $\Gamma$,
for each BW-vertex in $\Gamma_m$,
the two edges of label $m$ containing the BW-vertex
are oriented inward or outward at the BW-vertex
simultaneously
if each of the two edges is not a terminal edge
$($see Fig.~\ref{fig05}$($d$)$ and $($e$))$.
\end{lemma}

\begin{lemma}
\label{LemmaWithTerminal3}
{\rm (\cite[Lemma 3.2(1)]{ChartAppV})}
Let $\Gamma$ be a minimal chart,
and $m$ a label of $\Gamma$.
Let $G$ be a connected component of $\Gamma_m$.
If $1\le w(G)$, then $2\le w(G)$.
\end{lemma}

Let $D$ be a compact surface.
A simple arc $\alpha$ in $D$
is a {\it proper arc} of $D$
if $\alpha\cap\partial D=\partial \alpha$.

\begin{figure}[htb]
\centerline{\includegraphics{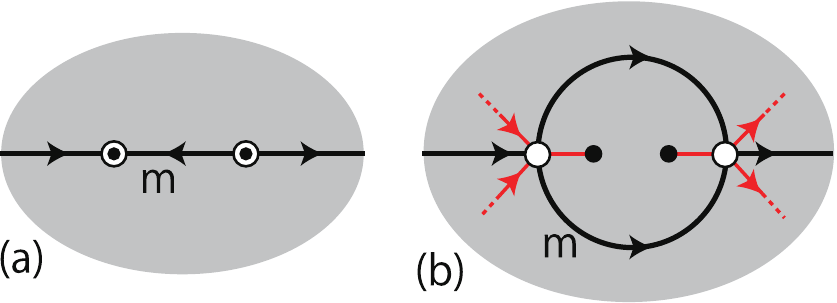}}
\caption{\label{fig06}
The gray region is the disk $D$,
and $m$ is a label.}
\end{figure}

\begin{lemma}
\label{LemmaTwoWhiteVerticesWithBoundaryTwoPoints}
Let $\Gamma$ be a minimal chart, 
and
$m$ a label of $\Gamma$.
Let $D$ be an admissible disk.
Suppose that
$\Gamma_m\cap \partial D$ is exactly two points.
If $w(\Gamma\cap D)=w(\Gamma_{m}\cap D)=2$,
then the disk $D$ contains one of the 
four pseudo charts
as shown in 
Fig.~\ref{fig03} and
Fig.~\ref{fig06}.
\end{lemma}

%%%%%%%%%%%%%
%%%%%%%%%%%%%%
%%%%%%%%%%%%

\begin{Proof}
Let $G$ be a connected component of $\Gamma_m \cap D$
with $w(G)\ge1$.
Since $\Gamma_m\cap \partial D$ consists of exactly two points,
there are three cases:
(i) $G\cap\partial D=\emptyset$,
(ii) $G\cap\partial D$ consists of exactly one point,
(iii) $G\cap\partial D$ consists of exactly two points.

{\bf Case (i).}
By Lemma~\ref{LemmaWithTerminal3},
we have $w(G)=2$. 
Let $e_1,e_2$ be internal edges (possibly terminal edges) 
of label $m$
with $e_1\cap \partial D\not=\emptyset$ for $i=1,2$.
Then $(e_1\cup e_2)\cap G=\emptyset$.
Hence neither $e_1$ nor $e_2$ contains a white vertex in $D$.
Let $N$ be a regular neighborhood of $(e_1\cup e_2)\cap D$ 
in $D$.
Then a connected component of $Cl(D-N)$ contains 
the connected component $G$,
say $D'$.
Moreover, the set $D'$ is an admissible disk 
with $\Gamma_m\cap \partial D'=\emptyset$, and
$w(\Gamma\cap D')=w(\Gamma_{m}\cap D')=2$.
Thus by Lemma~\ref{LemmaTwoWhiteVertices},
 the disk $D'$ contains the 
the pseudo chart
as shown in 
Fig.~\ref{fig03}(a),
so does $D$.

{\bf Case (ii).}
In a similar way of the proof of Lemma~\ref{LemmaAtMostOnePoint}, we can show $w(G)=2$.
Moreover by the similar way of the proof of Case (i),
we can show that 
 the disk $D$ contains the 
the  pseudo chart
as shown in 
Fig.~\ref{fig03}(b).

{\bf Case (iii).}
Let $e_1,e_2$ be internal edges (possibly terminal edges)
of label $m$ with $e_i\cap\partial D\not=\emptyset$
for $i=1,2$.
Since $G\cap\partial D$ consists of exactly two points,
we have $e_i\cap D \subset G$ for $i=1,2$.
Since $w(G)\ge1$,
each of $e_1,e_2$ contains a white vertex in $D$.
Let $w_1,w_2$ be white vertices in $D$
with $w_1\in e_1$ and $w_2\in e_2$.

We shall show $w_1\not=w_2$.
If $w_1=w_2$,
then the arc $(e_1\cup e_2)\cap D$ is a proper arc of $D$.
Let $N$ be a regular neighborhood of $(e_1\cup e_2)\cap D$ in
$D$.
Then $Cl(D-N)$ consists of two disks.
One of the two disks contains one white vertex,
say $E$.
Moreover $\Gamma_m\cap \partial E$ is at most one point.
Thus by Lemma~\ref{LemmaAtMostOnePoint},
we have $w(\Gamma\cap E)\ge2$.
Hence $w(\Gamma\cap D)\ge3$.
This contradicts $w(\Gamma\cap D)=2$.
Thus $w_1\not=w_2$.

Let $e_3$ be an internal edge (possibly terminal edge)
of label $m$ at $w_1$ 
different from $e_1$ but not middle at $w_1$.
By Assumption~\ref{AssumeTerminal},
the edge $e_3$ is not a terminal edge.
Since $D$ does not contain any loop 
by Lemma~\ref{LoopInAdmissibleDisk},
the edge $e_3$ is not a loop,
i.e. $e_3\ni w_2$.

If both of $w_1$ and $w_2$ are BW-vertices with respect to $\Gamma_m$,
then by Lemma~\ref{OriBWvertex}
the disk $D$ contains the pseudo chart
as shown in 
Fig.~\ref{fig06}(a).

If one of $w_1,w_2$ is not a BW-vertex 
with respect to $\Gamma_m$,
then there exists an internal edge $e_3'$
of label $m$ connecting $w_1$ and $w_2$ 
different from $e_3$.
Thus there exists a 2-angled disk of $\Gamma_m$
without feelers in $D$.
Hence by Lemma~\ref{Theorem2AngledDisk}
the disk $D$ contains the pseudo chart
as shown in 
Fig.~\ref{fig06}(b).
\end{Proof}

%%%%%%%%%%%%%
%%%%%%%%%%%%%

%\newpage
\section{3-angled disks with at most one white vertex
in their interiors}
\label{s:3AngledDiskOneWhiteVertex}

In this section we investigate 3-angled disks with
at most one white vertex in their interiors.

Let $\Gamma$ be a chart,
and $D$ a $k$-angled disk of $\Gamma_m$.
If any feeler of $D$ of label $m$ is a terminal edge,
then $D$ is called a {\it special} $k$-angled disk.

Let $\Gamma$ be a chart, 
$D$ a disk, and 
$G$ a pseudo chart with $G \subset D$.
Let $r:D\to D$ be a reflection of $D$, and $G^*$ the pseudo chart obtained from $G$ by changing the orientations of all of the edges.
Then the set $\{G,G^*, r(G), r(G^*)\}$ 
is called the {\it RO-family of the pseudo chart $G$}.

\begin{figure}[htb]
\centerline{\includegraphics{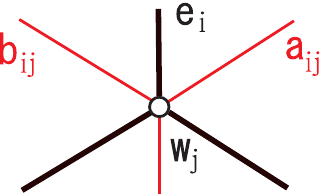}}
\caption{\label{fig07}
The three edges $a_{ij},e_i,b_{ij}$ are consecutive edges around the white vertex $w_j$.}
\end{figure}

In our argument,
we often need a name for an unnamed edge by using a given edge and a given white vertex.
For the convenience,
we use the following naming:
Let $e',e_i,e''$ be three consecutive edges containing  a white vertex $w_j$. Here, 
the two edges $e'$ and $e''$ are unnamed edges. 
There are six arcs in a neighborhood $U$ of the white vertex $w_j$. 
If the three arcs $e'\cap U$, $e_i \cap U$, $e'' \cap U$ lie anticlockwise around the white vertex $w_j$ in this order, 
then $e'$ and $e''$ are denoted by $a_{ij}$ and $b_{ij}$ 
respectively (see Fig.~\ref{fig07}).
There is a possibility $a_{ij}=b_{ij}$ if they are contained in a loop.

%%%%%%%%%%%%%%%%%%
%%%%%%%%%%%%%%%%%% Figure
%%%%%%%%%%%%%%%%%%

\begin{figure}[htb]
\centerline{\includegraphics{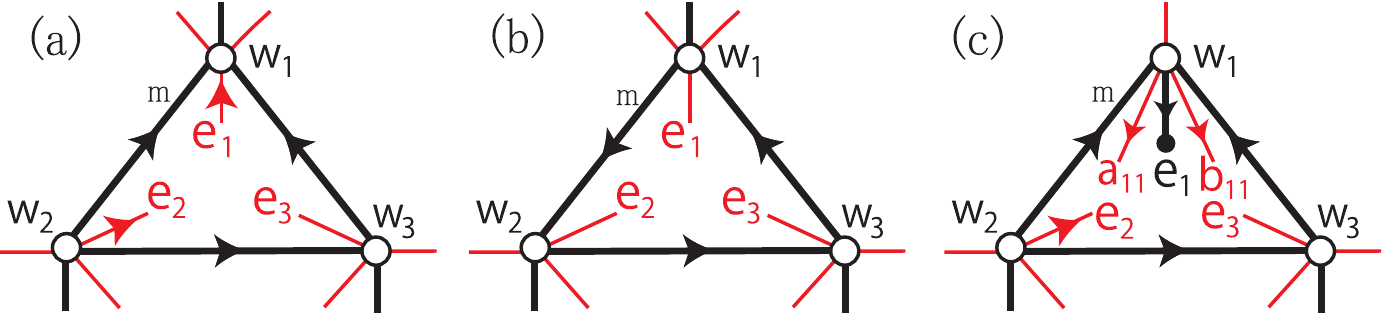}}
\caption{\label{fig08} 
The 3-angled disks (a) and (b) have no feelers,
 the other has one feeler, and
$m$ is a label.}
\end{figure}

We can show the following lemma easily:

\begin{lemma}
\label{ROfamily3AngledDisk}
{\rm (cf. \cite[Lemma 4.2]{ChartAppIII})}
Let $\Gamma$ be a minimal chart, and $m$ a label of $\Gamma$.
Let $D$ be a special $3$-angled disk of $\Gamma_m$
with at most one feeler.
Then a regular neighborhood of $D$ contains one of 
the RO-families of 
the three pseudo charts as shown in 
Fig.~\ref{fig08}.
\end{lemma}

%%%%%%%%%%%%%%%%%%
%%%%%%%%%%%%%%%%%% Figure
%%%%%%%%%%%%%%%%%%

We shall show the following lemma in this section.

\begin{figure}[htb]
\centerline{\includegraphics{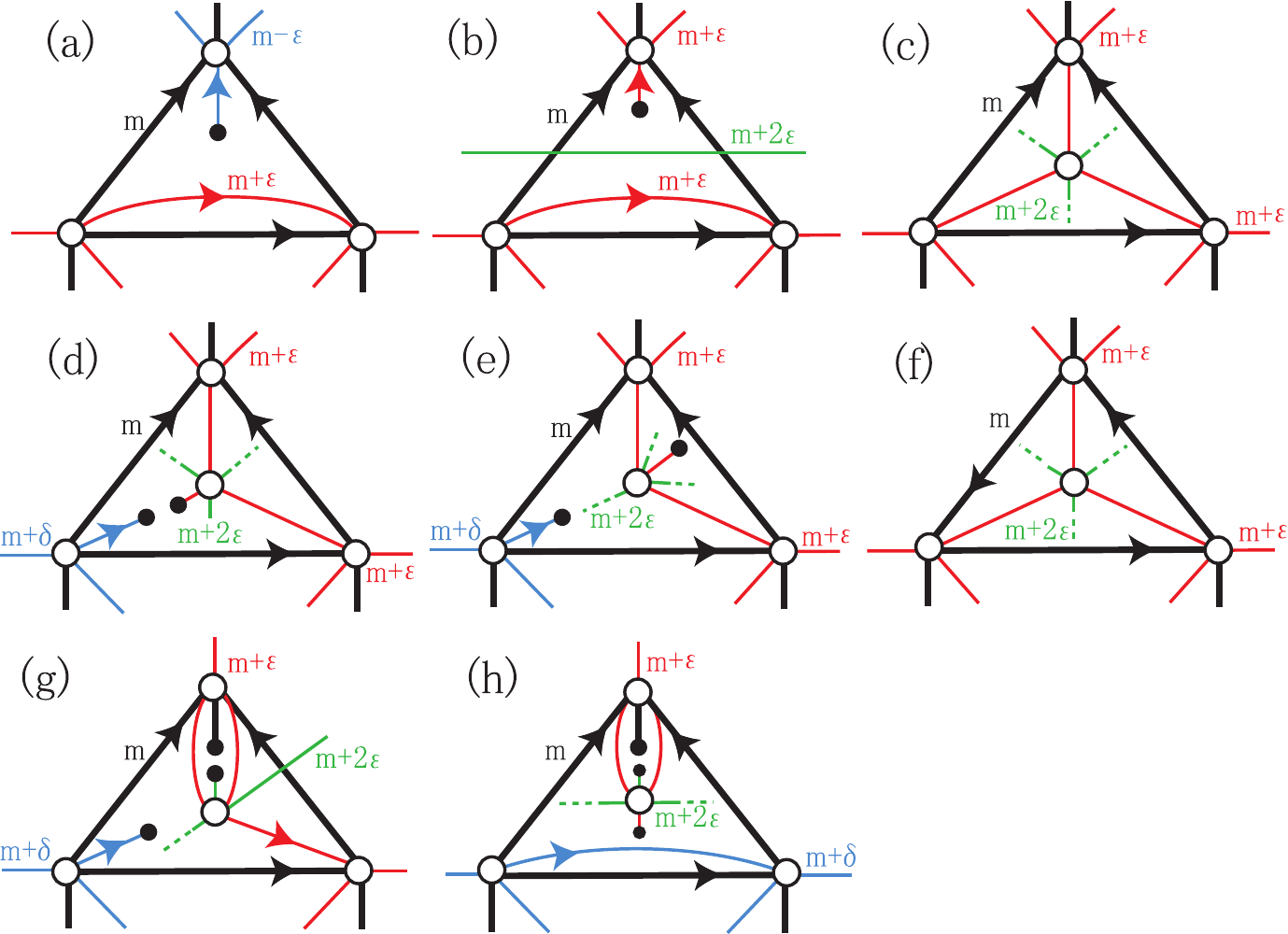}}
\caption{\label{fig09}
The 3-angled disks (g) and (h) have one feeler,
the others have no feelers,
$m$ is a label, $\varepsilon,\delta\in\{+1,-1\}$.}
\end{figure}

\begin{lemma}
{\rm (cf. \cite[Theorem 1.2]{ChartAppIII})}
\label{Theorem3AngledDisk}
Let $\Gamma$ be a minimal chart, and $m$ a label of $\Gamma$.
Let $D$ be a special $3$-angled disk of $\Gamma_m$
with at most two feelers.
Then we have the following:
\begin{enumerate}
\item[{\rm (a)}]
If $w(\Gamma\cap {\rm Int}D)=0$,
then a regular neighborhood of $D$ contains one of the RO-families of the two pseudo charts as shown in 
Fig.~\ref{fig09}$($a$)$ and $($b$)$.
\item[{\rm (b)}]
If $w(\Gamma\cap {\rm Int}D)=w(\Gamma_{m+\varepsilon}\cap {\rm Int}D)=1$ for some $\varepsilon\in\{+1,-1\}$,
then a regular neighborhood of $D$ contains one of the RO-families of the six pseudo charts as shown in 
Fig.~\ref{fig09}$($c$)$,$($d$)$,$($e$)$,$($f$)$,$($g$)$ 
and $($h$)$.
\end{enumerate}
\end{lemma}

%%%%%%%%%%%%%%%%%%
%%%%%%%%%%%%%%%%%% Figure
%%%%%%%%%%%%%%%%%%

To prove the above lemma,
we need five lemmata.

\begin{lemma}
\label{ROfamily3AngledDiskTypeA}
{\rm (\cite[Lemma 4.3]{ChartAppIII})}
Let $\Gamma$ be a minimal chart, and $m$ a label of $\Gamma$.
Let $D$ be a $3$-angled disk of $\Gamma_m$
as shown in Fig.~\ref{fig08}$($a$)$.
If $w(\Gamma\cap{\rm Int}D)=0$,
then a regular neighborhood of $D$ 
contains one of the RO-families of 
the two pseudo charts as shown in 
Fig.~\ref{fig09}$($a$)$ and $($b$)$.
\end{lemma}

\begin{lemma}
\label{ROfamily3AngledDiskTypeB}
{\rm (\cite[Lemma 4.5]{ChartAppIII})}
Let $\Gamma$ be a minimal chart, and $m$ a label of $\Gamma$.
Let $D$ be a $3$-angled disk of $\Gamma_m$
as shown in Fig.~\ref{fig08}$($b$)$.
Then $w(\Gamma\cap{\rm Int}D)\ge1$.
Moreover, if $w(\Gamma\cap{\rm Int}D)=1$,
then a regular neighborhood of $D$ 
contains one of the RO-family of 
the pseudo chart as shown in 
Fig.~\ref{fig09}$($f$)$.
\end{lemma}

\begin{lemma}
\label{ROfamily3AngledDiskTypeC}
{\rm (\cite[Lemma 5.2]{ChartAppIII})}
Let $\Gamma$ be a minimal chart, and $m$ a label of $\Gamma$.
Let $D$ be a $3$-angled disk of $\Gamma_m$
as shown in Fig.~\ref{fig08}$($c$)$.
Then $w(\Gamma\cap{\rm Int}D)\ge1$.
Moreover if $w(\Gamma\cap{\rm Int}D)=1$,
then a regular neighborhood of $D$ 
contains one of the RO-families of 
the two pseudo charts as shown in 
Fig.~\ref{fig09}$($g$)$ and $($h$)$.
\end{lemma}

\begin{lemma}
\label{ROfamily3AngledDiskTwoFeelers}
{\rm (\cite[Lemma 5.3]{ChartAppIII})}
Let $\Gamma$ be a minimal chart, and $m$ a label of $\Gamma$.
Let $D$ be a special $3$-angled disk of $\Gamma_m$
with two feelers.
Then $w(\Gamma\cap{\rm Int}D)\ge2$.
\end{lemma}

By Lemma~\ref{LemmaInsideLoopOutsideLoop},
we have the following lemma:

\begin{lemma}
\label{LoopInKAngledDisk}
Let $\Gamma$ be a minimal chart,
and $m$ a label of $m$.
Let $D$ be a $k$-angled disk of $\Gamma_m$.
If $w(\Gamma\cap {\rm Int}D)\le 2$,
then ${\rm Int}D$ does not contain any loop.
\end{lemma}

{\it Proof of Lemma~\ref{Theorem3AngledDisk}.}
By Lemma~\ref{ROfamily3AngledDiskTwoFeelers},
we can assume that the disk $D$ has at most one feeler.
Moreover by Lemma~\ref{ROfamily3AngledDisk},
a regular neighborhood of $D$ contains one of the
RO-families of the three pseudo charts as shown in 
Fig.~\ref{fig08}.

{\bf Statement (a).}
Suppose that $w(\Gamma\cap{\rm Int}D)=0$.
By Lemma~\ref{ROfamily3AngledDiskTypeB} and Lemma~\ref{ROfamily3AngledDiskTypeC},
the disk $D$ is a 3-angled disk 
as shown in Fig.~\ref{fig08}(a).
Thus by Lemma~\ref{ROfamily3AngledDiskTypeA},
a regular neighborhood of $D$ contains one of RO-families 
of the two pseudo charts as shown 
in Fig.~\ref{fig09}(a) and (b).

{\bf Statement (b).}
Suppose that $w(\Gamma\cap {\rm Int}D)=w(\Gamma_{m+\varepsilon}\cap {\rm Int}D)=1$.
If $D$ is a 3-angled disk as shown in Fig.~\ref{fig08}(b) or (c),
then by Lemma~\ref{ROfamily3AngledDiskTypeB} and
Lemma~\ref{ROfamily3AngledDiskTypeC}
a regular neighborhood of $D$ contains 
one of the RO-families of the three pseudo charts 
as shown in Fig.~\ref{fig09}$($f$)$,$($g$)$ and $($h$)$.

Now suppose that $D$ is a 3-angled disk as shown in Fig.~\ref{fig08}(a).
We use the notations as shown in 
Fig.~\ref{fig08}(a), where
\begin{enumerate}
\item[(1)] $e_3$ is an edge at $w_3$ in $D$ but 
not middle at $w_3$
(i.e. the edge $e_3$ is not a terminal edge 
by Assumption~\ref{AssumeTerminal}).
\end{enumerate}

Let $w$ be the white vertex in ${\rm Int}D$.
Since $w(\Gamma\cap {\rm Int}D)=w(\Gamma_{m+\varepsilon}\cap {\rm Int}D)=1$,
we have $w\in\Gamma_m\cap\Gamma_{m+\varepsilon}$
or
$w\in\Gamma_{m+\varepsilon}\cap\Gamma_{m+2\varepsilon}$.
If  $w\in\Gamma_m\cap\Gamma_{m+\varepsilon}$,
then there exists
a connected component $G$ of $\Gamma_m$ with $w(G)=1$.
This contradicts Lemma~\ref{LemmaWithTerminal3}.
Thus $w\in\Gamma_{m+\varepsilon}\cap\Gamma_{m+2\varepsilon}$.

Let $e,e',e''$ be internal edges (possibly terminal edges)
of label $m+\varepsilon$ at $w$.
By Lemma~\ref{LoopInKAngledDisk},
none of  $e,e',e''$ are loops.
By Assumption~\ref{AssumeTerminal},
\begin{enumerate}
\item[(2)]
at most one of  $e,e',e''$ is a terminal edge.
\end{enumerate}

If none of $e,e',e''$ are terminal edges
(i.e. each of $e,e',e''$ contains one of white vertices
$w_1,w_2,w_3$),
then a regular neighborhood of $D$ 
contains one of the RO-family 
of the pseudo chart as shown in Fig.~\ref{fig09}$($c$)$.

If one of  $e,e',e''$ is a terminal edge,
then by (2)
two of  $e,e',e''$ are not terminal edges.
Since the edge $e_3$ is not a terminal edge by (1),
one of $e,e',e''$ contains $w_3$ and one
of  $e,e',e''$ contains $w_1$ or $w_2$.
Hence a regular neighborhood of $D$ 
contains one of the RO-families 
of the two pseudo charts as shown 
in Fig.~\ref{fig09}$($d$)$ and $($e$)$.
{\hfill {$\square$}\vspace{1.5em}}

%%%%%%%%%%%%%
%%%%%%%%%%%%%

%%%%%%%%%%%%%
%%%%%%%%%%%%%

%\newpage
\section{2-angled disks without feelers}
\label{s:2-angledDisks}

\begin{figure}
\centerline{\includegraphics{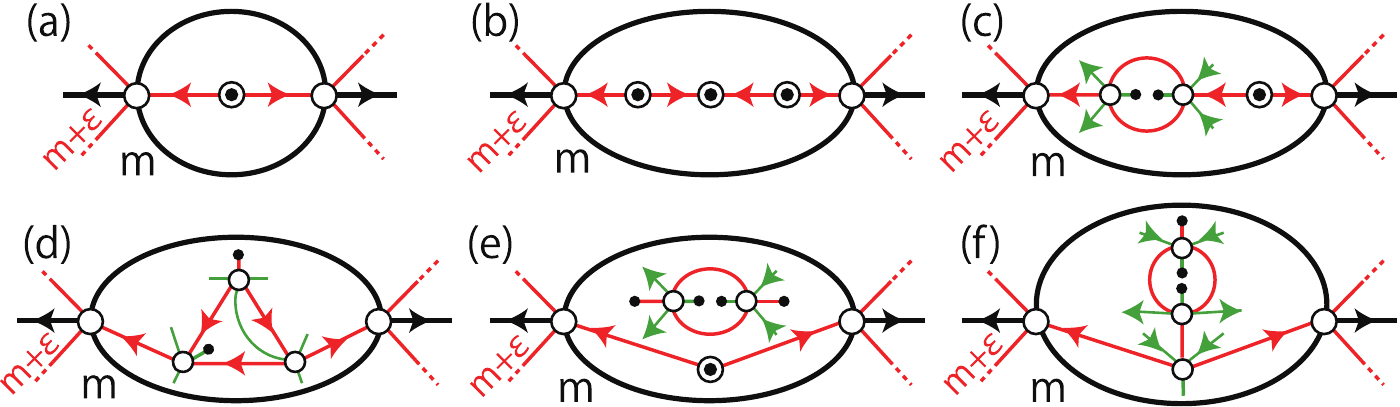}}
\caption{\label{fig10}
The 2-angled disk (a) has one white vertex in its interior.
The other 2-angeld disks have three white vertices in their interiors.}
\end{figure}

In this section we investigate a 2-angled disk without feelers.
In this section we shall show the following lemma:

\begin{lemma}
\label{Lemma2AngledDisks}
Let $\Gamma$ be a minimal chart,
and $m$ a label of $\Gamma$.
Let $D$ be a $2$-angled disk of $\Gamma_m$
without feelers, and
$w_1,w_2$ the white vertices in $\partial D$.
Let $e_1,e_2$ be the internal edges
$($possibly terminal edges$)$ of label $m$
at $w_1,w_2$, respectively,
such that $e_1\not\subset D$
and $e_2\not\subset D$.
Suppose that the two edges $e_1,e_2$ are oriented outward
$($resp. inward$)$ at $w_1,w_2$, 
respectively.
Then we have the following:
\begin{enumerate}
\item[{\rm (a)}]
$w(\Gamma\cap{\rm Int}D)\ge1$.
\item[{\rm (b)}]
If $w(\Gamma\cap {\rm Int}D)=1$,
then a regular neighborhood of $D$ contains one of the RO-family of the pseudo chart as shown in Fig.~\ref{fig10}$($a$)$.
\item[{\rm (c)}]
If $w_1,w_2\in\Gamma_{m+\varepsilon}$ for some $\varepsilon\in\{+1,-1\}$ and
if $w(\Gamma_{m+\varepsilon}\cap {\rm Int}D)=2$,
then $w(\Gamma\cap {\rm Int}D)\ge3$.
\item[{\rm (d)}]
If $w_1,w_2\in\Gamma_{m+\varepsilon}$ for some $\varepsilon\in\{+1,-1\}$ and
if $w(\Gamma\cap {\rm Int}D)=w(\Gamma_{m+\varepsilon}\cap {\rm Int}D)=3$,
then a regular neighborhood of $D$ contains 
one of RO-families of the five pseudo charts as shown in Fig.~\ref{fig10}$($b$)$,$($c$)$,$($d$)$,$($e$)$,$($f$)$.
\end{enumerate}
\end{lemma}

%%%%%%%%%%%%%%%%%%
%%%%%%%%%%%%%%%%%% Figure
%%%%%%%%%%%%%%%%%%

\begin{remark}
\label{Remark2angledDisk}
\upshape
Lemma~\ref{Lemma2AngledDisks}(a) is the same as
 \cite[Lemma 3.6(a)]{ChartAppVIII}.
\end{remark}

\begin{figure}
\centerline{\includegraphics{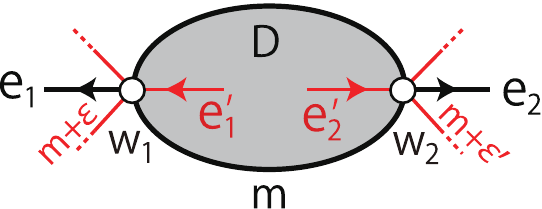}}
\caption{\label{fig11}
The gray region is a 2-angeld disk $D$ without feelers,
$m$ is a label, $\varepsilon,\varepsilon'\in\{+1,-1\}$.}
\end{figure}

From now on throughout this section,

\begin{enumerate}
\item[(i)]
we use the notations in Lemma~\ref{Lemma2AngledDisks}
as shown in Fig.~\ref{fig11},
\item[(ii)]
we may assume that
the two edges $e_1,e_2$ are oriented outward at $w_1,w_2$, 
respectively.
\end{enumerate}

%%%%%%%%%%%%%%%%%%
%%%%%%%%%%%%%%%%%% Figure
%%%%%%%%%%%%%%%%%%

\begin{lemma}
\label{E1'E2'NotTerminal}
Let $\Gamma,D,w_1,w_2$ be as above.
Let $e_1',e_2'$ be internal edges $($possibly terminal edges$)$ of label $m+\varepsilon, m+\varepsilon'$
at $w_1,w_2$, respectively,
with $e_1'\cup e_2'\subset D$
for some $\varepsilon,\varepsilon'\in\{+1,-1\}$.
Then neither $e_1'$ nor $e_2'$ is a terminal edge.
\end{lemma}

\begin{Proof}
We shall show that neither $e_1'$ nor $e_2'$ is middle at $w_1$ or $w_2$. 
If $e_1'$ is middle at $w_1$,
then by Condition (ii) in this section
the two internal edges in $\partial D$
are oriented from $w_2$ to $w_1$.
Hence there are three edges of label $m$ 
oriented outward at $w_2$.
This contradicts the definition of the chart.
Thus $e_1'$ is not middle at $w_1$.

Similarly we can show that $e_2'$ is not middle at $w_2$.
Hence by Assumption~\ref{AssumeTerminal},
neither $e_1'$ nor $e_2'$ is a terminal edge.
\end{Proof}

{\it Proof of Lemma~\ref{Lemma2AngledDisks}$($b$)$.}
Let $w$ be the white vertex in ${\rm Int}D$.
By Lemma~\ref{E1'E2'NotTerminal},
neither $e_1'$ nor $e_2'$ is a terminal edge.
Since $e_1'\not=e_2'$ by Condition (ii) of this section,
we have $e_1'\cap e_2'\ni w$. Thus
$w\in \Gamma_{m+\varepsilon}$ for some $\varepsilon\in\{+1,-1\}$, and
there exists a terminal edge of label $m+\varepsilon$ at $w$.
Hence a regular neighborhood of $D$ contains the pseudo chart as shown in Fig.~\ref{fig10}(a).
We complete the proof of Lemma~\ref{Lemma2AngledDisks}$($b$)$. {\hfill {$\square$}\vspace{1.5em}}

%%%%%%%%%%%%%%%%%%%%%%%%%%%%
%%%%%%%%%%%%%%%%%%%%%%%%%%%%

{\it Proof of Lemma~\ref{Lemma2AngledDisks}$($c$)$.}
Suppose $w(\Gamma\cap {\rm Int}D)=2$.
Let $N$ be a regular neighborhood of $\partial D$ in $D$.
Set $D'=Cl(D-N)$.
Then $D'$ is an admissible disk
such that $\Gamma_{m+\varepsilon}\cap\partial D'$ is 
exactly two points.
Moreover we have $w(\Gamma\cap D')=w(\Gamma_{m+\varepsilon}\cap D')=2$.
Hence 
by Lemma~\ref{LemmaTwoWhiteVerticesWithBoundaryTwoPoints},
the disk $D'$ contains one of the 
four pseudo charts
as shown in 
Fig.~\ref{fig03} and
Fig.~\ref{fig06}.

By Lemma~\ref{E1'E2'NotTerminal},
neither $e_1'$ nor $e_2'$ is a terminal edge.
Thus the disk $D'$ contains one of the 
two pseudo charts
as shown in 
Fig.~\ref{fig06}.
Hence
either $e_1'$ is oriented outward at $w_1$,
or $e_2'$  is oriented outward at $w_2$.
Thus 
either $e_1$ is oriented inward at $w_1$,
or $e_2$  is oriented inward at $w_2$.
However, this contradicts Condition (ii) of this section.
Therefore if $w(\Gamma_{m+\varepsilon}\cap {\rm Int}D)=2$,
then $w(\Gamma\cap {\rm Int}D)\ge3$.
We complete the proof of Lemma~\ref{Lemma2AngledDisks}$($c$)$.
{\hfill {$\square$}\vspace{1.5em}}

%%%%%%%%%%%%%%%%%%%%
%%%%%%%%%%%%%%%%%%%%
%%%%%%%%%%%%%%%%%%%%

{\it Proof of Lemma~\ref{Lemma2AngledDisks}$($d$)$.}
For each $i=1,2$,
let $G_i$ be the connected component of $\Gamma_{m+\varepsilon}\cap D$ with $G_i\supset e_i'$.
By Lemma~\ref{E1'E2'NotTerminal},
the edge $e_i'$ is not a terminal edge,
i.e. the edge $e_i'$ contains a white vertex in ${\rm Int}D$.
Thus $w(G_1\cap {\rm Int}D)\ge1$ and 
$w(G_2\cap {\rm Int}D)\ge1$.

{\bf Claim.} $G_1=G_2$.

{\it Proof of Claim.}
Suppose $G_1\not=G_2$.
Then $G_1\cap G_2=\emptyset$.
Since $w(\Gamma_{m+\varepsilon}\cap {\rm Int}D)=3$,
we have $w(G_1\cap {\rm Int}D)=1$ or $w(G_2\cap {\rm Int}D)=1$.

Without loss of generality we can assume 
$w(G_1\cap {\rm Int}D)=1$.
Let $w$ be the white vertex in $G_1$ different from $w_1$.
Let $e$ be an internal edge (possibly a terminal edge)
of label $m+\varepsilon$ at $w$ different from $e_1'$
such that $e$ is not middle at $w$.
By Assumption~\ref{AssumeTerminal},
the edge $e$ is not a terminal edge.
Since $w(G_1\cap {\rm Int}D)=1$,
the edge $e$ is a loop.
Hence by Lemma~\ref{LemmaInsideLoopOutsideLoop}
the associated disk $E$ of the loop $e$
contains at least two white vertices in its interior.
Thus 
$$w(\Gamma\cap {\rm Int}D)\ge w(G_1\cap {\rm Int}D)+w(G_2\cap {\rm Int}D)+w(\Gamma\cap{\rm Int}E)\ge1+1+2=4.$$
This contradicts $w(\Gamma\cap {\rm Int}D)=3$.
Hence $G_1=G_2$.
Thus Claim holds. \hfill {$\square$}\vspace{1.5em}

By Claim,
there exists a simple arc $L$ 
in $\Gamma_{m+\varepsilon}\cap D$
connecting $w_1$ and $w_2$.
That is $e_1'\cup e_2'\subset L$.
Since $w(\Gamma\cap {\rm Int}D)=3$,
there are three cases:
(i) $w({\rm Int}L)=3$,
(ii) $w({\rm Int}L)=2$,
(iii) $w({\rm Int}L)=1$.

{\bf Case (i).}
Let $w_3,w_4,w_5$ be the white vertices in ${\rm Int}L$
with $w_3\in e_1'$ and $w_5\in e_2'$.
If all of the white vertices $w_3,w_4,w_5$
are BW-vertices with respect to $\Gamma_{m+\varepsilon}$,
then by Lemma~\ref{OriBWvertex}
a regular neighborhood of $D$ contains the pseudo chart as shown in
Fig.~\ref{fig10}(b).

If one of the white vertices $w_3,w_4,w_5$
is not a BW-vertex with respect to $\Gamma_{m+\varepsilon}$,
then there exists an internal edge $e$ of label $m+\varepsilon$
such that $e\cap L=\partial e$.
Hence $\partial e=\{w_3,w_4\}$, or 
$\partial e=\{w_3,w_5\}$, or $\partial e=\{w_4,w_5\}$.

If $\partial e=\{w_3,w_4\}$ or $\partial e=\{w_4,w_5\}$,
then there exists a 2-angled disk of $\Gamma_{m+\varepsilon}$
without feelers in ${\rm Int}D$.
Hence by Lemma~\ref{Theorem2AngledDisk} and Lemma~\ref{OriBWvertex},
a regular neighborhood of $D$ contains one of the RO-family of the pseudo chart as shown in 
Fig.~\ref{fig10}(c).

If $\partial e=\{w_3,w_5\}$,
then there exists a special 3-angled disk $E$ of $\Gamma_{m+\varepsilon}$ with at most one feeler in ${\rm Int}D$.
Hence by Lemma~\ref{Theorem3AngledDisk}(a)
a regular neighborhood of $D$ contains one of the RO-families of the two pseudo charts as shown in 
Fig~\ref{fig09}(a),(b).
Thus there exists an internal edge $e'$ of label $m$ or $m+2\varepsilon$ in $E$.
If $\partial e'=\{w_3,w_5\}$, then the edge $e'$ is 
oriented inward at both vertices $w_3$ and $w_5$
by Condition (ii) of this section.
This contradicts the definition of the chart.
Hence $\partial e'=\{w_3,w_4\}$ or $\partial e'=\{w_4,w_5\}$.
Therefore a regular neighborhood of $D$ contains one of the RO-family of the pseudo chart as shown in 
Fig.~\ref{fig10}(d).

{\bf Case (ii).}
Let $w_3$ be the white vertex in ${\rm Int}D$
but not contained in $L$.
Let $e_3,e_3'$ be internal edges (possibly terminal edges)
of label $m+\varepsilon$ at $w_3$
not middle at $w_3$.
By Assumption~\ref{AssumeTerminal},
neither $e_3$ nor $e_3'$ is a terminal edge.

We shall show that
neither $e_3$ nor $e_3'$ is a loop.
If $e_3$ or $e_3'$ is a loop,
then the loop bounds the disk $E$ in $D$.
By Lemma~\ref{LemmaInsideLoopOutsideLoop},
we have $w(\Gamma\cap{\rm Int }E)\ge 2$. 
Hence $w(\Gamma\cap{\rm Int }D)\ge 
w({\rm Int}L)+w(\Gamma\cap E)\ge 2+3=5$.
This contradicts $w(\Gamma\cap{\rm Int }D)=3$.
Thus neither $e_3$ nor $e_3'$ is a loop.

Therefore each of $e_3$ and $e_3'$ is
an internal edge with two white vertices.
Thus each of $e_3$ and $e_3'$ contains a white vertex in
${\rm Int}L$.
Hence the simple arc $\widetilde{L}=e_1'\cup e_3\cup e_3'\cup e_2'$
is connecting $w_1$ and $w_2$ with 
$w({\rm Int}\widetilde{L})=3$.
Thus Case (ii) follows from Case (i).

{\bf Case (iii).}
The arc $L$ is the union $e_1'\cup e_2'$.
Let $N$ be a regular neighborhood of $L\cup\partial D$.
Then $Cl(D-N)$ consists of two disks.
One of the two disks contains a white vertex 
of $\Gamma_{m+\varepsilon}$,
say $D'$.
Then $\Gamma_{m+\varepsilon}\cap\partial D'$ is at most one point.
Thus by Lemma~\ref{LemmaAtMostOnePoint},
we have $w(\Gamma\cap D')\ge2$.
Hence the condition $w(\Gamma\cap {\rm Int}D)=w(\Gamma_{m+\varepsilon}\cap {\rm Int}D)=3$
implies 
$w(\Gamma\cap D')=w(\Gamma_{m+\varepsilon}\cap D')=2$.
Thus by Lemma~\ref{LemmaTwoWhiteVertices},
the disk $D'$ contains one of the two pseudo charts
as shown in Fig.~\ref{fig03}.
Therefore
a regular neighborhood of $D$ contains one of RO-families of the two pseudo charts
as shown in 
Fig.~\ref{fig10}(e),(f).
We complete the proof of Lemma~\ref{Lemma2AngledDisks}$($d$)$.
{\hfill {$\square$}\vspace{1.5em}}

%%%%%%%%%%%%%%%%%%%%%%%%%%
%%%%%%%%%%%%%%%%%%%%%%%%%%
%%%%%%%%%%%%%%%%%%%%%%%%%%

%\newpage
\section{A 3-angled disk without feelers}

\label{s:3-angledDisksTwoWhiteNoFeelers}

In this section we investigate a 3-angled disk without feelers and with two white vertices in its interior.

\begin{figure}[htb]
\centerline{\includegraphics{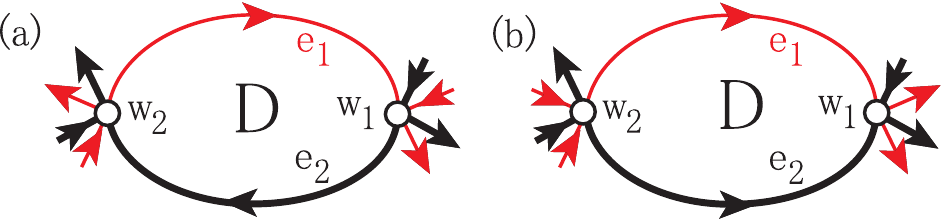}}
\caption{\label{fig12}
Lenses.}
\end{figure}

Let $\Gamma$ be a chart. 
Let $D$ be a disk 
such that 
\begin{enumerate}
\item[(1)] the boundary $\partial D$ consists of an internal edge $e_1$ of label $m$ and an internal edge $e_2$ of label ${m+1}$, and 
\item[(2)] any edge containing a white vertex in $e_1$ does not intersect the open disk Int$D$.
\end{enumerate}
Note that $\partial D$ may contain crossings.
Let $w_1$ and $w_2$ be the white vertices in $e_1$. 
If the disk $D$ satisfies one of the following conditions, then $D$ is called  {\it a lens of type $(m,m+1)$}
(see Fig.~\ref{fig12}):
\begin{enumerate}
	\item[(i)] Neither $e_1$ nor $e_2$ contains a middle arc. 
	\item[(ii)] One of the two edges $e_1$ and $e_2$ contains middle arcs at both white vertices $w_1$ and $w_2$ simultaneously.
\end{enumerate}

%%%%%%%%%%%%%%%%%% Figure
%%%%%%%%%%%%%%%%%%

\begin{lemma}{\rm (\cite[Theorem 1.1]{ChartApp1})}
\label{LensThreeWhiteVertex}
 There exist at least three white vertices
in the interior of the lens for any minimal chart.
\end{lemma}

\begin{lemma}{\rm (\cite[Corollary 1.3]{ChartAppII})}
\label{NoLens}
 There is no lens in any minimal chart with 
at most seven white vertices.
\end{lemma}

\begin{figure}
\centerline{\includegraphics{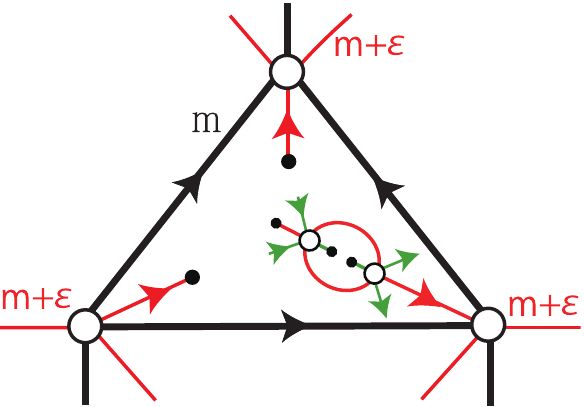}}
\caption{\label{fig13}
A 3-angeld disk without feelers where $m$ is a label, $\varepsilon\in\{+1,-1\}$.}
\end{figure}

\begin{lemma}
\label{3angledDiskNoFeelerTwoWhiteVertex}
Let $\Gamma$ be a minimal chart,
and $m$ a label of $\Gamma$.
Let $D$ be a $3$-angled disk of $\Gamma_m$
without feelers.
Suppose that all of the three white vertices $w_1,w_2,w_3$ on
 $\partial D$
are contained in $\Gamma_{m+\varepsilon}$ for
some $\varepsilon\in\{+1,-1\}$.
Let $e_1,e_2,e_3$ be internal edges 
$($possibly terminal edges$)$ of label $m+\varepsilon$
at $w_1,w_2,w_3$, respectively,
with $e_1\cup e_2\cup e_3\subset D$.
Suppose that $e_1,e_2$ are middle at $w_1,w_2$, respectively.
Then we have the following: 
\begin{enumerate}
\item[{\rm (a)}]
The edge $e_3$ is not middle at $w_3$.
\item[{\rm (b)}]
Suppose that $w(\Gamma\cap {\rm Int}D)=w(\Gamma_{m+\varepsilon}\cap {\rm Int}D)=2$.
If $e_1,e_3$ are oriented inward at $w_1,w_3$, respectively,
then either 
\begin{enumerate}
\item[{\rm (i)}] there exists a simple arc in $\Gamma_{m+\varepsilon}\cap D$ connecting $w_2$ and $w_3$, or
\item[{\rm (ii)}] a regular neighborhood of $D$ 
contains the pseudo chart
as shown in Fig.~\ref{fig13}.
\end{enumerate}
\end{enumerate}
\end{lemma}

%%%%%%%%%%%%%%%%%%
%%%%%%%%%%%%%%%%%% Figure
%%%%%%%%%%%%%%%%%%

\begin{Proof}
{\bf Statement (a).} 
Since $e_1,e_2$ are middle at $w_1,w_2$, respectively,
by Lemma~\ref{ROfamily3AngledDisk}
we can assume that $D$ is a 3-angled disk as shown in
Fig.~\ref{fig08}(a).
We use the orientations as shown in Fig.~\ref{fig08}(a).

Let $e',e''$ be internal edges of label $m$ at $w_3$ in
$\partial D$
with $e'\ni w_1$ and $e''\ni w_2$.
Then 
the edge $e'$ is oriented outward at $w_3$ and
$e''$ is oriented inward at $w_3$.
Thus the edge $e_3$ is not middle at $w_3$.

{\bf Statement (b).} 
We use the notations as shown in Fig.~\ref{fig08}(a).
It suffices to prove that 
if $\Gamma$ does not satisfy Condition (b)(i) of this lemma, 
then 
$\Gamma$ satisfies Condition (b)(ii).
Suppose that 
\begin{enumerate}
\item[$(*)$] there does not exist any simple arc 
in $\Gamma_{m+\varepsilon}\cap D$
connecting $w_2$ and $w_3$.
\end{enumerate}

We shall show that  the edge $e_3$ contains a white vertex in ${\rm Int}D$.
By Lemma~\ref{3angledDiskNoFeelerTwoWhiteVertex}(a)
and Assumption~\ref{AssumeTerminal},
the edge $e_3$ is not a terminal edge.
Since  by the condition of this lemma
the two edges $e_1,e_3$ are oriented inward at 
$w_1,w_3$, respectively,
we have $e_3\not=e_1$.
Moreover we have $e_3\not=e_2$ by Condition~$(*)$.
Thus the edge $e_3$ contains a white vertex in ${\rm Int}D$,
say $w_4$.

For the edge $e_2$,
there are two cases:
(i) $e_2$ is a terminal edge,
(ii) $e_2$ is not a terminal edge.

{\bf Case (i).}
Let $N$ be a regular neighborhood of $e_2\cup\partial D$
in $D$.
Set $D'=Cl(D-N)$.
Then $D'$ is an admissible disk 
such that $\Gamma_{m+\varepsilon}\cap\partial D$ is
exactly two points,
and $w(\Gamma\cap D')=w(\Gamma_{m+\varepsilon}\cap D')=2$.
Thus 
by Lemma~\ref{LemmaTwoWhiteVerticesWithBoundaryTwoPoints},
the disk $D'$ contains one of the four pseudo charts
as shown in Fig.~\ref{fig03} and 
Fig.~\ref{fig06}.

Since the edge $e_3$ contains the white vertex $w_4$ in
${\rm Int}D$,
the disk $D'$ contains one of the three pseudo charts
as shown in Fig.~\ref{fig03}(b) and 
Fig.~\ref{fig06}.
Moreover,
since the two edges $e_1$ and $e_3$ are oriented inward
at $w_1$ and $w_3$,
respectively,
the disk $D'$ contains the pseudo chart
as shown in Fig.~\ref{fig03}(b).
Therefore
the disk $D$ contains the pseudo chart as shown in
Fig.~\ref{fig13}.

{\bf Case (ii).}
Let $w_5$ be the white vertex in ${\rm Int}D$
different from $w_4$.
Since $e_2$ is not a terminal edge,
we have $e_2=e_1$ or $e_2\ni w_5$ by Condition~$(*)$.

If $e_2=e_1$, then there exists a lens $E$ in $D$ with $\partial E\supset e_2$.
Hence $w(\Gamma\cap{\rm Int }D)\ge3$ by Lemma~\ref{LensThreeWhiteVertex}.
This contradicts $w(\Gamma\cap{\rm Int}D)=2$. 
Thus $e_2\not=e_1$,
i.e.  $e_2\ni w_5$.

Let $G_1,G_2$ be the connected components of 
$\Gamma_{m+\varepsilon}\cap D$
with $G_1\supset e_2$ and $G_2\supset e_3$.
By Condition~$(*)$, we have $G_1\cap G_2=\emptyset$.
Since $e_2\ni w_5$ and $e_3\ni w_4$,
we have $w(G_1\cap{\rm Int}D)\ge1$ and $w(G_2\cap{\rm Int}D)\ge1$.
Thus $w(\Gamma\cap{\rm Int }D)=2$ implies that
$w(G_1\cap{\rm Int}D)=1$ and $w(G_2\cap{\rm Int}D)=1$.

One of $G_1,G_2$ does not contain the edge $e_1$,
say $G$.
Let $w$ be the white vertex in $G\cap {\rm Int}D$.
Let $e$ be an internal edge (possibly a terminal edge)
of label $m+\varepsilon$
at $w$ different from $e_2,e_3$ such that
$e$ is not middle at $w$.
By Assumption~\ref{AssumeTerminal},
the edge $e$ is not a terminal edge.
Thus $w(G\cap{\rm Int}D)=1$ implies that
the edge $e$ is a loop.
However this contradicts 
Lemma~\ref{LoopInKAngledDisk}. 
Therefore Case (ii) does not occur.
\end{Proof}

%%%%%%%%%%%%%%%%%%%%%%%%%%
%%%%%%%%%%%%%%%%%%%%%%%%%%
%%%%%%%%%%%%%%%%%%%%%%%%%%

%\newpage
\section{A special 3-angled disk with one feeler}

\label{s:3-angledDisksOneFeeler}

\begin{figure}
\centerline{\includegraphics{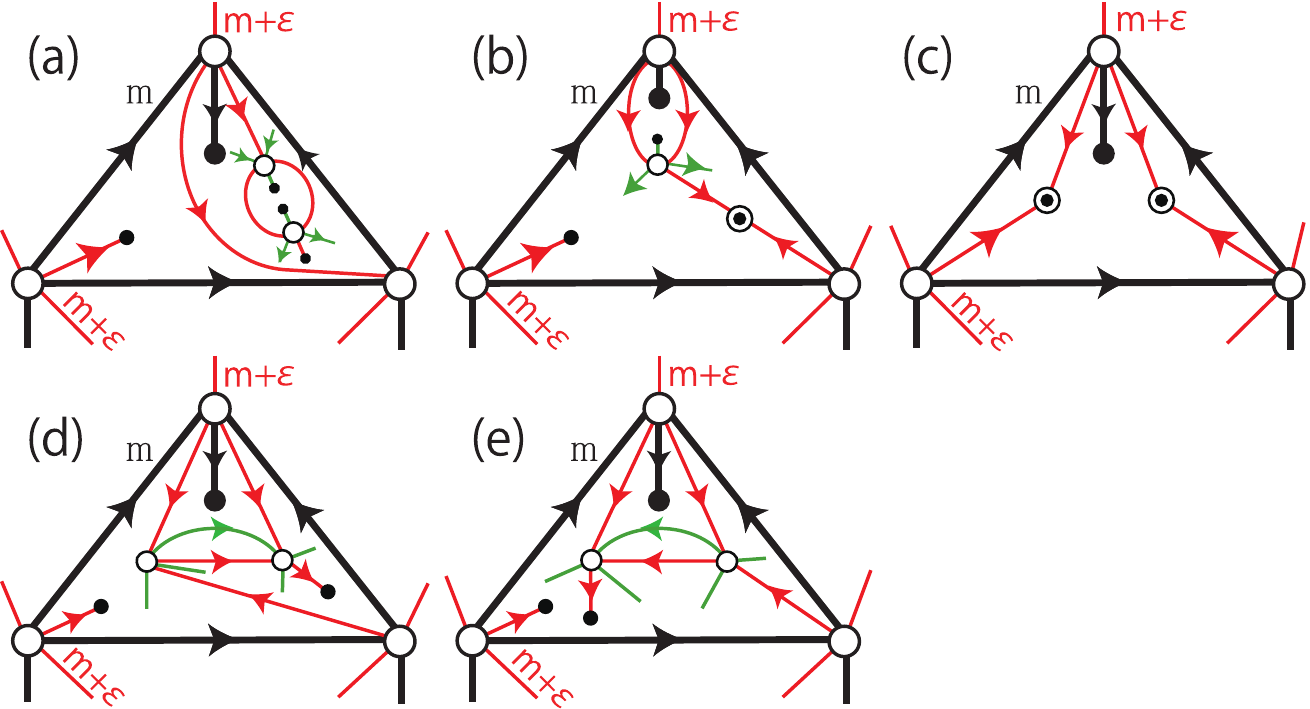}}
\caption{\label{fig14}
 3-angeld disks with one feeler and with two white vertices in their interiors.}
\end{figure}

In this section we investigate a special 3-angled disk with exactly one feeler and with at most three white vertices in its interior.
In this section we shall show the following two lemmata:

%%%%%%%%%%%%%%%%%%
%%%%%%%%%%%%%%%%%%

\begin{lemma}
\label{3angledDiskOneFeelerTwoWhiteVertex}
Let $\Gamma$ be a minimal chart,
and $m$ a label of $\Gamma$.
Let $D$ be a special $3$-angled disk of $\Gamma_m$
with exactly one feeler $e_1$,
and $w_1,w_2,w_3$ the three white vertices in 
 $\partial D$ with $w_1\in e_1$.
Suppose that all of $w_1,w_2,w_3$ 
are contained in $\Gamma_{m+\varepsilon}$ for
some $\varepsilon\in\{+1,-1\}$.
If $w(\Gamma\cap {\rm Int}D)=w(\Gamma_{m+\varepsilon}\cap {\rm Int}D)=2$,
then either 
\begin{enumerate}
\item[{\rm (a)}] there exists a simple arc in 
$(\Gamma_{m+\varepsilon}-w_1)\cap D$ 
connecting $w_2$ and $w_3$, or
\item[{\rm (b)}] a regular neighborhood of $D$ contains one of RO-families of the five pseudo charts
as shown in Fig.~\ref{fig14}.
\end{enumerate}
\end{lemma}

%%%%%%%%%%%%%%%%%%
%%%%%%%%%%%%%%%%%% Figure
%%%%%%%%%%%%%%%%%%

%%%%%%%%%%%%%%%%%%
%%%%%%%%%%%%%%%%%%

\begin{figure}
\centerline{\includegraphics{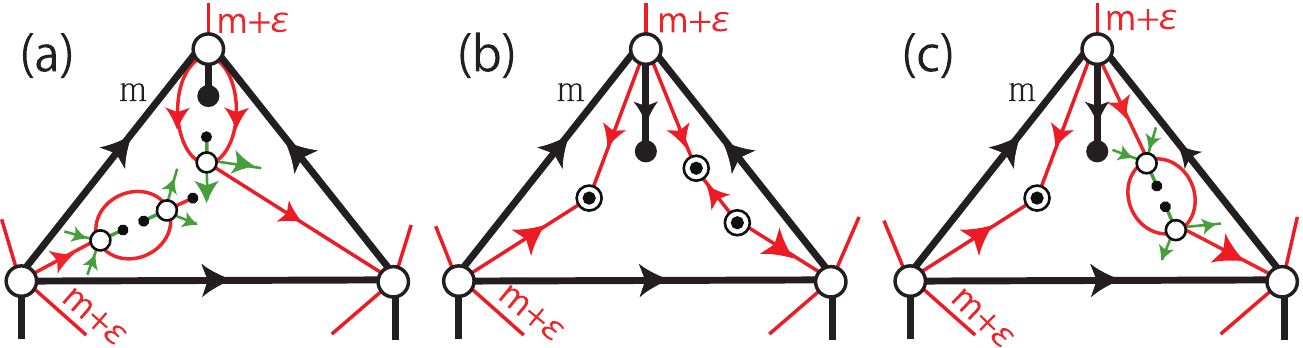}}
\caption{\label{fig15}
 3-angeld disks with one feeler and three white vertices in their interiors.}
\end{figure}

\begin{lemma}
\label{3angledDiskOneFeelerThreeWhiteVertex}
Let $\Gamma$ be a minimal chart,
and $m$ a label of $\Gamma$.
Let $D$ be a special $3$-angled disk of $\Gamma_m$
with exactly one feeler $e_1$,
and $w_1,w_2,w_3$ the three white vertices in 
 $\partial D$ with $w_1\in e_1$.
Suppose that all of $w_1,w_2,w_3$ 
are contained in $\Gamma_{m+\varepsilon}$ for
some $\varepsilon\in\{+1,-1\}$.
If $w(\Gamma\cap {\rm Int}D)=w(\Gamma_{m+\varepsilon}\cap {\rm Int}D)=3$ and 
if there exists an internal edge $e_2$
of label $m+\varepsilon$ at $w_2$ in $D$ such that 
$e_2$ is middle at $w_2$ but not a terminal edge,
then either 
\begin{enumerate}
\item[{\rm (a)}] there exists a simple arc in 
$(\Gamma_{m+\varepsilon}-w_1)\cap D$ 
connecting $w_2$ and $w_3$, or
\item[{\rm (b)}]  a regular neighborhood of $D$ contains one of RO-families of the three pseudo charts
as shown in Fig.~\ref{fig15}.
\end{enumerate}
\end{lemma}

%%%%%%%%%%%%%%%%%%
%%%%%%%%%%%%%%%%%% Figure
%%%%%%%%%%%%%%%%%%

From now on throughout this section,
we use the notations in Lemma~\ref{3angledDiskOneFeelerTwoWhiteVertex}.
Moreover,
by Lemma~\ref{ROfamily3AngledDisk}
we may assume that
\begin{enumerate}
\item[(i)] $D$ is a 3-angled disk as shown in Fig.~\ref{fig08}(c)
(we use the notations as shown in Fig.~\ref{fig08}(c)),
\item[(ii)] the feeler $e_1$ is oriented outward at $w_1$,
\item[(iii)]
the internal edge of label $m$ connecting $w_2$ and $w_3$
is oriented from $w_2$ to $w_3$.
\end{enumerate}

We can prove the following lemma easily: 
\begin{lemma}
\label{Remark3angledDiskOneFeeler}
Let $\Gamma,D$ be as above. 
Let $a_{11},b_{11},e_2,e_3$ be internal edges $($possibly terminal edges$)$ of label $m+\varepsilon$ at $w_1,w_1,w_2,w_3$ in $D$,
respectively
$($see Fig.~\ref{fig08}$($c$))$.
Then we have the following:
\begin{enumerate}
\item[{\rm (a)}]
Both of $a_{11},b_{11}$ are oriented outward at $w_1$,
and neither $a_{11}$ nor $b_{11}$ is a terminal edge.
\item[{\rm (b)}] The edge $e_2$ is oriented outward at $w_2$
 and middle at $w_2$. 
\item[{\rm (c)}]
The edge $e_3$ is not middle at $w_3$
$($i.e. $e_3$ is not a terminal edge$)$.
\end{enumerate}
\end{lemma}

\begin{lemma}
\label{a11WhiteVertex}
Let $\Gamma,D,b_{11}$ be as above.
If $w(\Gamma\cap {\rm Int}D)\le3$,
then the edge $b_{11}$ contains a white vertex in ${\rm Int}D$.
\end{lemma}

\begin{Proof}
Suppose that the edge $b_{11}$ does not contain 
any white vertex in ${\rm Int}D$.
Since the edge $b_{11}$ is not a terminal edge
by Lemma~\ref{Remark3angledDiskOneFeeler}(a),
we have $b_{11}=a_{11}$ or $b_{11}=e_2$ or $b_{11}=e_3$.
Since the edges $a_{11},b_{11},e_2$ are oriented
outward at $w_1,w_1,w_2$, respectively
by Lemma~\ref{Remark3angledDiskOneFeeler}(a) and (b),
we have $b_{11}=e_3$.
Hence there exists a lens $E$ in $D$ with 
$\partial E\supset b_{11}$.
Thus by Lemma~\ref{LensThreeWhiteVertex},
we have $w(\Gamma\cap{\rm Int}E)\ge3$.

By Lemma~\ref{Remark3angledDiskOneFeeler}(a),
the edge $a_{11}$ is not a terminal edge.
Moreover $a_{11}\not=e_2$ by Lemma~\ref{Remark3angledDiskOneFeeler}(b).
Thus $a_{11}$ contains a white vertex in ${\rm Int}D-E$, 
say $w_4$.
Hence $w(\Gamma\cap{\rm Int}D)\ge w(w_4)+w(\Gamma\cap{\rm Int}E)\ge1+3=4$.
This contradicts $w(\Gamma\cap {\rm Int}D)\le3$.
Therefore the edge $b_{11}$ contains a white vertex in ${\rm Int}D$.
\end{Proof}

%%%%%%%%%%%%%%%%%%%
%%%%%%%%%%%%%%%%%%%
%%%%%%%%%%%%%%%%%%%

{\it Proof of Lemma~\ref{3angledDiskOneFeelerTwoWhiteVertex}.}
It suffices to prove that 
if $\Gamma$ does not satisfy Condition~(a) of this lemma, 
then 
$\Gamma$ satisfies Condition~(b).
Suppose that 
\begin{enumerate}
\item[$(*)$] there does not exist any simple arc in 
$(\Gamma_{m+\varepsilon}-w_1)\cap D$ 
connecting $w_2$ and $w_3$.
\end{enumerate}

By Lemma~\ref{a11WhiteVertex},
 the edge $b_{11}$ contains a white vertex
in ${\rm Int}D$, say $w_4$.
Let $w_5$ be a white vertex in ${\rm Int}D$
different from $w_4$.

Now the edge $a_{11}$ is not a terminal edge by
Lemma~\ref{Remark3angledDiskOneFeeler}(a).
Moreover by Lemma~\ref{Remark3angledDiskOneFeeler}(b)
we have $a_{11}\not=e_2$.
Thus for the edge $a_{11}$,
there are three cases:
(i) $a_{11}=e_3$,
(ii) $a_{11}\ni w_4$,
(iii) $a_{11}\ni w_5$.

{\bf Case (i).}
Let $N$ be a regular neighborhood of 
$a_{11}\cup e_1\cup \partial D$ in $D$.
Then $Cl(D-N)$ consists of two disks.
One of the two disks contains the white vertex $w_4$,
say $D'$.
Then $D'$ is an admissible disk 
such that $\Gamma_{m+\varepsilon}\cap\partial D'$ is one point.
Since $w(\Gamma_{m+\varepsilon}\cap D')\ge1$,
we have $w(\Gamma\cap D')\ge2$ by Lemma~\ref{LemmaAtMostOnePoint}.
Hence the condition $w(\Gamma\cap {\rm Int}D)=w(\Gamma_{m+\varepsilon}\cap {\rm Int}D)=2$
implies that $w(\Gamma\cap D')=w(\Gamma_{m+\varepsilon}\cap D')=2$.
Thus by Lemma~\ref{LemmaTwoWhiteVertices},
the disk $D'$ contains the pseudo chart as shown
in Fig.~\ref{fig03}(b).
Hence the edge $e_2$ is a terminal edge,
and a regular neighborhood of $D$ contains the pseudo chart 
as shown in Fig.~\ref{fig14}(a).

{\bf Case (ii).}
Let $e_5,e_5'$ be internal edges (possibly terminal edges)
of label $m+\varepsilon$ at $w_5$
but not middle at $w_5$.
By Assumption~\ref{AssumeTerminal},
 neither $e_5$ nor $e_5'$ is a terminal edge.
Moreover, by Lemma~\ref{LoopInKAngledDisk},
 neither $e_5$ nor $e_5'$ is a loop.
Thus by Condition~$(*)$, we have either
($w_2,w_4\in e_5\cup e_5'$) or ($w_3,w_4\in e_5\cup e_5'$).

If $w_2,w_4\in e_5\cup e_5'$,
then the edge $e_3$ is a terminal edge by Condition~$(*)$.
This contradicts Lemma~\ref{Remark3angledDiskOneFeeler}(c).
Hence $w_3,w_4\in e_5\cup e_5'$.
Therefore the edge $e_2$ is a terminal edge by Condition~$(*)$.
Moreover by Lemma~\ref{OriBWvertex}
a regular neighborhood of $D$ contains the pseudo chart 
as shown in Fig.~\ref{fig14}(b).

{\bf Case (iii).}
Let $e_4,e_5$ be internal edges (possibly terminal edges)
of label $m+\varepsilon$ at $w_4,w_5$
different from $a_{11},b_{11}$, respectively, such that 
neither $e_4$ nor $e_5$ is middle at $w_4,w_5$.
By Assumption~\ref{AssumeTerminal},
neither $e_4$ nor $e_5$ is a terminal edge.
Moreover by Lemma~\ref{LoopInKAngledDisk},
neither $e_4$ nor $e_5$ is a loop.

For the edge $e_5$,
there are three cases:
$e_5\ni w_2$, $e_5\ni w_3$ or $e_5\ni w_4$.

If $e_5\ni w_2$ (i.e. $e_5=e_2$),
then $e_3\not\ni w_5$ by Condition~$(*)$.
Hence by Lemma~\ref{Remark3angledDiskOneFeeler}(c),
 we have $e_3\ni w_4$.
Moreover by Condition~$(*)$,
both of $w_4$ and $w_5$ are BW-vertices
with respect to $\Gamma_{m+\varepsilon}$.
Thus by Lemma~\ref{OriBWvertex}
a regular neighborhood of $D$ contains the pseudo chart 
as shown in Fig.~\ref{fig14}(c).

If $e_5\ni w_3$ (i.e. $e_5=e_3$),
then the edge $e_2$ is a terminal edge by Condition~$(*)$.
Since $e_4$ is not a loop,
we have $e_4\ni w_5$.
Thus by Lemma~\ref{Theorem3AngledDisk}(a),
a regular neighborhood of $D$ contains the pseudo chart 
as shown in Fig.~\ref{fig14}(d).

If $e_5\ni w_4$,
then $e_3\ni w_4$ or $e_3\ni w_5$ 
by Condition~$(*)$ and Lemma~\ref{Remark3angledDiskOneFeeler}(c).
If $e_3\ni w_4$,
then the edge $e_2$ is a terminal edge by Condition~$(*)$.
Thus  by Lemma~\ref{Theorem3AngledDisk}(a),
a regular neighborhood of $D$ contains the pseudo chart 
as shown in Fig.~\ref{fig14}(e). 
If $e_3\ni w_5$,
then $e_2$ is a terminal edge.
Thus by Lemma~\ref{Theorem3AngledDisk}(a), 
a regular neighborhood of $D$ contains the pseudo chart 
as shown in Fig.~\ref{fig14}(d).
We complete the proof of Lemma~\ref{3angledDiskOneFeelerTwoWhiteVertex}.
{\hfill {$\square$}\vspace{1.5em}}

%%%%%%%%%%%%%%%%%%%%%%
%%%%%%%%%%%%%%%%%%%%%%
%%%%%%%%%%%%%%%%%%%%%%

{\it Proof of Lemma~\ref{3angledDiskOneFeelerThreeWhiteVertex}.}
It suffices to prove that 
if $\Gamma$ does not satisfy Condition~(a), then 
$\Gamma$ satisfies Condition~(b).
Suppose that 
\begin{enumerate}
\item[$(*)$] there does not exist any simple arc in 
$(\Gamma_{m+\varepsilon}-w_1)\cap D$ 
connecting $w_2$ and $w_3$.
\end{enumerate}

{\bf Claim 1.} All of $e_2,e_3,a_{11},b_{11}$ contain
white vertices in ${\rm Int}D$.

{\it Proof of Claim $1$.}
By Lemma~\ref{a11WhiteVertex},
 the edge $b_{11}$ contains a white vertex
in ${\rm Int}D$, say $w'$. 

We shall show that $e_2$ contains
 a white vertex in ${\rm Int}D$.
First we have $e_2\not=e_3$ by Condition~$(*)$.
Moreover
we have $e_2\not=a_{11}$ and $e_2\not=b_{11}$ 
by Lemma~\ref{Remark3angledDiskOneFeeler}(a) and (b)
(because $e_2,a_{11},b_{11}$ are oriented outward at $w_2,w_1,w_1$, respectively).
Furthermore, $e_2$ is not a terminal edge 
by the condition of this lemma.
Hence the edge $e_2$ must contain a white vertex 
in ${\rm Int}D$, say $w''$.

We shall show $e_3\not=a_{11}$.
If $e_3=a_{11}$, then $e_3$ is a proper arc in $D$.
Let $N$ be a regular neighborhood of $e_1\cup e_3\cup\partial D$ in $D$.
Then $Cl(D-N)$ consists of two disks.
One of the two disks contains the white vertex $w'$,
say $D_1$. The other contains $w''$, say $D_2$.
Thus we have $w(\Gamma_{m+\varepsilon}\cap{\rm Int}D_i)\ge1$
for $i=1,2$.
Since the disk $D_i$ is an admissible disk for $i=1,2$ 
such that $\Gamma_{m+\varepsilon}\cap\partial D_i=$one point,
we have $w(\Gamma\cap{\rm Int}D_i)\ge2$
by Lemma~\ref{LemmaAtMostOnePoint}.
Hence $w(\Gamma\cap{\rm Int}D)\ge4$.
This contradicts $w(\Gamma\cap{\rm Int}D)=3$.
Thus $e_3\not=a_{11}$.

Finally we shall show that $e_3,a_{11}$ contain
 white vertices in ${\rm Int}D$.
Since $e_3\not=a_{11}$ and
since $b_{11},e_2$ contains white vertices in ${\rm Int}D$,
the four edges $e_3,a_{11},b_{11},e_2$ are different edges each other.
Since neither  $e_3$ nor $a_{11}$ is a terminal edge by Lemma~\ref{Remark3angledDiskOneFeeler}(a) and (c),
the edges $e_3,a_{11}$ contain white vertices in ${\rm Int}D$. 
Thus Claim~1 holds.
{\hfill {$\square$}\vspace{1.5em}}

\begin{figure}
\centerline{\includegraphics{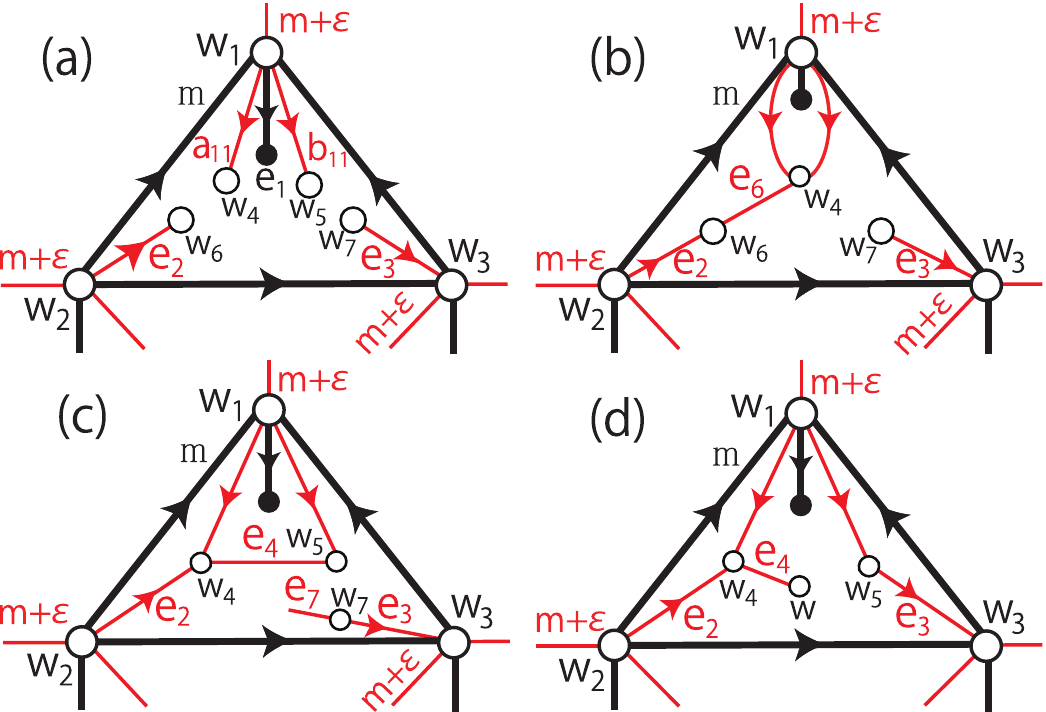}}
\caption{\label{fig16}
 3-angeld disks with one feeler and three white vertices in their interiors.}
\end{figure}

Let $w_4,w_5,w_6,w_7$ be white vertices in ${\rm Int}D$
with $w_4\in a_{11}$,  $w_5\in b_{11}$,
$w_6\in e_2$ and $w_7\in e_3$ (see Fig.~\ref{fig16}(a)).
It is possible that 
some of $w_4,w_5,w_6,w_7$ are the same vertices.
By Condition~$(*)$,
we have  
\begin{enumerate}
\item[(1)]
$w_6\not=w_7$.
\end{enumerate}

{\bf Claim 2.} None of $w_4,w_5,w_6,w_7$ are contained in any loop of label $m+\varepsilon$.

{\it Proof of Claim~$2$.} 
If one of $w_4,w_5,w_6,w_7$ is contained in a loop $\ell$ 
of label $m+\varepsilon$,
then by Lemma~\ref{LemmaInsideLoopOutsideLoop}
the loop $\ell$ bounds a disk $E$ in $D$
with $w(\Gamma\cap{\rm Int}E)\ge2$.
Hence $w(\Gamma\cap{\rm Int}D)\ge w(\{w_6,w_7\})+w(\Gamma\cap{\rm Int}E)\ge2+2=4$.
This contradicts $w(\Gamma\cap{\rm Int}D)=3$.
Thus Claim~2 holds.
{\hfill {$\square$}\vspace{1.5em}}

{\bf Claim 3.} $w_5\not=w_6$.

{\it Proof of Claim~$3$.}
Suppose $w_5=w_6$.
First we shall show $w_4\not=w_5$.
If $w_4=w_5$, then $w_4=w_5=w_6$.
Thus three edges $a_{11},b_{11},e_2$
are edges of label $m+\varepsilon$ oriented inward at $w_4$.
This contradicts the definition of a chart.
Hence $w_4\not=w_5$.

Let $N$ be a regular neighborhood of 
$e_1\cup e_2\cup b_{11}\cup \partial D$ in $D$.
Then $Cl(D-N)$ consists of two disks.
One of the two disks contains the white vertex $w_4$,
say $D_1$.
By (1),
the other contains the white vertex $w_7$, say $D_2$.
Thus $w(\Gamma_{m+\varepsilon}\cap {\rm Int}D_1)\ge1$ and 
 $w(\Gamma_{m+\varepsilon}\cap {\rm Int}D_2)\ge1$.
Moreover one of $\partial D_1$ and $\partial D_2$ 
contains exactly one point of $\Gamma_{m+\varepsilon}$.
Hence by Lemma~\ref{LemmaAtMostOnePoint},
we have $w(\Gamma\cap {\rm Int}D_1)\ge2$ or 
$w(\Gamma\cap {\rm Int}D_2)\ge2$.
Thus $w(\Gamma\cap {\rm Int}D)\ge w(w_5)+w(\Gamma\cap {\rm Int}D_1)+w(\Gamma\cap {\rm Int}D_2)\ge 1+3=4$.
This contradicts $w(\Gamma\cap{\rm Int}D)=3$.
Hence $w_5\not=w_6$. Thus Claim~3 holds.
{\hfill {$\square$}\vspace{1.5em}}

Since $w(\Gamma\cap{\rm Int}D)=3$,
two of $w_4,w_5,w_6,w_7$
are the same vertices.
Thus by (1) and Claim~3,
there are four cases:
(i) $w_4=w_5$,
(ii) $w_4=w_6$, (iii) $w_4=w_7$,
(iv) $w_5=w_7$.

{\bf Case (i).}
We shall show $w_4=w_7$.
If $w_4\not=w_7$,
then by (1) and Claim~3
the three white vertices $w_4(=w_5),w_6,w_7$ are different each other. 
Let $e_6$ be an internal edge (possibly a terminal edge)
of label $m+\varepsilon$ at $w_6$
different from $e_2$
such that  $e_6$ 
is not  middle at $w_6$.
Then by Assumption~\ref{AssumeTerminal},
the edge $e_6$ is not a terminal edge.
Since $e_6$ is not a loop by Claim~2, 
we have $e_6\in w_4$ by Condition~$(*)$
(see Fig.~\ref{fig16}(b)).

Let $e_7$ be an internal edge (possibly a terminal edge)
of label $m+\varepsilon$ at $w_7$
different from $e_3$
such that  $e_7$ 
is not middle at $w_7$.
Then by Assumption~\ref{AssumeTerminal},
the edge $e_7$ is not a terminal edge.
Thus by Condition~$(*)$,
the edge $e_7$ must be a loop.
However this contradicts Claim~2.
Therefore $w_4=w_7$.

Let $N$ be a regular neighborhood of 
$e_1\cup e_3\cup a_{11}\cup \partial D$ in $D$.
Then $Cl(D-N)$ consists of two disks.
One of the two disks contains the white vertex
$w_6$, say $D'$.
Then $D'$ is an admissible disk such that
$\Gamma_{m+\varepsilon}\cap\partial D'=$one point.
Thus by Lemma~\ref{LemmaAtMostOnePoint},
we have $w(\Gamma\cap D')\ge2$.
Since $w(\Gamma\cap {\rm Int}D)=w(\Gamma_{m+\varepsilon}\cap {\rm Int}D)=3$,
we have 
$w(\Gamma\cap D')=w(\Gamma_{m+\varepsilon}\cap D')=2$.
Hence by Lemma~\ref{LemmaTwoWhiteVertices},
the disk $D'$ contains the pseudo chart as shown 
in Fig.~\ref{fig03}(b).
Thus by Lemma~\ref{Theorem2AngledDisk},
a regular neighborhood of  $D$ contains the pseudo chart 
as shown in Fig.~\ref{fig15}(a).

{\bf Case (ii).}
Let $e_4$ be an internal edge (possibly a terminal edge)
 of label $m+\varepsilon$ at $w_4$
different from $e_2,a_{11}$.

We shall show that $e_4$ is a terminal edge.
If $e_4$ is not a terminal edge,
then $e_4$ contains a white vertex in ${\rm Int}D$
different from $w_4$,
say $w$.
There are two cases:
$w=w_5$ or $w\not=w_5$.

If $w=w_5$, then by Condition~$(*)$
the three white vertices 
$w_4,w(=w_5),w_7$ are different each other
(see Fig.~\ref{fig16}(c)).
Let $e_7$ be an internal edge (possibly a terminal edge)
 of label $m+\varepsilon$ at $w_7$
different from $e_3$
such that $e_7$ is not middle at $w_7$.
By Assumption~\ref{AssumeTerminal},
the edge $e_7$ is not a terminal edge.
Moreover we have $e_7\not\ni w_5$ by Condition~$(*)$.
Thus the edge $e_7$ must be a loop.
This contradicts Claim~2.

If $w\not=w_5$,
then by Claim~3 the three white vertices 
$w_4(=w_6),w,w_5$ are different each other.
Since $w(\Gamma\cap{\rm Int}D)=3$,
we have $w_5=w_7$ by Condition~$(*)$
(see Fig.~\ref{fig16}(d)).
By the similar way of the above,
there exists a loop of label $m+\varepsilon$
containing $w$.
This contradicts Claim~2.
Therefore $e_4$ is a terminal edge.

Let $N$ be a regular neighborhood of 
$e_1\cup e_2\cup e_4\cup a_{11}\cup \partial D$ in $D$.
Then $Cl(D-N)$ consists of two disks.
One of the two disks contains the white vertex
$w_5$, say $D_1$.
Let $D_2$ be the other disk.
Then $\Gamma_{m+\varepsilon}\cap \partial D_1=$two points,
and $\Gamma_{m+\varepsilon}\cap \partial D_2=\emptyset$.

We shall show $w(\Gamma_{m+\varepsilon}\cap {\rm Int}D_2)=0$.
If $w(\Gamma_{m+\varepsilon}\cap {\rm Int}D_2)\ge1$,
then by Lemma~\ref{LemmaAtMostOnePoint}
we have $w(\Gamma\cap{\rm Int}D_2)\ge 2$.
Thus 
$w(\Gamma\cap {\rm Int}D)\ge w(w_4)+w(\Gamma\cap {\rm Int}D_1)+w(\Gamma\cap {\rm Int}D_2)\ge 1+1+2=4$.
This contradicts $w(\Gamma\cap{\rm Int}D)=3$.
Hence $w(\Gamma_{m+\varepsilon}\cap {\rm Int}D_2)=0$.

Therefore 
the condition $w(\Gamma\cap {\rm Int}D)=w(\Gamma_{m+\varepsilon}\cap {\rm Int}D)=3$
implies that
$w(\Gamma\cap {\rm Int}D_2)=0$
and $w(\Gamma\cap {\rm Int}D_1)=w(\Gamma_{m+\varepsilon}\cap {\rm Int}D_1)=2$.
Thus by Lemma~\ref{LemmaTwoWhiteVerticesWithBoundaryTwoPoints},
a regular neighborhood of $D$ contains one of the two pseudo charts 
as shown in Fig.~\ref{fig15}(b),(c).

{\bf Case (iii).}
We shall show $w_4=w_5$.
If $w_4\not=w_5$,
then we can show $w(\Gamma\cap{\rm Int}D)\ge4$ 
by the similar way of the proof of Claim~3.
This contradicts  $w(\Gamma\cap{\rm Int}D)=3$.
Thus $w_4=w_5$.
Hence Case (iii) follows from Case (i).

{\bf Case (iv).}
We shall show $w_4=w_5$.
If $w_4\not=w_5$,
then by the similar way of the proof of Case (ii)
we can show that $w_5$ is contained in a terminal edge
of label $m+\varepsilon$.
Moreover we can show that 
there exists an admissible disk $D'$ in $D$
with $w_4\in D'$, 
$\Gamma_{m+\varepsilon}\cap\partial D'=(e_2\cup a_{11})\cap\partial D'=$two points 
and $w(\Gamma\cap D')=w(\Gamma_{m+\varepsilon}\cap D')=2$.
Thus by Lemma~\ref{LemmaTwoWhiteVerticesWithBoundaryTwoPoints},
the disk $D'$ contains one of the two pseudo charts
as shown in Fig.~\ref{fig06}.
Hence either $e_2$ is oriented inward at $w_3$,
or $a_{11}$ is oriented inward at $w_1$.
However this contradicts
Lemma~\ref{Remark3angledDiskOneFeeler}(a),(b).
Therefore $w_4=w_5$.
Hence Case (iv) follows from Case (i).

Therefore we complete the proof of Lemma~\ref{3angledDiskOneFeelerThreeWhiteVertex}.
{\hfill {$\square$}\vspace{1.5em}}

%%%%%%%%%%%%%%%%%%
%%%%%%%%%%%%%%%%%% Figure
%%%%%%%%%%%%%%%%%%

%%%%%%%%%%%%%%%%%%%%%%%%
%%%%%%%%%%%%%%%%%%%%%%%%
%\newpage

\section{IO-Calculation}

\label{s:IOC}

In this section,
we review IO-Calculation.

Let $\Gamma$ be a chart,
 and $v$ a vertex. 
Let $\alpha$ be a short arc of $\Gamma$ in a small neighborhood of $v$ such that $v$ is an endpoint of $\alpha$. 
If the arc $\alpha$ is oriented to $v$, then $\alpha$ is called {\it an inward arc}, 
and otherwise $\alpha$ is called {\it an outward arc}.

\begin{figure}
\centerline{\includegraphics{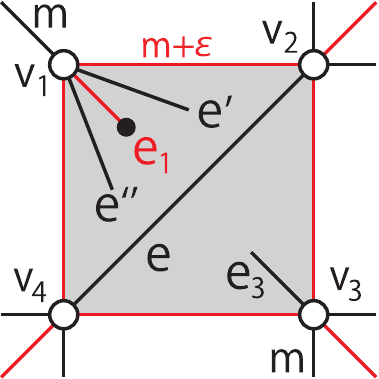}}
\caption{\label{fig17} The gray region is the 4-angled disk $F$, $m$ is a label, $\varepsilon\in\{+1,-1\}$.}
\end{figure}

Let $\Gamma$ be an $n$-chart. 
Let $F$ be a closed domain with $\partial F\subset \Gamma_{m-1}\cup\Gamma_{m}\cup \Gamma_{m+1}$ for some label $m$ of $\Gamma$, where $\Gamma_0=\emptyset$ and $\Gamma_{n}=\emptyset$. 
By Condition (iii) for charts,
in a small neighborhood of each white vertex, there are three inward arcs and three outward arcs.
Also in a small neighborhood of each black vertex, there exists only one inward arc or one outward arc.
We often use the following fact, 
when we fix (inward or outward) arcs 
near white vertices and black vertices: 
\begin{enumerate}
\item[$(*)$]
{\it The number of inward arcs contained in $F\cap \Gamma_m$ is equal to the number of outward arcs in $F\cap \Gamma_m$.
}
\end{enumerate}
When we use this fact, 
we say that we use {\it IO-Calculation with respect to $\Gamma_m$ in $F$}.
For example, in a minimal chart $\Gamma$, 
consider the pseudo chart as shown in Fig.~\ref{fig17} 
where
\begin{enumerate}
\item[(1)] $F$ is a $4$-angled disk of $\Gamma_{m+\varepsilon}$ with exactly one feeler $e_1$ for some $\varepsilon\in\{+1,-1\}$,
\item[(2)]  $v_1,v_2,v_3,v_4$ are white vertices in $\partial F$ with $v_1\in e_1$ and $v_1,v_3\in \Gamma_m$,
\item[(3)] $e$ is an internal edge of label $m$ or $m+2\varepsilon$ with $\partial e=\{v_2,v_4\}$ and $e\subset F$, 
\item[(4)] $e',e'',e_3$ are internal edges 
(possibly terminal edges) of label $m$ 
at $v_1,v_1,v_3$ in $F$, respectively.
\end{enumerate}
Then we can show that $w(\Gamma\cap{\rm Int}F)\ge1$.
Suppose $w(\Gamma\cap{\rm Int}F)=0$.
Since $e_1$ is a terminal edge,
the edge $e_1$ is middle at $v_1$ 
by Assumption~\ref{AssumeTerminal}.
Thus neither $e'$ nor $e''$ is middle at $v_1$.
Hence by Assumption~\ref{AssumeTerminal}
\begin{enumerate}
\item[(5)] neither  $e'$ nor $e''$ is a terminal edge. 
\end{enumerate}
If necessary we change orientations of all the edges,
we can assume that $e_1$ is oriented inward at $v_1$.
Thus  
\begin{enumerate}
\item[(6)] both of $e'$ and $e''$ are oriented inward at $v_1$. 
\end{enumerate}

If $e_3$ is a terminal edge oriented inward at $v_3$
and if $e$ is of label $m$,
then
by (5) and (6) 
the number of inward arcs in $F\cap \Gamma_m$ is four,  
but the number of outward arcs in $F\cap \Gamma_m$ is two. 
This contradicts the fact $(*)$. 
If $e_3$ is oriented inward at $v_3$ but not a terminal edge
and if $e$ is of label $m$,
then
by (5) and (6) 
the number of inward arcs in $F\cap \Gamma_m$ is four,  
but the number of outward arcs in $F\cap \Gamma_m$ is one. 
This contradicts the fact~$(*)$. 
Similarly for the others cases
 we have the same contradiction.
Thus $w(\Gamma\cap{\rm Int}F)\ge1$.
Instead of the above argument, 
\begin{enumerate}
\item[]
{\it we have $w(\Gamma\cap{\rm Int}F)\ge1$ 
by IO-Calculation with respect to $\Gamma_{m}$ in $F$.}
\end{enumerate}

%%%%%%%%%%%%%%%%%%
%%%%%%%%%%%%%%%%%% Figure
%%%%%%%%%%%%%%%%%%

%%%%%%%%%%%%%%%%%%%%%%%%%%
%%%%%%%%%%%%%%%%%%%%%%%%%%
%%%%%%%%%%%%%%%%%%%%%%%%%%

%\newpage
\section{Disk Lemma}

\label{s:DiskLemma}

In this section we review New Disk Lemma.
We shall give a condition of non-minimal charts.

\begin{figure}[htb]
\centerline{\includegraphics{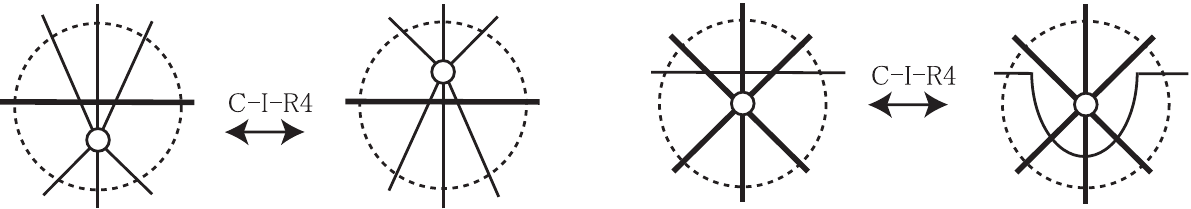}}
\caption{\label{fig18} 
C-moves keeping thicken figures fixed.}
\end{figure}

Let $\Gamma$ and $\Gamma^\prime $ be C-move equivalent charts. 
Suppose that a pseudo chart $X$ of $\Gamma$ is also a pseudo chart of $\Gamma^\prime$. 
Then we say that 
$\Gamma$ is modified to $\Gamma^\prime$ by {\it C-moves keeping $X$ fixed}.
In Fig.~\ref{fig18},
we give examples of C-moves keeping pseudo charts  fixed.

%%%%%%%%%%%%%%%%%%
%%%%%%%%%%%%%%%%%% Figure
%%%%%%%%%%%%%%%%%%

Let $\Gamma$ be a chart, and $D$ a disk.
Let $\alpha$ be a simple arc in $\partial D$,
and $\gamma$ a simple arc in an internal edge of label $k$.
The simple arc $\gamma$ is called
a {\it {$(D,\alpha)$-arc}} of label $k$
provided that 
$\partial \gamma \subset $Int$\alpha$
and
Int$\gamma\subset $Int$D$. 
If there is no $(D,\alpha)$-arc in $\Gamma$,
then the chart $\Gamma$ is said to be
$(D,\alpha)$-{\it arc free}.

\begin{figure}
\centerline{\includegraphics{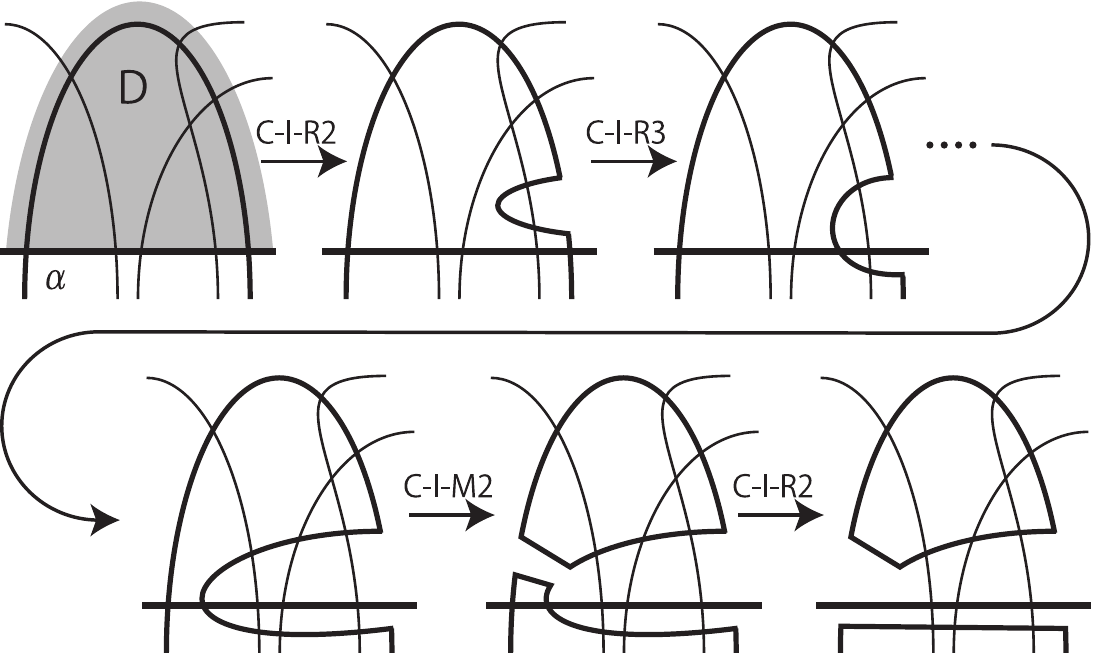}}
\caption{\label{fig19}
The gray region is the disk $D$.}
\end{figure}

\begin{lemma}
$($New Disk Lemma$)$
{\em (\cite[Lemma 7.1]{Chart33},
cf. \cite[Lemma 3.2]{ChartApp1})} 
\label{NewDiskLemma}
Let $\Gamma$ be a chart and
$D$ a disk 
whose interior does not contain 
a white vertex nor a black vertex of $\Gamma$.
Let $\alpha$ be a simple arc in $\partial D$ 
such that ${\rm Int}\alpha$ does not contain 
a white vertex nor a black vertex of $\Gamma$.
Let $V$ be a regular neighborhood of $\alpha$. 
Suppose that 
the arc $\alpha$ is contained in 
an internal edge of some label $k$ of $\Gamma$.
Then by applying C-I-M2 moves, C-I-R2 moves, 
and C-I-R3 moves in $V$, 
there exists 
a $(D,\alpha)$-arc free chart $\Gamma'$ 
obtained from the chart $\Gamma$ 
keeping $\alpha$ fixed 
$($see Fig.~\ref{fig19}$)$.
\end{lemma}

%%%%%%%%%%%%%%%%%%
%%%%%%%%%%%%%%%%%% Figure
%%%%%%%%%%%%%%%%%%

\begin{lemma}
\label{EdgeLabelKEdgeK+1}
Let $\Gamma$ be a chart,
and $k$ a label of $\Gamma$.
Let $D$ be a special $s$-angled disk of $\Gamma_k$
with $w(\Gamma\cap{\rm Int}D)=0$.
Let $e$ be an internal edges of label $k+\delta$
such that $e$ is a proper arc of $D$
for some $\delta\in\{+1,-1\}$.
Let $V$ be a regular neighborhood of $e$.
Then $\Gamma$ is modified to a chart by C-moves in $V$
keeping $e$ fixed without increasing 
the complexity of the chart
such that the edge $e$ intersects at most one point
for each internal edge of label $k-\delta$.
\end{lemma}

\begin{Proof}
Since $w(\Gamma\cap{\rm Int}D)=0$,
any terminal edge in $D$
intersects the boundary $\partial D$.
Let $X$ be the union of all the terminal edges in $D$.
Let $N$ be a regular neighborhood of $X$
in $D$.
Then $Cl(D-N)$ is a disk without black vertices. 
Since the edge $e$ is a proper arc of $D$,
the edge $e$ divides the disk $Cl(D-N)$ into two disks.
Let $D_1,D_2$ be the two disks.
By New Disk Lemma(Lemma~\ref{NewDiskLemma}),
we can assume that $\Gamma$ is  $(D_1,\widetilde{e})$-arc free,
where $\widetilde{e}=e\cap D_1$.

We shall show that 
for any internal edge $e'$ of label $k-\delta$,
the intersection $e\cap e'$ consists of at most two points.
If there exists an internal edge $e'$ of
label $k-\delta$
such that $e\cap e'$ 
contains at least three points,
then there exists a proper arc $\gamma$ of $D_1$
such that $\gamma\subset e'$ 
and $\partial \gamma\subset e$
(because $e'$ is a proper edge of $D$).
Thus $\gamma$ is a $(D_1,\widetilde{e})$-arc 
of label $k-\delta$.
This contradicts the fact that
$\Gamma$ is  $(D_1,\widetilde{e})$-arc free.
Hence for each internal edge $e'$ of label $k-\delta$,
the intersection $e\cap e'$ consists of
at most two points.

Hence
if there exists an internal edge $e'$ of label $k-\delta$
such that $e\cap e'$ consists of two points,
then $e'\cap D_2$ is a
$(D_2,\widetilde{e})$-arc of label $k-\delta$
(because $\Gamma$ is $(D_1,\widetilde{e})$-arc free).
Again by applying New Disk Lemma(Lemma~\ref{NewDiskLemma})
for the disk $D_2$,
we can assume that $\Gamma$ is  $(D_2,\widetilde{e})$-arc free.
Thus we can show that
for each internal edge $e'$ of label $k-\delta$,
the intersection $e\cap e'$ consists of
at most one point.
\end{Proof}

%%%%%%%%%%%%%%%%%%%%
%%%%%%%%%%%%%%%%%%%%
%%%%%%%%%%%%%%%%%%%%

The following lemma is shown in \cite[Lemma 18.24 (E)]{BraidThree}.
Thus we omit the proof.

\begin{figure}
\centerline{\includegraphics{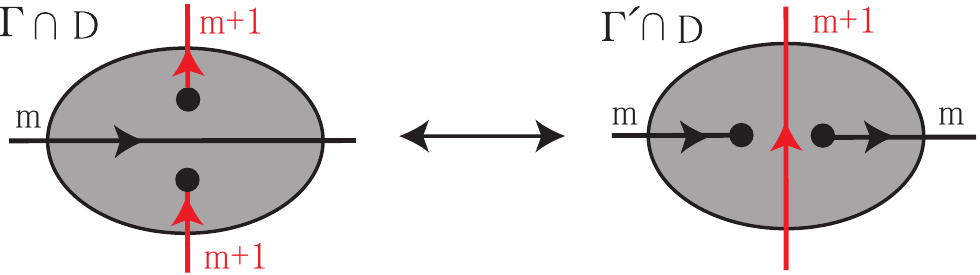}}
\caption{\label{fig20}
The gray region is the disk $D$.}
\end{figure}

\begin{lemma}{\rm $($Cut Edge Lemma$)$}
\label{CutEdgeLemma}
Let $\Gamma$ and $\Gamma'$ be charts, and $D$ a disk with $\Gamma \cap D^c=\Gamma' \cap D^c$. 
If $\Gamma \cap D$ and $\Gamma' \cap D$ are pseudo charts as shown in Fig.~\ref{fig20}, 
then $\Gamma$ is C-move equivalent to $\Gamma'$. 
\end{lemma}

%%%%%%%%%%%%%%%%%%
%%%%%%%%%%%%%%%%%% Figure
%%%%%%%%%%%%%%%%%%

%%%%%%%%%%%%%%%%%%%%
%%%%%%%%%%%%%%%%%%%%
%%%%%%%%%%%%%%%%%%%%

\begin{figure}
\centerline{\includegraphics{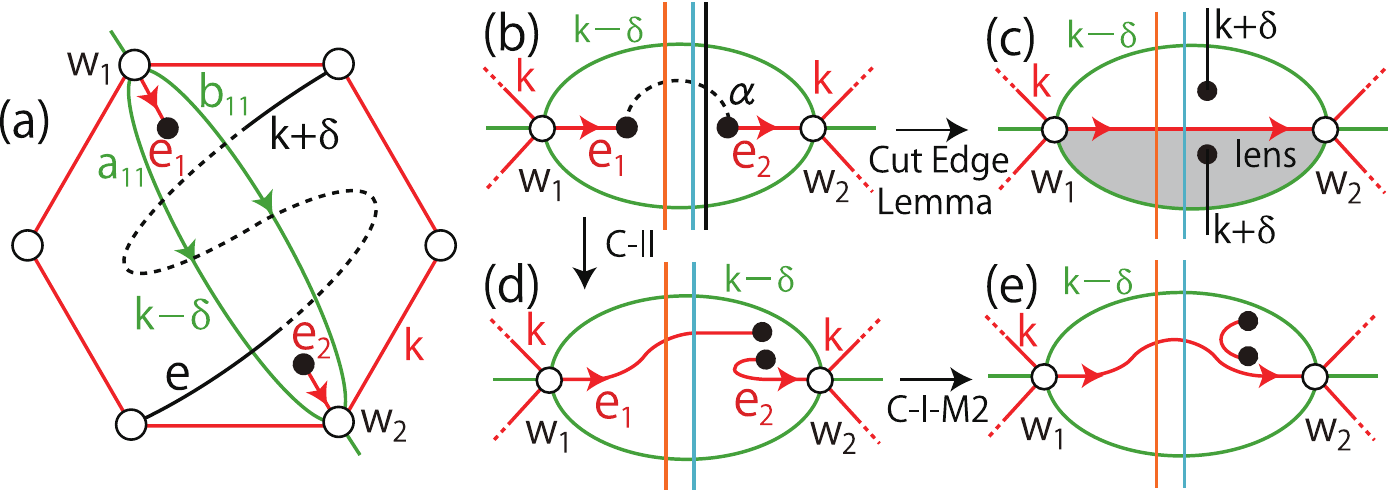}}
\caption{\label{fig21}
(a) A $k$-angled disk of $\Gamma_k$ containing a 2-angled disk.
(b),(c),(d),(e) 2-angled disks without feelers.
The gray region is a lens, $k$ is a label, $\delta\in\{+1,-1\}$.}
\end{figure}

\begin{lemma}
\label{NotMinimal2AngledDisk}
Let $\Gamma$ be a chart,
and $k$ a label of $\Gamma$.
Let $D$ be a special $s$-angled disk of $\Gamma_k$
with $w(\Gamma\cap{\rm Int}D)=0$.
Let $e_1,e_2$ be two feelers of $D$, and
$w_1,w_2$ the white vertices in $e_1,e_2$,
respectively.
Suppose that there exist 
two internal edges $a_{11},b_{11}$ of label $k-\delta$ 
in $D$ with $a_{11}\cap b_{11}=\{w_1, w_2\}$ 
for some $\delta\in\{+1,-1\}$
$($see Fig.~\ref{fig21}$($a$))$.
Let $D'$ be the $2$-angled disk of $\Gamma_{k-\delta}$ in $D$ with $\partial D'=a_{11}\cup b_{11}$.
If $D'$ intersects at most one internal edge of label $k+\delta$,
then $\Gamma$ is not minimal.
\end{lemma}

\begin{Proof}
Suppose that  $\Gamma$ is minimal.
Without loss of generality
we may assume that the feeler $e_1$ is oriented outward at $w_1$.
Since the feeler $e_1$ is middle at $w_1$ by
Assumption~\ref{AssumeTerminal},
both of $a_{11}$ and $b_{11}$ are oriented 
from $w_1$ to $w_2$.
Hence the feeler $e_2$ is oriented inward at $w_2$
(see Fig.~\ref{fig21}(a)).

Since $D'$ intersects at most one internal edge $e$
of label $k+\delta$,
by Lemma~\ref{EdgeLabelKEdgeK+1}
we can assume that each of $e\cap a_{11}$ and 
$e\cap b_{11}$ consists of
at most one point.
That is, the intersection $\Gamma_{k+\delta}\cap D'$ 
consists at most
one proper arc of $D'$.
 
Suppose that $\Gamma_{k+\delta}\cap D'$ is one proper arc of $D'$.
Let $\alpha$ be an arc in ${\rm Int}D'$ connecting
the two black vertices of $e_1$ and $e_2$
such that $\Gamma_{k+\delta}\cap\alpha$ consists of one point
(see Fig.~\ref{fig21}(b)).
Then we apply C-II moves to the two feelers $e_1$ and $e_2$
along the arc $\alpha$ such that
we can assume that $\Gamma\cap{\rm Int}\alpha=\Gamma_{k+\delta}\cap{\rm Int}\alpha=$ one point.
And then we can apply 
Cut Edge Lemma(Lemma~\ref{CutEdgeLemma})
 between $e_1$ and $e_2$
(see Fig.~\ref{fig21}(c)).
Then we obtain a lens $E$ in $D$ with $E\supset a_{11}$.
Thus by Lemma~\ref{LensThreeWhiteVertex}
we have $w(\Gamma\cap{\rm Int}E)\ge3$.
Hence $w(\Gamma\cap{\rm Int}D)\ge3$.
This contradicts $w(\Gamma\cap{\rm Int}D)=0$.
Thus $\Gamma_{k+\delta}\cap D'=\emptyset$.

Let $\alpha$ be an arc in ${\rm Int}D'$ connecting
the two black vertices of $e_1$ and $e_2$
with $\Gamma_{k+\delta}\cap\alpha=\emptyset$
(see Fig.~\ref{fig21}(b)).
By using C-II moves along $\alpha$,
we can assume that $\Gamma\cap {\rm Int}\alpha=\emptyset$
(see Fig.~\ref{fig21}(d)).
Then by applying a C-I-M2 move between $e_1$ and $e_2$,
we obtain a free edge of label $k$
(see Fig.~\ref{fig21}(e)).
Hence the complexity of the chart decreases.
Therefore $\Gamma$ is not minimal.
\end{Proof}

%%%%%%%%%%%%%%%%%%
%%%%%%%%%%%%%%%%%% Figure
%%%%%%%%%%%%%%%%%%

%%%%%%%%%%%%%%%%%%%%%%%%%%
%%%%%%%%%%%%%%%%%%%%%%%%%%
%%%%%%%%%%%%%%%%%%%%%%%%%%

%\newpage
\section{ $4$-angled disks of $\Gamma_k$
with a proper edge of label $k+\delta$}

\label{s:ProperEdge}

In this section,
we investigate a $4$-angled disk of $\Gamma_k$
with a proper edge of label $k+\delta$
for some $\delta\in\{+1,-1\}$.

Let $\Gamma$ be a chart,
and $k$ a label of $\Gamma$.
Let $D$ be a $s$-angled disk of $\Gamma_k$,
and $\delta\in\{+1,-1\}$.
Let $e$ be an internal edge of label $k+\delta$
which is a proper edge of $D$.
Let $w_1,w_2$ be the white vertices in $e$.
If any edge of label $k$ containing $w_1$ or $w_2$
does not intersect the open disk ${\rm Int}D$,
then we say that the edge $e$
is a {\it nice} edge with respect to $D$.

\begin{lemma}
\label{NiceEdge4AngledDisk}
Let $\Gamma$ be a minimal chart,
and $k$ a label of $\Gamma$.
Let $D$ be a special $4$-angled disk of $\Gamma_k$
with a nice edge.
If $w(\Gamma\cap{\rm Int}D)=0$,
then $D$ has no feelers.
\end{lemma}

\begin{Proof}
Let $w_1,w_2,w_3,w_4$ be the four white vertices situated on $\partial D$ in this order.
Let $e$ be a nice edge with respect to $D$.

Suppose that $D$ has a feeler $e_1$ with $e_1\ni w_1$.
Since $D$ has a nice edge $e$,
there are two cases:
(i) $D$ has one feeler,
(ii) $D$ has two feelers.

{\bf Case (i).}
If $D$ has exactly one feeler $e_1$,
then there are two internal edges $e',e''$ 
(possibly terminal edges) of label $k+\delta$ at $w_1$
in $D$ but not middle at $w_1$
for some $\delta\in\{+1,-1\}$.

If necessary we reflect the chart $\Gamma$,
we may assume that $\partial e=\{w_2,w_4\}$ or
$\partial e=\{w_3,w_4\}$.
If $\partial e=\{w_2,w_4\}$ and if $w_3\in\Gamma_{k+\delta}$,
then we have $w(\Gamma\cap{\rm Int}D)\ge1$
by considering as $F=D$, $m=k+\delta$ and $\varepsilon=-\delta$ in the example of IO-Calculation in Section~\ref{s:IOC}.
This contradicts $w(\Gamma\cap{\rm Int}D)=0$.
If $\partial e=\{w_2,w_4\}$ and if $w_3\in\Gamma_{k-\delta}$,
then we have $w(\Gamma\cap{\rm Int}D)\ge1$
by IO-Calculation with respect to $\Gamma_{k+\delta}$
in $D$.
This contradicts $w(\Gamma\cap{\rm Int}D)=0$.

Similarly if  $\partial e=\{w_3,w_4\}$,
then we have the same contradiction
by IO-Calculation with respect to $\Gamma_{k+\delta}$
in $D$.
Therefore Case (i) does not occur.

{\bf Case (ii).}
If $D$ has exactly two feelers,
then one of $w_2,w_3,w_4$ is contained in a feeler $e'$.
Hence there are three cases:
($w_2\in e'$, $\partial e=\{w_3,w_4\}$) or
($w_3\in e'$, $\partial e=\{w_2,w_4\}$) or
($w_4\in e'$, $\partial e=\{w_2,w_3\}$).

Let $a_{11},b_{11}$ be internal edges 
(possibly terminal edges) of label $k+\delta$ 
at $w_1$ in $D$ for some $\delta\in\{+1,-1\}$.
Since the feeler $e_1$ is middle at $w_1$ by
Assumption~\ref{AssumeTerminal},
neither $a_{11}$ nor $b_{11}$ is middle at $w_1$.
Thus by Assumption~\ref{AssumeTerminal},
neither $a_{11}$ nor $b_{11}$ is a terminal edge.
Hence the condition $w(\Gamma\cap{\rm Int}D)=0$
implies $a_{11}\cap b_{11}\ni w_2$ or 
 $a_{11}\cap b_{11}\ni w_3$ or
 $a_{11}\cap b_{11}\ni w_4$.
Thus by Lemma~\ref{NotMinimal2AngledDisk}
the chart $\Gamma$ is not minimal.
This contradicts the fact that $\Gamma$ is minimal.
Therefore Case (ii) does not occur.
Thus $D$ has no feelers.
\end{Proof}

\begin{figure}
\centerline{\includegraphics{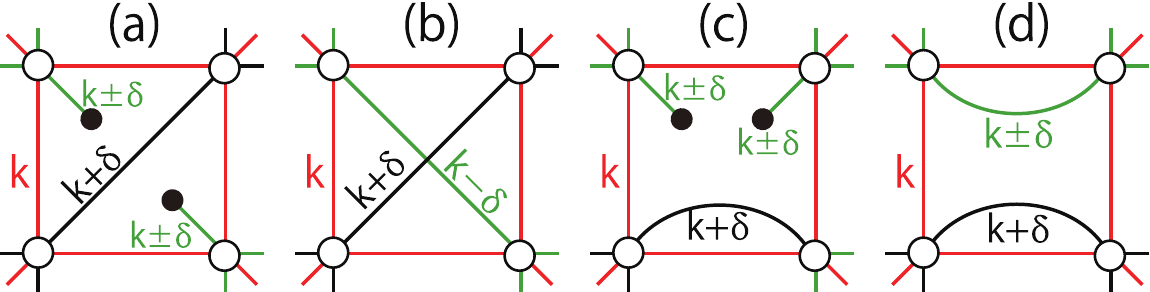}}
\caption{\label{fig22}
 4-angeld disks with nice edges where $k$ is a label,
$\delta\in\{+1,-1\}$.}
\end{figure}

\begin{lemma}
\label{NiceEdge4AngledDiskFigure}
Let $\Gamma$ be a minimal chart,
and $k$ a label of $\Gamma$.
Let $D$ be a special $4$-angled disk of $\Gamma_k$
with a nice edge of label $k+\delta$ 
for some $\delta\in\{+1,-1\}$.
If $w(\Gamma\cap{\rm Int}D)=0$,
then a regular neighborhood of $D$ contains one of 
the four pseudo charts
as shown in Fig.~\ref{fig22}
by C-moves in $D$ keeping $\Gamma_k\cup\Gamma_{k+\delta}$
fixed without increasing the complexity of $\Gamma$.
\end{lemma}

%%%%%%%%%%%%%%%%%%
%%%%%%%%%%%%%%%%%% Figure
%%%%%%%%%%%%%%%%%%

\begin{Proof}
Let $e$ be a nice edge of label $k+\delta$.
By Lemma~\ref{NiceEdge4AngledDisk},
the disk $D$ has no feelers.
Let $w_1,w_2,w_3,w_4$ be the four white vertices situated on
$\partial D$ in this order.
If necessary we reflect the chart $\Gamma$,
we can assume that 
(i) $\partial e=\{w_2,w_4\}$ or 
(ii) $\partial e=\{w_3,w_4\}$.

{\bf Case (i).}
If both of $w_1$ and $w_3$ are contained in terminal edges in $D$,
then a regular neighborhood of $D$ contains 
the pseudo chart
as shown in Fig.~\ref{fig22}(a).

If one of $w_1$ and $w_3$ is not contained in a terminal edge in $D$,
then the condition $w(\Gamma\cap{\rm Int}D)=0$ implies that
there exists an internal edge $e'$ 
connecting $w_1$ and $w_3$.
Since the edge $e'$ intersects the nice edge $e$ of 
label $k+\delta$,
the label of the edge $e'$ must be $k-\delta$.
Hence by Lemma~\ref{EdgeLabelKEdgeK+1},
we can assume that $e\cap e'$ is one point.
Thus a regular neighborhood of $D$ contains 
the pseudo chart
as shown in Fig.~\ref{fig22}(b).

{\bf Case (ii).}
If both of $w_1$ and $w_2$ are contained in terminal edges in $D$,
then a regular neighborhood of $D$ contains 
the pseudo chart
as shown in Fig.~\ref{fig22}(c).

If one of $w_1$ and $w_2$ is not contained in a terminal edge in $D$,
then there exists an internal edge $e'$ 
connecting $w_1$ and $w_2$.
If the label of the edge $e'$ is $k+\delta$,
then  $e\cap e'=\emptyset$.
Thus a regular neighborhood of $D$ contains 
the pseudo chart
as shown in Fig.~\ref{fig22}(d).

If the label of the edge $e'$ is $k-\delta$,
then by Lemma~\ref{EdgeLabelKEdgeK+1}
we can assume that $e\cap e'=\emptyset$.
Thus a regular neighborhood of $D$ contains 
the pseudo chart
as shown in Fig.~\ref{fig22}(d).
\end{Proof}

%%%%%%%%%%%%%%%%%%%%%%%%
%%%%%%%%%%%%%%%%%%%%%%%%
%%%%%%%%%%%%%%%%%%%%%%%%

\begin{figure}
\centerline{\includegraphics{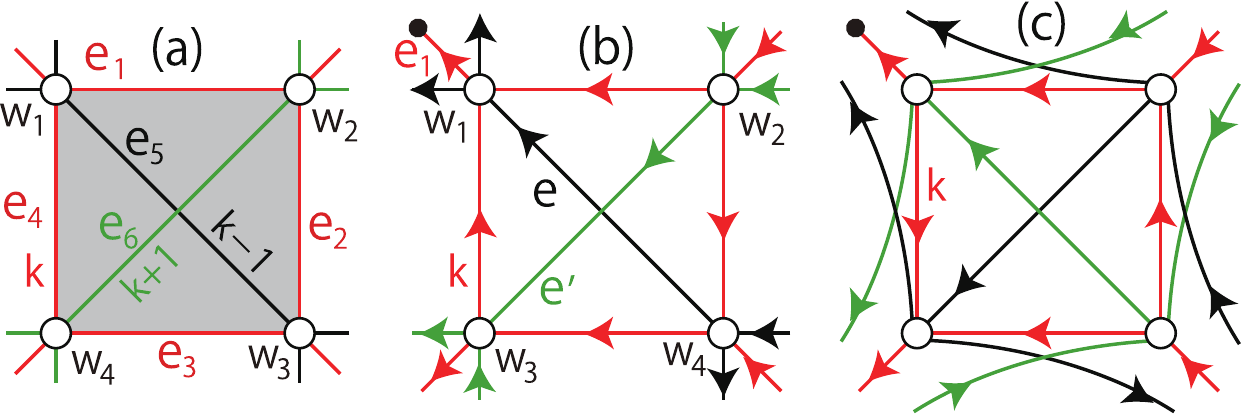}}
\caption{\label{fig23}
(a) The gray region is the M4-disk.
(b),(c) 4-angled disks with nice edges.}
\end{figure}

Let $\Gamma$ be a chart and 
$k$ a label of $\Gamma$.
If a disk $D$ satisfies the following three conditions,
then $D$ is called an 
{\it M4-disk of label $k$}
 (see Fig.~\ref{fig23}(a)).
\begin{enumerate}
\item[(i)] 
$\partial D$ consists of four internal edges 
$e_1,e_2,e_3,e_4$ of label $k$ 
situated on $\partial D$ 
in this order.
\item[(ii)] 
Set $w_1=e_1\cap e_4,w_2=e_1\cap e_2,
w_3=e_2\cap e_3,w_4=e_3\cap e_4$.
Then
\begin{enumerate}
\item[(a)] 
$D\cap\Gamma_{k-1}$ consists of an internal edge $e_5$ 
connecting $w_1$ and $w_3$, 
and
\item[(b)] 
$D\cap\Gamma_{k+1}$ consists of an internal edge $e_6$ 
connecting $w_2$ and $w_4$. 
\end{enumerate}
\item[(iii)] 
Int$D$ does not contain any white vertex.
\end{enumerate}
We call the union 
$X=\cup_{i=1}^6 e_i$ 
the {\it M4-pseudo chart} for the disk $D$.

%%%%%%%%%%%%%%%%%%
%%%%%%%%%%%%%%%%%% Figure
%%%%%%%%%%%%%%%%%%

\begin{lemma}
\label{M4-diskLemma}
{\rm (\cite[Lemma 7.3]{Chart33})}
Let $\Gamma$ be a chart, and 
$k$ a label of $\Gamma$. 
Suppose that $D$ is an M4-disk 
of label $k$ 
with an M4-pseudo chart $X$. 
Then by deforming $\Gamma$ 
in a regular neighbourhood of $D$ 
without increasing the complexity of $\Gamma$, 
the chart $\Gamma$ is 
C-move equivalent to a chart $\Gamma'$ 
with
$D\cap(\cup_{i=k-2}^{k+2}\Gamma'_i)=X$.
\end{lemma}

%%%%%%%%%%%%%%%%%
%%%%%%%%%%%%%%%%%

\begin{lemma}
\label{4AngledDiskCrossNotMinimal}
Let $\Gamma$ be a chart, and $k$ a label of $\Gamma$.
Let $D$ be a $4$-angled disk of $\Gamma_k$ 
without feelers with $w(\Gamma\cap{\rm Int}D)=0$.
Let $w_1,w_2,w_3,w_4$ be the four white vertices situated on
$\partial D$ in this order.
Suppose that there exist two internal edges $e,e'$ in $D$
with $\partial e=\{w_1,w_3\}$ and 
$\partial e'=\{w_2,w_4\}$.
If $\Gamma$ satisfies one of the following two conditions,
then $\Gamma$ is not minimal.
\begin{enumerate}
\item[{\rm (a)}] 
One of  $w_1$ and $w_3$ is contained in a terminal edge of 
label $k$, and the edge $e'$ is middle at $w_2$ or $w_4$.
\item[{\rm (b)}] 
One of  $w_2$ and $w_4$ is contained in a terminal edge of 
label $k$, and the edge $e$ is middle at $w_1$ or $w_3$.
\end{enumerate}
\end{lemma}

\begin{Proof}
We shall show only that
if $\Gamma$ satisfies Condition~(a),
then $\Gamma$ is not minimal.

Suppose that $\Gamma$ satisfies  Condition~(a)
and $\Gamma$ is minimal.
By Lemma~\ref{NiceEdge4AngledDiskFigure},
we can assume that
a regular neighborhood of $D$ contains the pseudo chart
as shown in Fig.~\ref{fig22}(b).
That is the intersection $e\cap e'$ is one crossing.

Without loss of generality we can assume that
the white vertex $w_1$ is contained in 
a terminal edge $e_1$ of label $k$, and
\begin{enumerate}
\item[(1)] 
the edge $e'$ is middle at $w_2$. 
\end{enumerate}
Moreover we can assume that
the edge $e_1$ is oriented outward at $w_1$.
Since $e_1$ is middle at $w_1$ by Assumption~\ref{AssumeTerminal},
we have that
the edge $e$ is oriented inward at $w_1$, and
the two internal edges of label $k$ at $w_1$
are oriented inward at $w_1$.
Thus by (1)
the edge $e'$ is oriented outward at $w_2$,
and the two internal edges of label $k$ at $w_2$ in
$\partial D$ are oriented outward at $w_2$.
Therefore we have the orientations of the other edges
as shown in Fig.~\ref{fig23}(b).

Now the set $X=e\cup e'\cup\partial D$
is an M4-pseudo chart for the disk $D$.
Thus by Lemma~\ref{M4-diskLemma},
we can assume that $D\cap(\cup_{i=k-2}^{k+2}\Gamma_i)=X$.
By applying a C-I-M4 move
in a neighborhood of $D$,
we obtain a terminal edge of label $k$ at $w_1$
but not middle at $w_1$
(see Fig.~\ref{fig23}(c)).
This contradicts Assumption~\ref{AssumeTerminal}.
Hence if $\Gamma$ satisfies  Condition~(a),
then $\Gamma$ is not minimal.
\end{Proof}

%%%%%%%%%%%%%%%%%%%%%%%%%%
%%%%%%%%%%%%%%%%%%%%%%%%%%
%%%%%%%%%%%%%%%%%%%%%%%%%%

%\newpage
\section{ $5$-angled disks of $\Gamma_k$
with a proper edge of label $k+\delta$}

\label{s:ProperEdge5AngledDisk}

In this section we investigate a $5$-angled disk of $\Gamma_k$ with a proper edge of label $k+\delta$
for some $\delta\in\{+1,-1\}$.

\begin{lemma}
\label{NiceEdge5AngledDiskAtMostOneFeeler}
Let $\Gamma$ be a minimal chart,
and $k$ a label of $\Gamma$.
Let $D$ be a special $5$-angled disk of $\Gamma_k$
with a nice edge. 
If $w(\Gamma\cap{\rm Int}D)=0$,
then $D$ has at most one feeler.
\end{lemma}

\begin{figure}
\centerline{\includegraphics{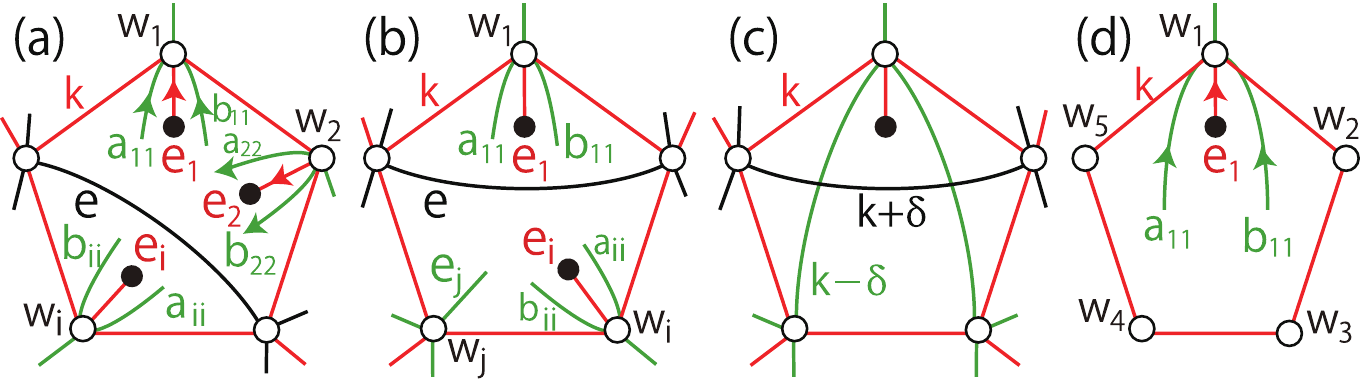}}
\caption{\label{fig24}
 5-angeld disks with a nice edge where $k$ is a label, $\delta\in\{+1,-1\}$.}
\end{figure}

\begin{Proof}
Let $w_1,w_2,\cdots,w_5$ be the five white vertices
situated on $\partial D$ in this order.
Suppose that $D$ has at least two feelers.
Since $D$ has a nice edge,
the disk $D$ has at most three feelers.
Thus there are two cases:
(i) $D$ has three feelers (see Fig.~\ref{fig24}(a)),
(ii) $D$ has two feelers (see Fig.~\ref{fig24}(b)).

{\bf Case (i).}
Since $D$ has three feelers,
we can assume that $w_1,w_2$ are contained in feelers
$e_1,e_2$, respectively.
Moreover there exists $i\in\{3,4,5\}$
such that the white vertex $w_i$ is contained in
the third feeler $e_i$ (see Fig.~\ref{fig24}(a)).

Let $a_{11},b_{11},a_{22},b_{22},a_{ii},b_{ii}$
be internal edges (possibly terminal edges)
of label $k-1$ or $k+1$ 
at $w_1,w_1,w_2,w_2,w_i,w_i$ in $D$, respectively.
Since the three edges $e_1,e_2,e_i$ are middle 
at $w_1,w_2,w_i$, respectively 
by Assumption~\ref{AssumeTerminal},
\begin{enumerate}
\item[(1)] none of $a_{11},b_{11},a_{22},b_{22},a_{ii},b_{ii}$ are terminal edges.
\end{enumerate}

If necessary we change the orientations of all the edges,
we can assume that the edge $e_1$ is oriented inward at $w_1$.
Since $e_1$ is middle at $w_1$ by Assumption~\ref{AssumeTerminal}, we can show that
\begin{enumerate}
\item[(2)] both of $a_{11},b_{11}$ are oriented inward at $w_1$, and
\item[(3)] both of $a_{22},b_{22}$ are oriented outward at $w_2$.
\end{enumerate}

If $e_i$ is oriented inward at $w_i$,
then we can show that 
both of $a_{ii},b_{ii}$ are oriented inward at $w_i$.
Thus by (1),(2) and (3)
we have $w(\Gamma\cap{\rm Int}D)\ge1$
by IO-Calculation with respect to $\Gamma_{k-1}$ or $\Gamma_{k+1}$ in $D$.
This contradicts $w(\Gamma\cap{\rm Int}D)=0$.

If $e_i$ is oriented outward at $w_i$,
then we can show that 
both of $a_{ii},b_{ii}$ are oriented outward at $w_i$.
Thus we have the same contradiction.
Therefore Case (i) does not occur.

{\bf Case (ii).}
Without loss of generality
we can assume that 
$w_1$ is contained in a feeler $e_1$.
Since $D$ has exactly two feelers,
there exists $i\in\{2,3,4,5\}$
such that the white vertex $w_i$ is contained in a feeler $e_i$.
Moreover
since $D$ has a nice edge,
there exists $j\in\{2,3,4,5\}$
with $j\not=i$ such that 
the white vertex $w_j$ is not contained in the nice edge.
Thus there exists an internal edge $e_j$
(possibly a terminal edge) of label $k+\delta$ at $w_j$ 
in $D$
for some $\delta\in\{+1,-1\}$
(see Fig.~\ref{fig24}(b)).

Let $a_{11},b_{11},a_{ii},b_{ii}$
be internal edges (possibly terminal edges)
of label $k-1$ or $k+1$ 
at $w_1,w_1,w_i,w_i$ in $D$, respectively.
Similarly we can show that 
none of $a_{11},b_{11},a_{ii},b_{ii}$ are terminal edges.
Moreover we can assume that
both of $a_{11}$ and $b_{11}$ are oriented inward at $w_1$.
Thus the condition $w(\Gamma\cap{\rm Int}D)=0$ implies that
we have $a_{11}\cap b_{11}\ni w_i$ and
the edge $e_j$ is a terminal edge.
Hence by Lemma~\ref{NotMinimal2AngledDisk}
the chart $\Gamma$ is not minimal.
This contradicts the fact that 
$\Gamma$ is minimal.
Therefore Case (ii) does not occur.
Hence $D$ has at most one feeler.
\end{Proof}

%%%%%%%%%%%%%%%%%%
%%%%%%%%%%%%%%%%%% Figure
%%%%%%%%%%%%%%%%%%

%%%%%%%%%%%%%%%%%
%%%%%%%%%%%%%%%%%
%%%%%%%%%%%%%%%%%

\begin{lemma}
\label{NiceEdge5AngledDiskOneFeeler}
Let $\Gamma$ be a minimal chart,
and $k$ a label of $\Gamma$.
Let $D$ be a special $5$-angled disk of $\Gamma_k$
with a nice edge of label $k+\delta$ 
for some $\delta\in\{+1,-1\}$.
If $w(\Gamma\cap{\rm Int}D)=0$
and if $D$ has exactly one feeler, then
a regular neighborhood of $D$ contains  
 the pseudo chart
as shown in Fig.~\ref{fig24}$($c$)$
by C-moves in $D$ keeping $\Gamma_k\cup\Gamma_{k+\delta}$
fixed without increasing the complexity of $\Gamma$.
\end{lemma}

\begin{Proof}
Let $w_1,w_2,\cdots,w_5$ be the five white vertices
situated on $\partial D$  in this order.
Without loss of generality we can assume that
the white vertex $w_1$ is contained in a feeler $e_1$.
Let $a_{11},b_{11}$
be internal edges (possibly terminal edges)
of label $k+\varepsilon$ 
at $w_1$ in $D$ for some $\varepsilon\in\{+1,-1\}$
such that $a_{11},e_1,b_{11}$ lie anticlockwise around
the white vertex $w_1$
(see Fig.~\ref{fig24}(d)).

If necessary we change the orientations
of all the edges,
we can assume that 
$e_1$ is oriented inward at $w_1$.
Thus we can show that
both of $a_{11}$ and $b_{11}$ are oriented inward at $w_1$,
and neither $a_{11}$ nor $b_{11}$ is a terminal edge.

If $a_{11}\ni w_5$ (resp. $b_{11}\ni w_2$),
then there exists a lens $E$ in $D$
with $\partial E\supset a_{11}$ 
(resp. $\partial E\supset b_{11}$).
Hence by Lemma~\ref{LensThreeWhiteVertex},
we have $w(\Gamma\cap{\rm Int}E)\ge3$.
Thus $w(\Gamma\cap{\rm Int}D)\ge3$.
This contradicts $w(\Gamma\cap{\rm Int}D)=0$.
Hence $a_{11}\not\ni w_5$ and $b_{11}\not\ni w_2$.
Thus $a_{11}\ni w_4$, $b_{11}\ni w_3$.

Since $D$ has a nice edge $e$ of label $k+\delta$,
we have 
$\partial e=\{w_2,w_5\}$.
Thus the nice edge $e$ of label $k+\delta$ intersects
the edges $a_{11}$ and $b_{11}$ of label $k+\varepsilon$.
Hence $\varepsilon=-\delta$,
i.e. the label of $a_{11}$ and $b_{11}$ must be $k-\delta$.
Thus by Lemma~\ref{EdgeLabelKEdgeK+1}
we can assume that 
each of $e\cap a_{11}$ and $e\cap b_{11}$
is one point.
Therefore a regular neighborhood of $D$ contains
 the pseudo chart
as shown in Fig.~\ref{fig24}(c).
\end{Proof}

%%%%%%%%%%%%%%%%%%%%%%
%%%%%%%%%%%%%%%%%%%%%%
%%%%%%%%%%%%%%%%%%%%%%

%\newpage
\section{ $7$-angled disks of $\Gamma_k$
with a proper edge of label $k+\delta$}

\label{s:ProperEdge7AngledDisk}

In this section we investigate a $7$-angled disk of $\Gamma_k$ with a proper edge of label $k+\delta$
for some $\delta\in\{+1,-1\}$.

\begin{figure}
\centerline{\includegraphics{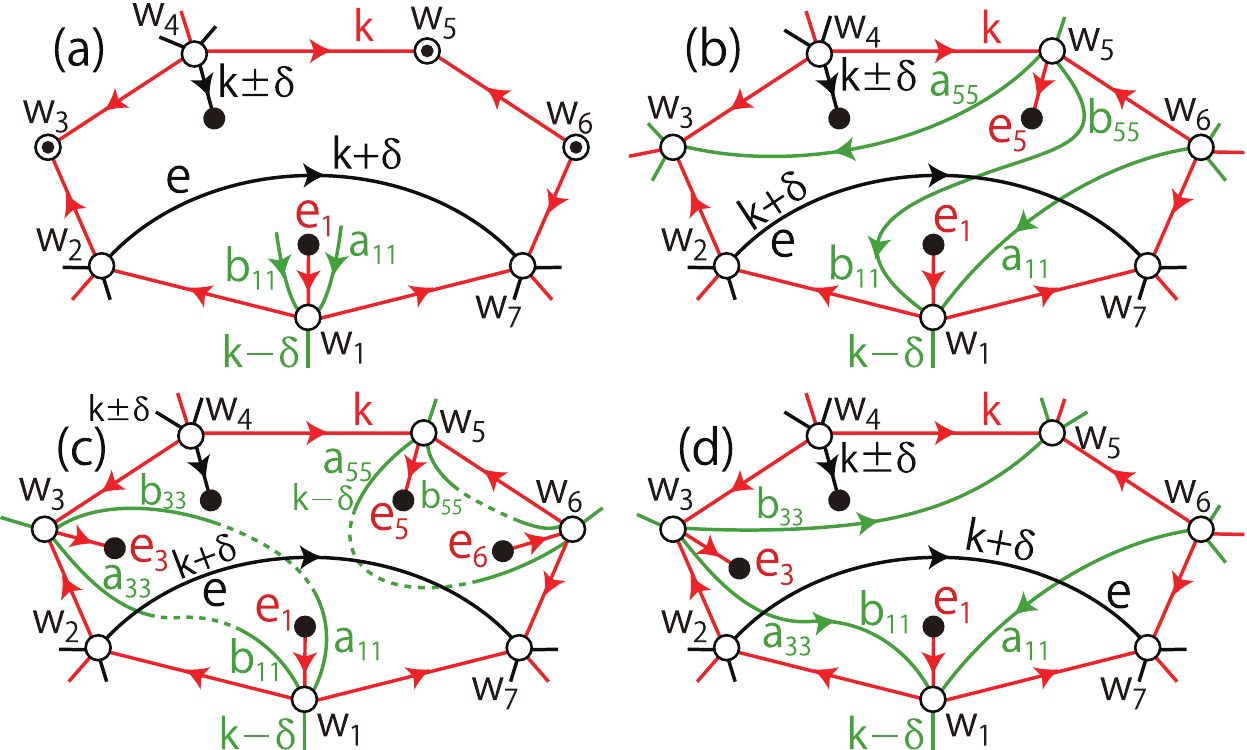}}
\caption{\label{fig25}
 7-angeld disks with a nice edge.}
\end{figure}

\begin{lemma}
\label{NiceEdge7AngledDisk}
Let $\Gamma$ be a minimal chart,
and $k$ a label of $\Gamma$.
Let $D$ be a special $7$-angled disk of $\Gamma_k$,
and $w_1,w_2,\cdots,w_7$ the seven white vertices situated on
$\partial D$ in this order.
Suppose that 
\begin{enumerate}
\item[{\rm (a)}]
$D$ has a feeler $e_1$ at $w_1$,
\item[{\rm (b)}]
$D$ has a nice edge $e$ of label $k+\delta$
with $\partial e=\{w_2,w_7\}$ 
for some $\delta\in\{+1,-1\}$,
\item[{\rm (c)}] $w_1,w_3,w_5,w_6$ are white vertices in
$\Gamma_{k-\delta}$,
\item[{\rm (d)}]
$D$ has a terminal edge of label $k\pm\delta$ at $w_4$,
\item[{\rm (e)}]
$w_3,w_5,w_6$ are BW-vertices with respect to $\Gamma_k$
$($see Fig.~\ref{fig25}$($a$))$.
\end{enumerate}
If $w(\Gamma\cap{\rm Int}D)=0$,
then $D$ has at most two feelers.
Moreover if $D$ has two feelers, then
a regular neighborhood of $D$ contains one of 
the RO-family of the pseudo chart
as shown in Fig.~\ref{fig25}$($b$)$
by C-moves in $D$ without increasing 
the complexity of $\Gamma$.
\end{lemma}

%%%%%%%%%%%%%%%%%%
%%%%%%%%%%%%%%%%%% Figure
%%%%%%%%%%%%%%%%%%

\begin{Proof}
We shall show that the white vertex $w_4$ is not contained in
a feeler.
If $w_4$ is contained in a feeler,
then by Condition (d) of this lemma
the white vertex $w_4$ is contained in two terminal edges
of label $k$ and $k\pm\delta$.
However this contradicts Assumption~\ref{AssumeTerminal}.
Hence 
\begin{enumerate}
\item[(1)]  $w_4$ is not contained in a feeler.
\end{enumerate}

Without loss of generality,
we can assume that the feeler $e_1$ is oriented inward
at $w_1$.
Thus by Assumption~\ref{AssumeTerminal},
\begin{enumerate}
\item[(2)] the two internal edges of label $k$ at $w_1$
are oriented outward at $w_1$.
\end{enumerate}

If necessary we reflect the chart $\Gamma$,
we can assume that the nice edge $e$ is oriented 
from $w_2$ to
$w_7$.
Hence by (2), the internal edge of label $k$
connecting $w_2$ to $w_3$ is oriented from $w_2$ to $w_3$.
Thus by Lemma~\ref{OriBWvertex},
the internal edge of label $k$
connecting $w_3$ to $w_4$ is oriented from $w_4$ to $w_3$.
Since $D$ has a terminal edge $e_4$ of label $k\pm\delta$ at $w_4$
by (1) and Condition~(d) of this lemma,
the edge $e_4$ is oriented outward at $w_4$
by Assumption~\ref{AssumeTerminal}.
Moreover the internal edge of label $k$ 
connecting $w_4$ to $w_5$ is oriented from $w_4$ to $w_5$.
Therefore by Lemma~\ref{OriBWvertex},
a regular neighborhood of $D$ contains the pseudo chart
as shown in Fig.~\ref{fig25}(a).

Let $e_3,e_5,e_6$ be the terminal edges of label $k$
at $w_3,w_5,w_6$, respectively.
Let $a_{ii},b_{ii}$ $(i=1,3,5,6)$ be
internal edges (possibly terminal edges) 
of label $k-\delta$ at $w_i$
such that $a_{11},b_{11},a_{66},b_{66}$ are oriented inward 
at $w_1,w_1,w_6,w_6$, respectively, and
 $a_{33},b_{33},a_{55},b_{55}$ are oriented outward 
at $w_3,w_3,w_5,w_5$, respectively.
We can show that 
\begin{enumerate}
\item[(3)]
none of the eight edges $a_{11},b_{11},a_{33},b_{33},a_{55},b_{55},a_{66},b_{66}$ are terminal edges.
\end{enumerate}

{\bf Claim~1.} The disk $D$ has at most two feelers.

{\it Proof of Claim~$1$.} Suppose that $D$ has at least three feelers. Since $D$ contains the feeler $e_1$ at $w_1$,
there are two cases:
(i) $e_3\cup e_5\cup e_6\subset D$,
(ii) only two of $e_3,e_5,e_6$ are contained in $D$.

{\bf Case (i).}
By (3), we have $a_{11}\cap b_{11}\ni w_3$ and
$a_{66}\cap b_{66}\ni w_5$
(see Fig.~\ref{fig25}(c)).
Thus there exists a 2-angled disk $E$ of $\Gamma_{k-\delta}$
with $\partial E=a_{11}\cup b_{11}$ such that
$E$ intersects only one internal edge $e$ of 
label $k+\delta$.
Hence by Lemma~\ref{NotMinimal2AngledDisk}
the chart $\Gamma$ is not minimal.
This contradicts the fact that $\Gamma$ is minimal.
Thus Case (i) does not occur.

{\bf Case (ii).}
If $e_6\subset D$,
then there are four internal edges $a_{11},b_{11},a_{66},b_{66}$ are oriented inward at $w_1,w_1,w_6,w_6$ in $D$,
respectively.
Thus we have $w(\Gamma\cap{\rm Int}D)\ge1$
by (3) and by IO-Calculation with respect to $\Gamma_{k-\delta}$ in $D$.
This contradicts $w(\Gamma\cap{\rm Int}D)=0$.

If $e_6\not\subset D$,
then $e_3\cup e_5\subset D$.
Thus  there are four internal edges $a_{33},b_{33},a_{55},b_{55}$ are oriented outward at $w_3,w_3,w_5,w_5$ in $D$,
respectively.
Hence we have $w(\Gamma\cap{\rm Int}D)\ge1$
by (3) and by IO-Calculation with respect to $\Gamma_{k-\delta}$ in $D$.
This contradicts $w(\Gamma\cap{\rm Int}D)=0$.
Thus Case (ii) does not occur.

Therefore both of the two cases do not occur.
Hence $D$ has at most two feelers. Thus Claim~1 holds.
{\hfill {$\square$}\vspace{1.5em}}

We shall show that
if $D$ has two feelers,
then 
a regular neighborhood of $D$ contains one of 
the RO-family of the pseudo chart
as shown in Fig.~\ref{fig25}(b)
by C-moves in $D$ without increasing 
the complexity of $\Gamma$.

Suppose that $D$ contains two feelers.
Since $D$ contains the feeler $e_1$,
the disk $D$ contains only one of $e_3,e_5,e_6$.
 
{\bf Claim 2.} The disk $D$ contains the feeler $e_5$.

{\it Proof of Claim~$2$.}
Now $a_{11},e_1,b_{11}$ lie anticlockwise around the white
vertex $w_1$.

Suppose that $D$ contains the edge $e_3$.
Then $a_{33},e_3,b_{33}$ lie anticlockwise around the white
vertex $w_3$.
By (3), the edge $b_{33}$ is not a terminal edge.
Thus there are two cases:
$b_{33}\ni w_1$ or $b_{33}\ni w_5$.

If $b_{33}\ni w_1$,
then $a_{33}\ni w_1$.
Hence $a_{33}\cap b_{33}\ni w_1$.
By the similar way of the proof of Case (i) of Claim~1,
we have the same contradiction.

If $b_{33}\ni w_5$, then
$a_{33}=b_{11}$ and $a_{11}\ni w_6$
(see Fig.~\ref{fig25}(d)).
Thus by Lemma~\ref{EdgeLabelKEdgeK+1},
we can assume that
$b_{33}\cap\Gamma_{k+\delta}=\emptyset$.
By C-II moves, we can move the black vertex in $e_3$
near the white vertex $w_5$ along the edge $b_{33}$.
By applying a C-I-M2 move between $e_3$ and
the internal edge of label $k$ at $w_5$,
we obtain a new terminal edge of label $k$ at $w_6$.
However the terminal edge is not middle at $w_6$.
This contradicts Assumption~\ref{AssumeTerminal}.
Thus $D$ does not contain $e_3$.

Suppose that $D$ contains the edge $e_6$.
Then the four internal edges $a_{11},b_{11},a_{66},b_{66}$ are oriented inward at $w_1,w_1,w_6,w_6$ in $D$,
respectively.
Hence we have $w(\Gamma\cap{\rm Int}D)\ge1$
by (3) and by IO-Calculation with respect to $\Gamma_{k-\delta}$ in $D$.
This contradicts $w(\Gamma\cap{\rm Int}D)=0$.
Thus $D$ does not contain the edge $e_6$.
Therefore $D$ contains the edge $e_5$.
Thus Claim~2 holds.
{\hfill {$\square$}\vspace{1.5em}}

By Claim~2, the disk $D$ contains two feelers $e_1,e_5$.
Now $a_{11},e_1,b_{11}$ lie anticlockwise around the white
vertex $w_1$, and 
$a_{55},e_5,b_{55}$ lie anticlockwise around the white
vertex $w_5$.
By (3), the edge $a_{55}$ is not a terminal edge.
Thus there are two cases:
$a_{55}\ni w_1$ or $a_{55}\ni w_3$.

If $a_{55}\ni w_1$, then  $b_{55}\ni w_1$.
Hence $a_{55}\cap b_{55}\ni w_1$.
By the similar way of the proof of Case (i) of Claim~1,
we have the same contradiction.

If $a_{55}\ni w_3$,
then $b_{55}=b_{11}$ and $a_{11}\ni w_6$.
Thus by Lemma~\ref{EdgeLabelKEdgeK+1},
we can assume that
$a_{55}\cap\Gamma_{k+\delta}=\emptyset$,
$a_{11}\cap\Gamma_{k+\delta}=a_{11}\cap e=$one crossing, and
$b_{11}\cap\Gamma_{k+\delta}=b_{11}\cap e=$one crossing.
Therefore a regular neighborhood of $D$ contains
the pseudo chart as shown in 
Fig.~\ref{fig25}(b).
Thus we complete the proof of Lemma~\ref{NiceEdge7AngledDisk}.
\end{Proof}

The following two lemmata will be used in the next section.

\begin{figure}
\centerline{\includegraphics{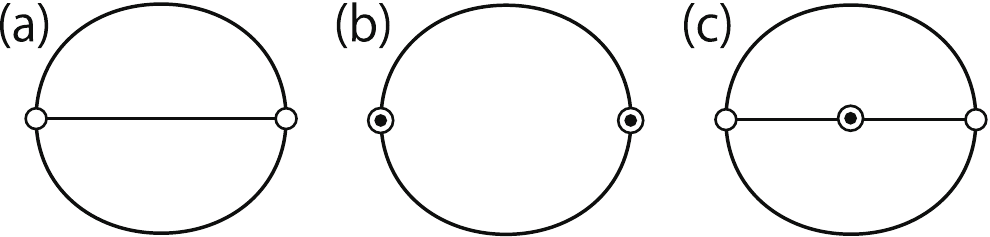}}
\caption{\label{fig26}
The graphs (a) and (b) have two white vertices,
the other has three white vertices.}
\end{figure}

\begin{lemma}
\label{LemmaWithTerminal3B}
{\rm (\cite[Lemma 3.2(2)]{ChartAppV})}
Let $\Gamma$ be a minimal chart,
and $m$ a label of $\Gamma$.
Let $G$ be a connected component of $\Gamma_m$.
If $1\le w(G)\le 3$
and $G$ does not contain any loop, 
then $G$ is one of three graphs as shown 
in Fig.~\ref{fig26}.
\end{lemma}

%%%%%%%%%%%%%%%%%%
%%%%%%%%%%%%%%%%%% Figure
%%%%%%%%%%%%%%%%%%

\begin{figure}
\centerline{\includegraphics{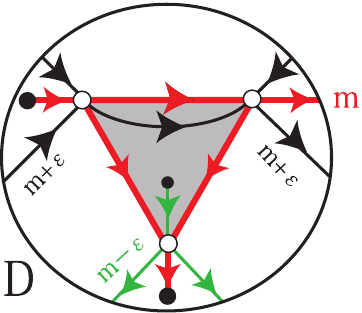}}
\caption{\label{fig27}
The gray region is the 3-angled disk $D_1$. 
The thick lines are edges of label $m$,
and $\varepsilon\in\{+1,-1\}$.}
\end{figure}

\begin{lemma}$($Triangle Lemma$)$
{\rm (\cite[Lemma 8.3(2)]{ChartAppIV})}
\label{LemmaTriangle}
 For a minimal chart $\Gamma$, 
if there exists a $3$-angled disk $D_1$ of $\Gamma_m$ without feelers in a disk $D$ as shown in 
Fig.~\ref{fig27},
then $w(\Gamma\cap${\rm Int}$D_1)\ge1$.
\end{lemma}

%%%%%%%%%%%%%%%%%%
%%%%%%%%%%%%%%%%%% Figure
%%%%%%%%%%%%%%%%%%

%\newpage
%%%%%%%%%%%%%%%%%%%%%%%%%%
%%%%%%%%%%%%%%%%%%%%%%%%%%
%%%%%%%%%%%%%%%%%%%%%%%%%%
\section{A chart of type $(4,3)$}
\label{s:Chart43}

In this section we investigate a minimal chart of type
$(m;4,3)$.

Suppose that there exists a minimal chart $\Gamma$ of type
$(m;4,3)$.
By Lemma~\ref{LemmaNoLoop},
the chart $\Gamma$ has no loop.
Since $w(\Gamma_{m+2})=3$, by Lemma~~\ref{LemmaWithTerminal3B}
the set $\Gamma_{m+2}$ contains the graph as shown 
in Fig.~\ref{fig26}(c), say $G$.
Let $e_1$ be the terminal edge of label $m+2$ in $G$,
and $w_1$ the white vertex in $e_1$.
The graph $G$ divides $S^2$ into three disks.
One of the three disks is a 2-angled disk, say $D_1$.
One of the three disks contains the edge $e_1$,
say $D_2$.
Let $D_3$ be the third disk.

Since $D_2$ is a 3-angled disk with one feeler $e_1$,
by Lemma~\ref{Theorem3AngledDisk}(cf. Lemma~\ref{ROfamily3AngledDiskTypeC})
we have 
\begin{enumerate}
\item[(a)] $w(\Gamma\cap{\rm Int}D_2)\ge1$.
\end{enumerate}
Without loss of generality we can assume that
the edge $e_1$ is oriented outward at $w_1$.
Since $e_1$ is middle at $w_1$ 
by Assumption~\ref{AssumeTerminal},
the two internal edges of label $m+2$ at $w_1$
are oriented inward at $w_1$.
Hence the 2-angled disk $D_1$
satisfies the condition of Lemma~\ref{Lemma2AngledDisks}.
Thus by Lemma~\ref{Lemma2AngledDisks}(a),
we have 
\begin{enumerate}
\item[(b)] $w(\Gamma\cap{\rm Int}D_1)\ge1$.
\end{enumerate}

\begin{figure}
\centerline{\includegraphics{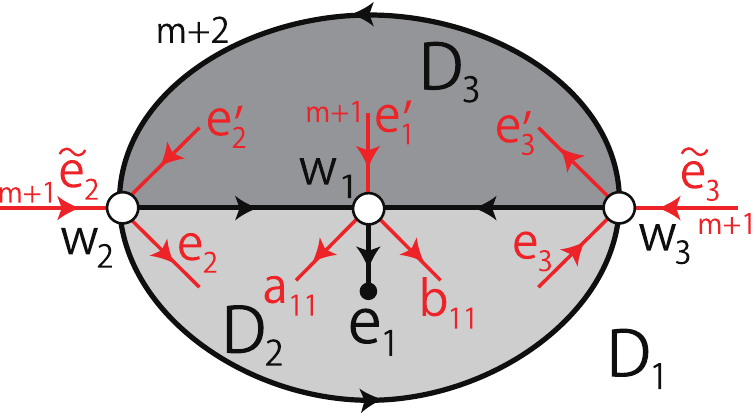}}
\caption{The light gray region is the 3-angled disk $D_2$
with a feeler $e_1$. The dark gray region is the 3-angled disk $D_3$ without feelers.
\label{fig28}
}
\end{figure}

Let $w_2,w_3$ be the two white vertices in $\Gamma_{m+2}$
different from $w_1$.
The intersection $D_1\cap D_2$ is an internal edge 
of label $m+2$ connecting $w_2$ and $w_3$.
Without loss of generality
we can assume that
the edge $D_1\cap D_2$ is oriented from $w_2$ to $w_3$.
By looking around $w_2$,
the edge $D_1\cap D_3$ is oriented from $w_3$ to $w_2$.
Therefore the chart $\Gamma$
contains the pseudo chart as shown in Fig.~\ref{fig28}.

%%%%%%%%%%%%%%%%%%
%%%%%%%%%%%%%%%%%% Figure
%%%%%%%%%%%%%%%%%%

Since $\Gamma$ is of type $(m;4,3)$,
we have 
\begin{enumerate}
\item[(c)] $w(\Gamma\cap{\rm Int}D_1)+w(\Gamma\cap{\rm Int}D_2)+w(\Gamma\cap{\rm Int}D_3)=4$.
\end{enumerate}
Moreover the four white vertices in $\Gamma-\Gamma_{m+2}$
are contained in $\Gamma_m\cap\Gamma_{m+1}$.
Thus 
we have
\begin{enumerate}
\item[(d)]
$w(\Gamma\cap{\rm Int}D_i)=w(\Gamma_{m+1}\cap{\rm Int}D_i)$ for $i=1,2,3$.
\end{enumerate}

By (d) and Lemma~\ref{Lemma2AngledDisks}(c)
we have the following lemma:
\begin{lemma}
There does not exist a minimal chart of type $(m;4,3)$
with $w(\Gamma\cap{\rm Int}D_1)=2$.
\end{lemma}

Therefore by (a),(b),(c) and the above lemma,
there are four cases (see Table~1).

\begin{center}
$\begin{array}{|c||c|c|c|c|}
\hline
w(\Gamma\cap{\rm Int}D_1) & 1 & 1 & 1 & 3\\
\hline
w(\Gamma\cap{\rm Int}D_2) & 1 & 2 & 3 & 1\\ 
\hline
w(\Gamma\cap{\rm Int}D_3) & 2 & 1 & 0 & 0 \\
\hline
& {\rm Sect.~\ref{s:D1OneD2One}}  &  {\rm Sect.~\ref{s:D1OneD2Two}} &  {\rm Sect.~\ref{s:D1OneD2Three}} &
{\rm Sect.~\ref{s:D1ThreeD2OneD3Zero}} \\
\hline
\end{array}$
\vspace{2mm}

{Table~$1$. There are four cases for a minimal chart $\Gamma$ of type $(m;4,3)$.}
\end{center}

From now on throughout this paper,
we use the notations as shown in Fig.~\ref{fig28},
where
$e_2,e_2',\widetilde{e_2}$ are internal edges (possibly terminal edges) of label $m+1$ at $w_2$,
and $e_3,e_3',\widetilde{e_3}$ are internal edges (possibly terminal edges) of label $m+1$ at $w_3$.

%%%%%%%%%%%%%%%%%%%%%%
%%%%%%%%%%%%%%%%%%%%%%
%%%%%%%%%%%%%%%%%%%%%%

\begin{lemma}
\label{D1ThreeD}
If there exists a minimal chart of type $(m;4,3)$
with $w(\Gamma\cap{\rm Int}D_1)=3$
such that $\widetilde{e_2}\cap \widetilde{e_3}$
is one white vertex,
then $e_2'=e_3'$ and
there exists a $3$-angled disk $D$ of $\Gamma_{m+1}$
without feelers
with $\partial D=e_2'\cup\widetilde{e_2}\cup \widetilde{e_3}$
and $w(\Gamma\cap{\rm Int}D)=0$.
\end{lemma}

\begin{figure}
\centerline{\includegraphics{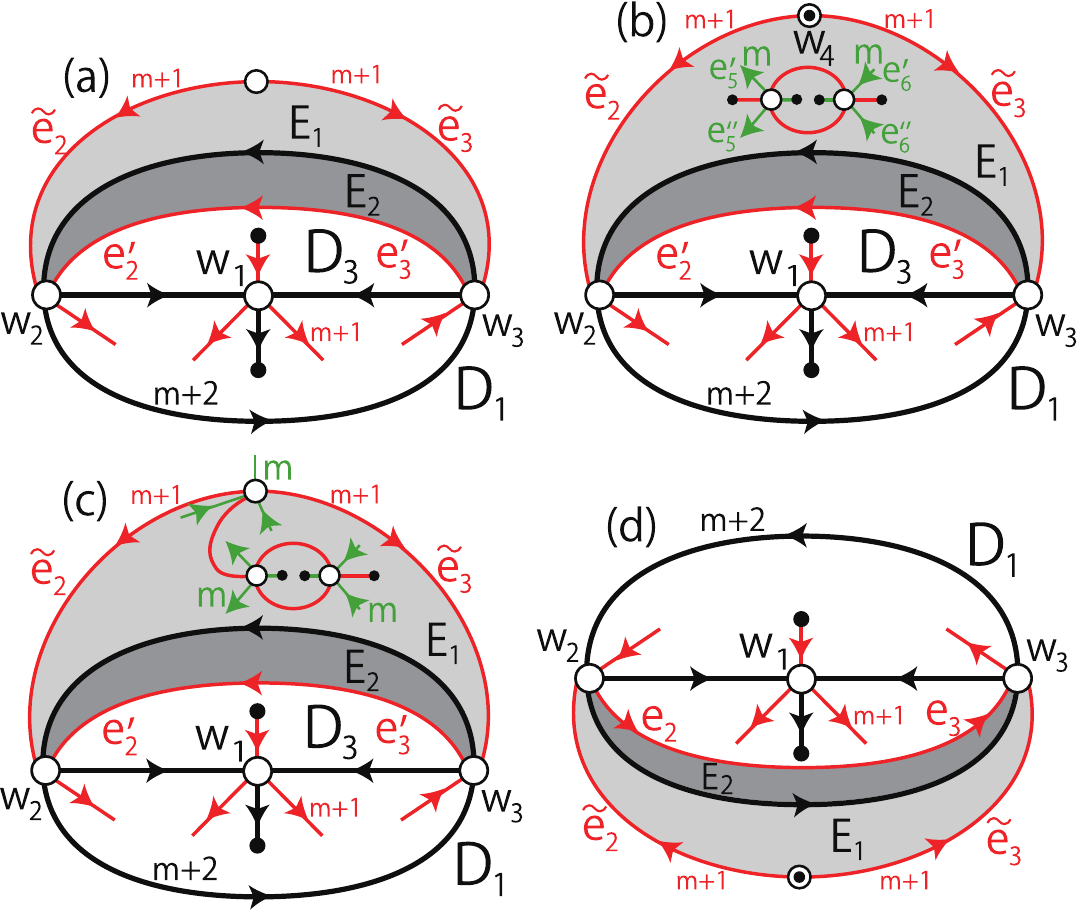}}
\caption{\label{fig29}
The light gray region is $E_1$,
the dark gray region is $E_2$, and $D=E_1\cup E_2$.}
\end{figure}

\begin{Proof}
Since $w(\Gamma\cap{\rm Int}D_1)=3$,
by Table~1 we have $w(\Gamma\cap{\rm Int}D_3)=0$.
Thus by Lemma~\ref{Theorem3AngledDisk}(a),
a regular neighborhood of $D_3$ contains
the pseudo chart as shown in Fig.~\ref{fig09}(b).
Hence $e_2'=e_3'$.

Let $D$ be the 3-angled disk of  $\Gamma_{m+1}$
with at most one feeler
and $\partial D=e_2'\cup\widetilde{e_2}\cup \widetilde{e_3}$.
The disk $D$ is divided by an internal edge of label $m+2$
into two disks.
One of the two disks is contained in $D_1$, say $E_1$.
The other is contained in $D_3$, say $E_2$
(see Fig.~\ref{fig29}(a)).
Since $w(\Gamma\cap{\rm Int}D_3)=0$,
we have $w(\Gamma\cap{\rm Int}E_2)=0$.

We shall show that $w(\Gamma\cap{\rm Int}E_1)=0$.
Suppose $w(\Gamma\cap{\rm Int}E_1)\ge1$.
By Condition~(d) in this section,
we have $w(\Gamma\cap{\rm Int}D_1)=w(\Gamma_{m+1}\cap{\rm Int}D_1)=3$.
Since $\widetilde{e_2}\cap \widetilde{e_3}$
is one white vertex,
by Lemma~\ref{Lemma2AngledDisks}(d)
a regular neighborhood of $D_1$ contains
one of the two pseudo charts as shown in 
Fig.~\ref{fig10}(e),(f).
Thus the condition  $w(\Gamma\cap{\rm Int}E_1)\ge1$
implies that there exists a 2-angled disk $E$
of $\Gamma_{m+1}$ in $E_1$
(see Fig.~\ref{fig29}(b),(c)).

{\bf Case (i).}
Suppose that $\Gamma$ contains the pseudo chart
as shown in Fig.~\ref{fig29}(b).
Let $w_4$ be the white vertex $\widetilde{e_2}\cap \widetilde{e_3}$,
and $e_4$ the terminal edge of label $m+1$ at $w_4$.
There are two cases:
$e_4\subset D$ or $e_4\not\subset D$.

If $e_4\subset D$, then we have
$w(\Gamma\cap({\rm Int} D-E))\ge1$ by IO-Calculation with respect
to $\Gamma_m$ in $Cl(D-E)$.
Thus $w(\Gamma\cap{\rm Int}D)\ge3$.
Hence by Condition~(a) of this section
$$\begin{array}{l}
w(\Gamma\cap{\rm Int}D_1)+w(\Gamma\cap{\rm Int}D_2)+w(\Gamma\cap{\rm Int}D_3)\vspace{2mm}\\
\ge w(w_4)+w(\Gamma\cap{\rm Int}D)+w(\Gamma\cap{\rm Int}D_2)\ge1+3+1=5.
\end{array}$$
This contradicts Condition~(c) of this section.

If $e_4\not\subset D$,
then there exists a lens in $D$
(i.e. $e_5'=e_6'$ and $e_5''=e_6''$).
Thus by Lemma~\ref{LensThreeWhiteVertex}
we have  $w(\Gamma\cap{\rm Int}D)\ge3$.
Hence we have the same contradiction.
Thus Case (i) does not occur.

{\bf Case (ii).}
Suppose that $\Gamma$ contains the pseudo chart
as shown in Fig.~\ref{fig29}(c).
Then by the similar way of the proof of Case (i),
we have $w(\Gamma\cap({\rm Int}D-E))\ge1$.
Thus we have the same contradiction.
Hence Case (ii) does not occur.

Therefore $w(\Gamma\cap{\rm Int}E_1)=0$.
Thus $w(\Gamma\cap{\rm Int}D)=w(\Gamma\cap{\rm Int}E_1)+w(\Gamma\cap{\rm Int}E_2)=0$.
Moreover by Lemma~\ref{Theorem3AngledDisk}(a)
the disk $D$ has no feelers.
\end{Proof}

%%%%%%%%%%%%%%%%%%
%%%%%%%%%%%%%%%%%% Figure
%%%%%%%%%%%%%%%%%%

%%%%%%%%%%%%%%%%%%%%%
%%%%%%%%%%%%%%%%%%%%%

\begin{lemma}
\label{D1OneD}
If there exists a minimal chart of type $(m;4,3)$
with $e_2'=e_3'$
such that $\widetilde{e_2}\cap \widetilde{e_3}$
is one white vertex,
then there exists a $3$-angled disk $D$ of $\Gamma_{m+1}$
without feelers
with $\partial D=e_2'\cup\widetilde{e_2}\cup \widetilde{e_3}$
and $w(\Gamma\cap{\rm Int}D)=0$.
\end{lemma}

\begin{Proof}
By Table~1, we have $w(\Gamma\cap{\rm Int}D_1)=1$ or
$w(\Gamma\cap{\rm Int}D_1)=3$.
If $w(\Gamma\cap{\rm Int}D_1)=3$,
we have the desired result from Lemma~\ref{D1ThreeD}.

Now suppose that $w(\Gamma\cap{\rm Int}D_1)=1$.
By Lemma~\ref{Lemma2AngledDisks}(b),
a regular neighborhood of $D_1$ contains
the pseudo chart as shown in Fig.~\ref{fig10}(a).
Let $w_4$ be the white vertex in ${\rm Int}D_1$.
Then $w_4\in \widetilde{e_2}\cap \widetilde{e_3}$.
Let $D$ be the 3-angled disk of  $\Gamma_{m+1}$
with at most one feeler
and $\partial D=e_2'\cup\widetilde{e_2}\cup \widetilde{e_3}$.

The disk $D$ is divided by an internal edge of label $m+2$
into two disks.
One of the two disks is contained in $D_1$, say $E_1$.
The other is contained in $D_3$, say $E_2$
(see Fig.~\ref{fig29}(a)).
Since $w_4\in\partial E_1$ and $w(\Gamma\cap{\rm Int}D_1)=1$,
we have $w(\Gamma\cap{\rm Int}E_1)=0$.

We shall show that $w(\Gamma\cap{\rm Int}E_2)=0$.
That is, we shall show that the following three cases do not occur:
(i)  $w(\Gamma\cap{\rm Int}E_2)=1$,
(ii)  $w(\Gamma\cap{\rm Int}E_2)=2$,
(iii)  $w(\Gamma\cap{\rm Int}E_2)\ge3$.

{\bf Case (i).}
There exists a connected component of $\Gamma_{m+1}$
with exactly one white vertex in $E_2$.
This contradicts Lemma~\ref{LemmaWithTerminal3}.
Hence Case (i) does not occur.

{\bf Case (ii).}
By Condition~(d) of this section,
we have $w(\Gamma\cap{\rm Int}E_2)=w(\Gamma_{m+1}\cap{\rm Int}E_2)=2$.
Thus there exists a disk $E_2'$ in ${\rm Int}E_2$
with $w(\Gamma\cap E_2')=w(\Gamma_{m+1}\cap E_2')=2$ and $\Gamma_{m+1}\cap\partial E_2'=\emptyset$.
Hence by Lemma~\ref{LemmaTwoWhiteVertices}
the disk $E_2'$ contains the pseudo chart as shown in
Fig.~\ref{fig03}(a)
(cf. Fig.~\ref{fig29}(b)).
By the similar way of the proof of Case (i) in 
Lemma~\ref{D1ThreeD}, 
we can show that $w(\Gamma\cap{\rm Int}D)\ge3$
and
 we have the same contradiction.
Therefore Case (ii) does not occur.

{\bf Case (iii).}
The condition $E_2\subset D_3$ implies 
$w(\Gamma\cap{\rm Int}D_3)\ge3$.
Thus by Condition~(a) of this section,
$$w(\Gamma\cap{\rm Int}D_1)+w(\Gamma\cap{\rm Int}D_2)+w(\Gamma\cap{\rm Int}D_3)\ge 1+1+3=5.$$
This contradicts Condition~(c) of this section.
Hence Case (iii) does not occur.

Therefore $w(\Gamma\cap{\rm Int}E_2)=0$.
Hence $w(\Gamma\cap{\rm Int}D)=w(\Gamma\cap{\rm Int}E_1)+w(\Gamma\cap{\rm Int}E_2)=0$.
Moreover by Lemma~\ref{Theorem3AngledDisk}(a)
the disk $D$ has no feelers.
\end{Proof}

\begin{lemma}
\label{D1OneD2}
If there exists a minimal chart of type $(m;4,3)$
with $w(\Gamma\cap{\rm Int}D_1)=1$ and $e_2=e_3$,
then there exists a $3$-angled disk $D$ of $\Gamma_{m+1}$
without feelers
with $\partial D=e_2\cup\widetilde{e_2}\cup \widetilde{e_3}$
and $w(\Gamma\cap{\rm Int}D)=0$.
\end{lemma}

\begin{Proof}
Since $w(\Gamma\cap{\rm Int}D_1)=1$,
by Lemma~\ref{Lemma2AngledDisks}(b)
a regular neighborhood of $D_1$ contains
the pseudo chart as shown in Fig.~\ref{fig10}(a).
Let $w_4$ be the white vertex in ${\rm Int}D_1$.
Then $w_4\in \widetilde{e_2}\cap \widetilde{e_3}$.
Let $D$ be the 3-angled disk of  $\Gamma_{m+1}$
with at most one feeler
and $\partial D=e_2\cup\widetilde{e_2}\cup \widetilde{e_3}$.

The disk $D$ is divided by an internal edge of label $m+2$
into two disks.
One of the two disks is contained in $D_1$, say $E_1$.
The other is contained in $D_2$, say $E_2$
(see Fig.~\ref{fig29}(d)).
Since $w_4\in\partial E_1$ and $w(\Gamma\cap{\rm Int}D_1)=1$,
we have $w(\Gamma\cap{\rm Int}E_1)=0$.

The edge $e_2$ divides the disk $D_2$ into two disks.
One of the two disks is the disk $E_2$.
Let $E_3$ be the other disk.
Then $E_3$ contains the terminal edge of label $m+2$ 
at $w_1$.
Thus by IO-Calculation with respect to $\Gamma_{m+1}$ in $E_3$,
we can show that 
\begin{enumerate}
\item[(1)]
$w(\Gamma\cap{\rm Int}E_3)\ge1$.
\end{enumerate}

We shall show that $w(\Gamma\cap{\rm Int}E_2)=0$.
That is, we shall show that the following three cases do not occur:
(i)  $w(\Gamma\cap{\rm Int}E_2)=1$,
(ii)  $w(\Gamma\cap{\rm Int}E_2)=2$,
(iii)  $w(\Gamma\cap{\rm Int}E_2)\ge3$.

{\bf Case (i).}
There exists a connected component of $\Gamma_{m+1}$
with exactly one white vertex in $E_2$.
This contradicts Lemma~\ref{LemmaWithTerminal3}.
Hence Case (i) does not occur.

{\bf Case (ii).}
By the similar way of the proof of Case (ii) in
Lemma~\ref{D1OneD}
we can show that $w(\Gamma\cap{\rm Int}D)\ge3$.
Hence by (1) we have
$$\begin{array}{l}
w(\Gamma\cap{\rm Int}D_1)+w(\Gamma\cap{\rm Int}D_2)+w(\Gamma\cap{\rm Int}D_3)\vspace{2mm}\\
\ge w(w_4)+w(\Gamma\cap{\rm Int}D)+
w(\Gamma\cap{\rm Int}E_3)\ge1+3+1=5.
\end{array}$$
This contradicts Condition~(c) of this section.
Therefore Case (ii) does not occur.

{\bf Case (iii).}
Since $E_2\subset D_2$, by (1) we have  
$w(\Gamma\cap{\rm Int}D_2)=w(\Gamma\cap{\rm Int}E_2)+w(\Gamma\cap{\rm Int}E_3)\ge3+1=4$.
Thus 
$$w(\Gamma\cap{\rm Int}D_1)+w(\Gamma\cap{\rm Int}D_2)+w(\Gamma\cap{\rm Int}D_3)\ge 1+4+0=5.$$
This contradicts Condition~(c) of this section.
Hence Case (iii) does not occur.

Therefore $w(\Gamma\cap{\rm Int}E_2)=0$.
Hence $w(\Gamma\cap{\rm Int}D)=w(\Gamma\cap{\rm Int}E_1)+w(\Gamma\cap{\rm Int}E_2)=0$.
Moreover by Lemma~\ref{Theorem3AngledDisk}(a)
the disk $D$ has no feelers.
\end{Proof}

\begin{corollary}
\label{D1OneDNotMinimal}
There does not exist a minimal chart of type $(m;4,3)$
with $w(\Gamma\cap{\rm Int}D_1)=1$, $e_2=e_3$ 
and $e_2'=e_3'$.
\end{corollary}

\begin{Proof}
Suppose that there exists a minimal chart of type $(m;4,3)$
with $w(\Gamma\cap{\rm Int}D_1)=1$, $e_2=e_3$ 
and $e_2'=e_3'$.
By Lemma~\ref{Lemma2AngledDisks}(b),
a regular neighborhood of $D_1$ contains
the pseudo chart as shown in Fig.~\ref{fig10}(a).
Thus $\widetilde{e_2}\cap \widetilde{e_3}$ is one white vertex, say $w_4$,
and 
\begin{enumerate}
\item[(1)]
$w_4$ is contained in a terminal edge of label $m+1$.
\end{enumerate}

Let $D,D'$ be the 3-angled disks of $\Gamma_{m+1}$
with $\partial D=e_2\cup \widetilde{e_2}\cup \widetilde{e_3}$
and  $\partial D'=e_2'\cup \widetilde{e_2}\cup \widetilde{e_3}$
such that $D\cap D'=\widetilde{e_2}\cup \widetilde{e_3}$.
By (1),
one of $D$ and $D'$ has a feeler.
However this contradicts
Lemma~\ref{D1OneD} and Lemma~\ref{D1OneD2}.
Therefore we complete the proof of 
Corollary~\ref{D1OneDNotMinimal}.
\end{Proof}

\begin{corollary}
\label{D1OneDNotMinimal2}
If there exists a minimal chart of type $(m;4,3)$
with $w(\Gamma\cap{\rm Int}D_1)=1$ and $e_2=e_3$,
then $e_3'$ is not a terminal edge.
\end{corollary}

\begin{Proof}
Since $w(\Gamma\cap{\rm Int}D_1)=1$ and $e_2=e_3$,
by Lemma~\ref{D1OneD2}
there exists a $3$-angled disk $D$ of $\Gamma_{m+1}$
without feelers
with $\partial D=e_2\cup\widetilde{e_2}\cup \widetilde{e_3}$
and $w(\Gamma\cap{\rm Int}D)=0$.

If $e_3'$ is a terminal edge, then
by Triangle Lemma(Lemma~\ref{LemmaTriangle})
we have $w(\Gamma\cap{\rm Int}D)\ge1$.
This is a contradiction.
Hence $e_3'$ is not a terminal edge.
\end{Proof}

Similarly we can show the following corollary:

\begin{corollary}
\label{D1OneDNotMinimal3}
If there exists a minimal chart of type $(m;4,3)$
with $w(\Gamma\cap{\rm Int}D_1)=1$ and $e_2'=e_3'$, 
then $e_2$ is not a terminal edge.
\end{corollary}

%\newpage
%%%%%%%%%%%%%%%%%%%%%%%%%%
%%%%%%%%%%%%%%%%%%%%%%%%%%
%%%%%%%%%%%%%%%%%%%%%%%%%%
\section{Case of a chart $\Gamma$ of type $(4,3)$
with $w(\Gamma\cap{\rm Int}D_1)=3$}
\label{s:D1ThreeD2OneD3Zero}

In this section, we use the notations as shown in Fig.~\ref{fig28}.
In this section
we shall show the following proposition.

%%%%%%%%%%%%%%%%%%%%%
%%%%%%%%%%%%%%%%%%%%%
\begin{proposition}
\label{D1ThreeD2OneD3Zero}
There does not exist a minimal chart of type $(m;4,3)$
with $w(\Gamma\cap{\rm Int}D_1)=3$.
\end{proposition}

\begin{figure}
\centerline{\includegraphics{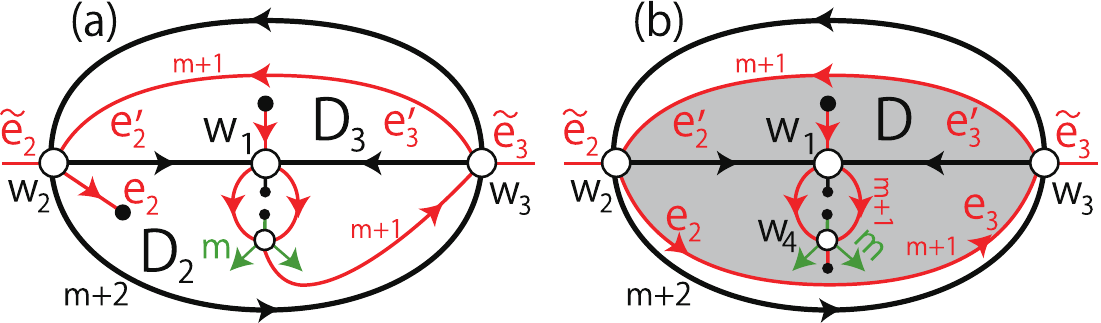}}
\caption{The gray region is the 2-angled disk $D$.
\label{fig30}
}
\end{figure}

\begin{lemma}
\label{D2OneD3Zero}
If there exists a minimal chart of type $(m;4,3)$
with $w(\Gamma\cap{\rm Int}D_1)=3$,
$w(\Gamma\cap{\rm Int}D_2)=1$ and
$w(\Gamma\cap{\rm Int}D_3)=0$,
then the chart $\Gamma$ contains the pseudo chart
as shown in Fig.~\ref{fig30}$($a$)$.
\end{lemma}

%%%%%%%%%%%%%%%%%%
%%%%%%%%%%%%%%%%%% Figure
%%%%%%%%%%%%%%%%%%

\begin{Proof}
Since $w(\Gamma\cap{\rm Int}D_3)=0$,
by Lemma~\ref{Theorem3AngledDisk}(a)
a regular neighborhood of $D_3$
contains the pseudo chart as shown in Fig.~\ref{fig09}(b).
Hence we have $e_2'=e_3'$.

Since $w(\Gamma\cap{\rm Int}D_2)=1$ and since
$D_2$ has one feeler,
by Lemma~\ref{Theorem3AngledDisk}(b)
a regular neighborhood of $D_2$
contains one of the two pseudo charts as shown in
Fig.~\ref{fig09}(g),(h).

If a regular neighborhood of $D_2$
contains  the  pseudo chart as shown in
Fig.~\ref{fig09}(h),
then we have $e_2=e_3$
(see Fig.~\ref{fig30}(b)).
Moreover there exists a 2-angled disk $E$ of $\Gamma_{m+1}$
in $D_2$.
Let $w_4$ be the white vertex in ${\rm Int}D_2$.
Let $D$ be the 2-angled disk of $\Gamma_{m+1}$ in $D_2\cup D_3$ with $\partial D=e_2\cup e_2'$.
Then we have $w(\Gamma\cap({\rm Int}D-E))\ge1$ by IO-Calculation with respect to $\Gamma_m$
in $Cl(D-E)$.
Hence 
$$\begin{array}{rl}
1=1+0& =w(\Gamma\cap{\rm Int}D_2)+w(\Gamma\cap{\rm Int}D_3)\vspace{2mm}\\
& \ge w(\Gamma\cap({\rm Int}D-E))+w(w_4)\ge 1+1=2.
\end{array}$$
This is a contradiction.

Thus a regular neighborhood of $D_2$
contains  the  pseudo chart as shown in
Fig.~\ref{fig09}(g).
Therefore $\Gamma$ contains the pseudo chart
as shown in Fig.~\ref{fig30}(a).
\end{Proof}

%%%%%%%%%%%%%%%%%%%%%%%
%%%%%%%%%%%%%%%%%%%%%%%

{\it Proof of Proposition~\ref{D1ThreeD2OneD3Zero}.}
Suppose that there exists a minimal chart $\Gamma$
of type $(m;4,3)$ with $w(\Gamma\cap{\rm Int}D_1)=3$.
Then by Table~1 in Section~\ref{s:Chart43},
we have $w(\Gamma\cap{\rm Int}D_2)=1$ and 
$w(\Gamma\cap{\rm Int}D_3)=0$.
Thus by Lemma~\ref{D2OneD3Zero},
the chart $\Gamma$ contains the pseudo chart
as shown in Fig.~\ref{fig30}(a).
By Condition~(d) of Section~\ref{s:Chart43},
we have $w(\Gamma\cap{\rm Int}D_1)=w(\Gamma_{m+1}\cap{\rm Int}D_1)=3$.
Hence by Lemma~\ref{Lemma2AngledDisks}(d),
a regular neighborhood of $D_1$
contains one of the five pseudo charts as shown in
Fig.~\ref{fig10}(b),(c),(d),(e),(f).

\begin{figure}
\centerline{\includegraphics{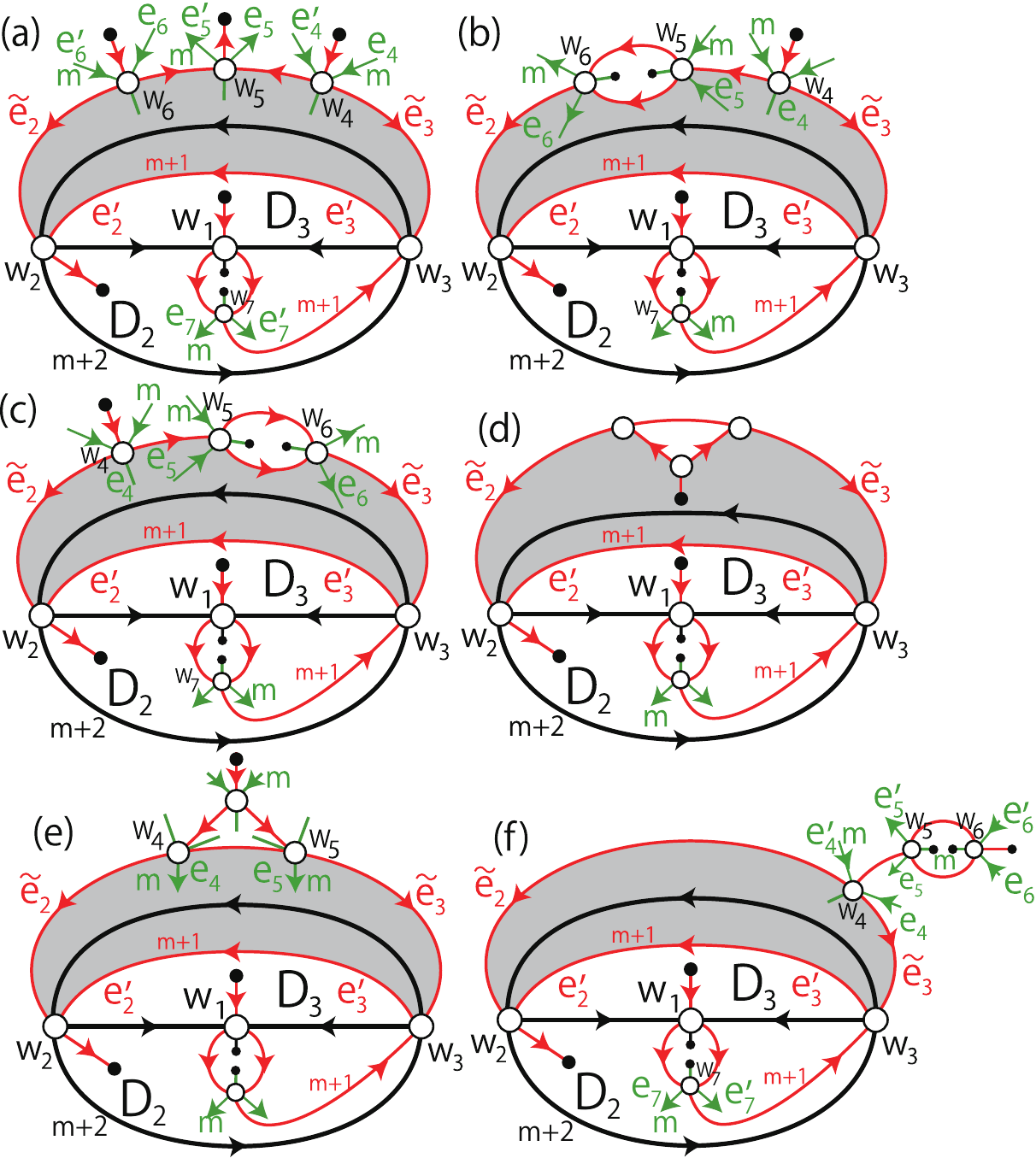}}
\caption{The gray regions are the disk $D$.
\label{fig31}}
\end{figure}

{\bf Case (i).}
Suppose that a regular neighborhood of $D_1$
contains  the  pseudo chart as shown in
Fig.~\ref{fig10}(b).
Let $D$ be the 5-angled disk of $\Gamma_{m+1}$
with $\partial D\supset e_2'\cup\widetilde{e_2}\cup \widetilde{e_3}$ and $w(\Gamma\cap{\rm Int}D)=0$.
Then $D$ has a nice edge of label $m+2$
connecting $w_2$ and $w_3$.
Thus by Lemma~\ref{NiceEdge5AngledDiskAtMostOneFeeler} and
Lemma~\ref{NiceEdge5AngledDiskOneFeeler},
the disk $D$ has no feelers
(see Fig.~\ref{fig31}(a)).

We use the notations as shown in Fig.~\ref{fig31}(a), where 
$w_4,w_5,w_6$ are the white vertices in $\partial D$,
$w_7$ are the white vertex in ${\rm Int}D_2$, and
 $e_i,e_i'$ $(i=4,5,6,7)$ are internal edges (possibly terminal edges)
of label $m$ at $w_i$
such that 
$e_4,e_4',e_6,e_6'$ are oriented inward at
$w_4,w_4,w_6,w_6$, respectively, and
$e_5,e_5',e_7,e_7'$ are oriented outward at
$w_5,w_5,w_7,w_7$, respectively.
Then we can show that for $i=4,5,6,7$
neither $e_i$ nor $e_i'$ is a terminal edge.
Let $E$ be the 2-angled disk of $\Gamma_{m+1}$ in $D_2$.
Since $S^2-(D\cup E)$ does not contain any white vertices,
the edge $e_5$ must be equal to $e_4'$ or $e_6'$.

If $e_5=e_4'$, then there exists a lens $E'$ with $\partial E'\supset e_5$. This contradicts  Lemma~\ref{NoLens}.
If $e_5=e_6'$, then $e_5'=e_6$ and
there exists a lens $E'$ with $\partial E'\supset e_5'$.
This contradicts Lemma~\ref{NoLens}.
Hence Case (i) does not occur.

{\bf Case (ii).}
Suppose that a regular neighborhood of $D_1$
contains the pseudo chart as shown in
Fig.~\ref{fig10}(c).
Let $D$ be the special 5-angled disk of $\Gamma_{m+1}$
with $\partial D\supset e_2'\cup\widetilde{e_2}\cup \widetilde{e_3}$ and $w(\Gamma\cap{\rm Int}D)=0$.
Then $D$ has a nice edge of label $m+2$
connecting $w_2$ and $w_3$.
Thus by Lemma~\ref{NiceEdge5AngledDiskAtMostOneFeeler} and
Lemma~\ref{NiceEdge5AngledDiskOneFeeler},
the disk $D$ has no feelers
(see Fig.~\ref{fig31}(b) and (c)).

Let $w_4,w_5,w_6$ be the white vertices in $\partial D$
different from $w_2$ and $w_3$ such that
$w_4$ is contained in a terminal edge of label $m+1$.
Let $e_4,e_5,e_6$ be internal edges
(possibly terminal edges) of label $m$ at $w_4,w_5,w_6$
in $D$, respectively.
Then we can show that neither $e_5$ nor $e_6$ 
is a terminal edge.
Thus the condition $w(\Gamma\cap{\rm Int}D)=0$ implies 
that $e_5=e_6$ and the edge $e_4$ is a terminal edge.
Hence there exists a lens in $D$.
This contradicts Lemma~\ref{NoLens}.
Thus Case (ii) does not occur.

{\bf Case (iii).}
Suppose that a regular neighborhood of $D_1$
contains  the  pseudo chart as shown in
Fig.~\ref{fig10}(d).
Then there exists a special 4-angled disk or 
5-angled disk of $\Gamma_{m+1}$, say $D$, such that
 $\partial D\supset e_2'\cup\widetilde{e_2}\cup \widetilde{e_3}$ and $w(\Gamma\cap{\rm Int}D)=0$
(see Fig.~\ref{fig31}(d) and (e)).

If $D$ is a 5-angled disk 
(see Fig.~\ref{fig31}(d)),
then $D$ has exactly one feeler
and a nice edge of label $m+2$ connecting $w_2$ and $w_3$.
However this contradicts 
Lemma~\ref{NiceEdge5AngledDiskOneFeeler}.
Hence $D$ is a 4-angled disk 
(see Fig.~\ref{fig31}(e)).

Let $w_4,w_5$ be the white vertices in $\partial D$
different from $w_2$ and $w_3$.
Let $e_4,e_5$ be internal edges
(possibly terminal edges) of label $m$ at $w_4,w_5$
in $D$, respectively.
Then $e_4,e_5$ are oriented outward at $w_4,w_5$,
respectively.
Thus $e_4\not=e_5$.
Hence  the condition $w(\Gamma\cap{\rm Int}D)=0$ implies 
that both of $e_4$ and $e_5$ are terminal edges.
However we can show that one of $e_4$ and $e_5$
is not middle at $w_4$ or $w_5$.
This contradicts Assumption~\ref{AssumeTerminal}. 
Hence Case (iii) does not occur.

{\bf Case (iv).}
Suppose that a regular neighborhood of $D_1$
contains the pseudo chart as shown in
Fig.~\ref{fig10}(e).
Then $\widetilde{e_2}\cap \widetilde{e_3}$ is one white vertex, say $w_4$.
By Lemma~\ref{D1ThreeD},
there exists a 3-angled disk $D$ of $\Gamma_{m+1}$ without
 feelers such that
 $\partial D=e_2'\cup\widetilde{e_2}\cup \widetilde{e_3}$
and $w(\Gamma\cap{\rm Int}D)=0$.
Thus by Lemma~\ref{Theorem3AngledDisk}(a),
a regular neighborhood of $D$ contains the pseudo chart
as shown in Fig.~\ref{fig09}(a).
Since the white vertex $w_2$ is contained 
in the terminal edge $e_2$
of label $m+1$ and since the white vertex $w_4$ is contained in a terminal edge
of label $m+1$,
 by Triangle Lemma(Lemma~\ref{LemmaTriangle})
we have  $w(\Gamma\cap{\rm Int}D)\ge1$.
This is a contradiction.
Hence Case (iv) does not occur.

{\bf Case (v).}
Suppose that a regular neighborhood of $D_1$
contains  the pseudo chart as shown in
Fig.~\ref{fig10}(f).
Then $\widetilde{e_2}\cap \widetilde{e_3}$ is one white vertex, say $w_4$.
By Lemma~\ref{D1ThreeD},
there exists a 3-angled disk $D$ of $\Gamma_{m+1}$ without
 feelers such that
 $\partial D=e_2'\cup\widetilde{e_2}\cup \widetilde{e_3}$
and $w(\Gamma\cap{\rm Int}D)=0$
(see Fig.~\ref{fig31}(f)).

We use the notations as shown in Fig.~\ref{fig31}(f),
where $w_4,w_5,w_6$ are white vertices in ${\rm Int}D_1$,
$w_7$ is the white vertex in ${\rm Int}D_2$,
and $e_i,e_i'$ $(i=4,5,6,7)$ are internal edges
(possibly terminal edges) of label $m$ at $w_i$
such that $e_4,e_4',e_6,e_6'$ are oriented inward at 
$w_4,w_4,w_6,w_6$, respectively, and
 $e_5,e_5',e_7,e_7'$ are oriented outward at 
$w_5,w_5,w_7,w_7$, respectively.
By the similar way of Case (i),
we can show that $e_5=e_4$ or $e_5=e_6$.
Thus there exists a lens $E$ with $\partial E\supset e_5$.
This contradicts Lemma~\ref{NoLens}.
Hence Case (v) does not occur.

Therefore all the five cases do not occur.
Hence there does not exist a minimal chart $\Gamma$
of type $(m;4,3)$ with $w(\Gamma\cap{\rm Int}D_1)=3$.
We complete the proof of Proposition~\ref{D1ThreeD2OneD3Zero}.
{\hfill {$\square$}\vspace{1.5em}}

%%%%%%%%%%%%%%%%%%
%%%%%%%%%%%%%%%%%% Figure
%%%%%%%%%%%%%%%%%%

%\newpage
%%%%%%%%%%%%%%%%%%%%%%%%%%
%%%%%%%%%%%%%%%%%%%%%%%%%%
%%%%%%%%%%%%%%%%%%%%%%%%%%
\section{Case of a chart $\Gamma$ 
with $w(\Gamma\cap{\rm Int}D_1)=1$}
\label{s:D1OneD2One}

In this section
we shall show that there does not exist a minimal chart of type $(m;4,3)$
with $w(\Gamma\cap{\rm Int}D_1)=1$ and
$w(\Gamma\cap{\rm Int}D_2)=1$.

\begin{figure}
\centerline{\includegraphics{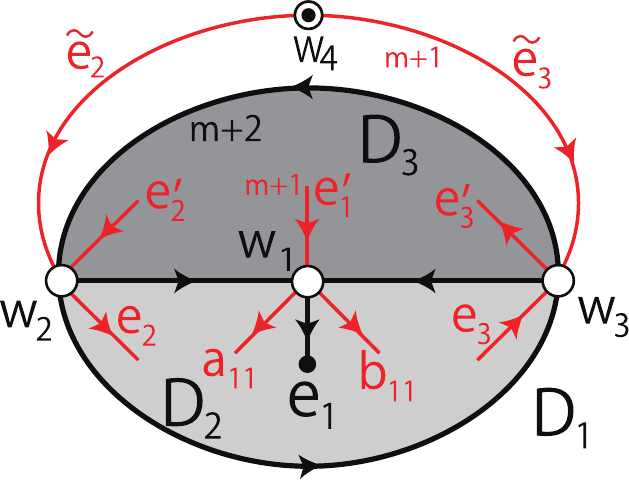}}
\caption{The light gray region is $D_2$. The dark gray region is $D_3$.
\label{fig32} }
\end{figure}

From now on throughout in this paper,
 we use the notations as shown in Fig.~\ref{fig28}.
Moreover we may assume that
\begin{enumerate}
\item[(e)] $w(\Gamma\cap{\rm Int}D_1)=1$.
\end{enumerate}
By Lemma~\ref{Lemma2AngledDisks}(b),
a regular neighborhood of $D_1$ contains
the pseudo chart as shown in Fig.~\ref{fig10}(a)
(see Fig.~\ref{fig32}).
Let $w_4$ be the white vertex in ${\rm Int}D_1$.
Then $w_4\in \widetilde{e_2}\cap \widetilde{e_3}$.
Also we use the notations as shown in Fig.~\ref{fig32}.

%%%%%%%%%%%%%%%%%%
%%%%%%%%%%%%%%%%%% Figure
%%%%%%%%%%%%%%%%%%

\begin{lemma}
\label{D1OneD4Angled}
There does not exist a minimal chart of type $(m;4,3)$
with $w(\Gamma\cap{\rm Int}D_1)=1$ such that
$e_2'\cap e_3'$ is one white vertex.
\end{lemma}

\begin{figure}
\centerline{\includegraphics{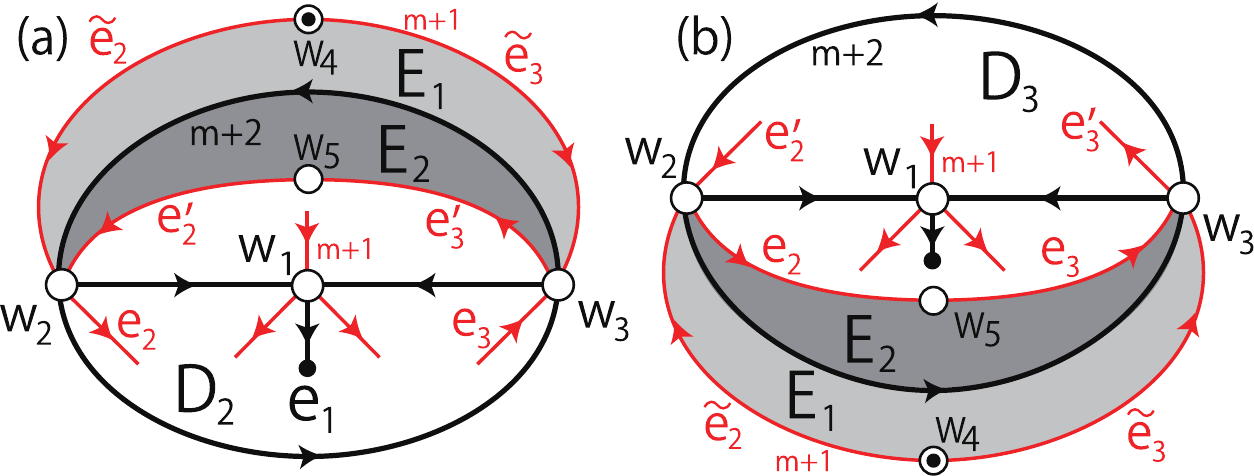}}
\caption{The light gray region is $E_1$, the dark gray region is $E_2$, and $D=E_1\cup E_2$.
\label{fig33} }
\end{figure}

\begin{Proof}
Suppose that there exists a minimal chart of type $(m;4,3)$
with $w(\Gamma\cap{\rm Int}D_1)=1$ such that
\begin{enumerate}
\item[(1)] $e_2'\cap e_3'$ is one white vertex.
\end{enumerate}

\clearpage

Let $D$ be the 4-angled disk of $\Gamma_{m+1}$
with $\partial D=e_2'\cup e_3'\cup\widetilde{e_2}\cup \widetilde{e_3}$ and $w_1\not\in D$
(see Fig.~\ref{fig33}(a)).
The disk $D$ is divided by an internal edge of label $m+2$
into two disks.
One of the two disks is contained in $D_1$, say $E_1$.
The other is contained in $D_3$, say $E_2$.
Since $w_4\in\partial E_1$ and $w(\Gamma\cap{\rm Int}D_1)=1$,
we have $w(\Gamma\cap{\rm Int}E_1)=0$.

Next we shall show that $w(\Gamma\cap{\rm Int}E_2)=0$.
Suppose that $w(\Gamma\cap{\rm Int}E_2)\ge1$.
By Condition~(d) of Section~\ref{s:Chart43},
we have $w(\Gamma_{m+1}\cap{\rm Int}E_2)\ge1$.
Thus there exists a disk $E_2'$ in ${\rm Int}E_2$
with $w(\Gamma_{m+1}\cap E_2')\ge1$ such that 
$\Gamma_{m+1}\cap\partial E_2'$ is at most one point.
Hence by Lemma~\ref{LemmaAtMostOnePoint}
we have $w(\Gamma\cap E_2')\ge2$.
Thus $w(\Gamma\cap{\rm Int}E_2)\ge2$.
Hence by (1) and Condition~(a) 
of Section~\ref{s:Chart43},
we have 
$$\begin{array}{l}
w(\Gamma\cap{\rm Int}D_1)+w(\Gamma\cap{\rm Int}D_2)+w(\Gamma\cap{\rm Int}D_3)\vspace{2mm}\\
\ge w(w_4)+w(\Gamma\cap{\rm Int}D_2)+ w(\Gamma\cap{\rm Int}E_2)+w(e_2'\cap e_3')\vspace{2mm}\\
\ge1+1+2+1=5.
\end{array}$$
This contradicts Condition~(c) of 
Section~\ref{s:Chart43}.
Therefore $w(\Gamma\cap{\rm Int}E_2)=0$.
Hence $w(\Gamma\cap{\rm Int}D)=w(\Gamma\cap{\rm Int}E_1)+w(\Gamma\cap{\rm Int}E_2)=0$.

Now the internal edge of label $m+2$ in $D$
is a nice edge with respect to $D$.
Since $w(\Gamma\cap{\rm Int}D)=0$,
 by Lemma~\ref{NiceEdge4AngledDisk}
the 4-angled disk $D$ has no feelers.

Let $w_5$ be the white vertex $e_2'\cap e_3'$,
and $e_5$ an internal edge (possibly a terminal edge) 
of label $m$ at $w_5$ in $D$.
Since $e_2'$ is oriented outward at $w_5$
and since $e_3'$ is oriented inward at $w_5$
(see Fig.~\ref{fig33}(a)),
the edge $e_5$ is not middle at $w_5$.
Thus by Assumption~\ref{AssumeTerminal},
the edge $e_5$ is not a terminal edge.
Hence we have $e_5\ni w_4$.
Since the internal edge of label $m+2$ in $D$
is middle at $w_2$ and
since the white vertex $w_4$ is contained in
a terminal edge of label $m+1$,
by Lemma~\ref{4AngledDiskCrossNotMinimal}
the chart $\Gamma$ is not minimal.
This contradicts the fact that $\Gamma$ is minimal.
Therefore we complete the proof of Lemma~\ref{D1OneD4Angled}.
\end{Proof}

%%%%%%%%%%%%%%%%%%
%%%%%%%%%%%%%%%%%% Figure
%%%%%%%%%%%%%%%%%%

\begin{lemma}
\label{D1OneD4Angled2}
There does not exist a minimal chart of type $(m;4,3)$
with $w(\Gamma\cap{\rm Int}D_1)=1$ such that
$e_2\cap e_3$ is one white vertex.
\end{lemma}

\begin{Proof}
Suppose that there exists a minimal chart of type $(m;4,3)$
with $w(\Gamma\cap{\rm Int}D_1)=1$ such that
\begin{enumerate}
\item[(1)] $e_2\cap e_3$ is one white vertex.
\end{enumerate}

Let $D$ be the 4-angled disk of  $\Gamma_{m+1}$
with $\partial D=e_2\cup e_3\cup\widetilde{e_2}\cup \widetilde{e_3}$ and $w_1\not\in D$
(see Fig.~\ref{fig33}(b)).
The disk $D$ is divided by an internal edge of label $m+2$
into two disks.
One of the two disks is contained in $D_1$, say $E_1$.
The other is contained in $D_2$, say $E_2$.
Since $w_4\in\partial E_1$ and $w(\Gamma\cap{\rm Int}D_1)=1$,
we have $w(\Gamma\cap{\rm Int}E_1)=0$.

Next we shall show that $w(\Gamma\cap{\rm Int}E_2)=0$.
Suppose that $w(\Gamma\cap{\rm Int}E_2)\ge1$.
By the similar way of the proof of Lemma~\ref{D1OneD4Angled},
we can show that $w(\Gamma\cap{\rm Int}E_2)\ge2$.

The union $e_2\cup e_3$ divides the disk $D_2$ 
into two disks.
One of the two disks is the disk $E_2$.
Let $E_3$ be the other disk.
By IO-Calculation with respect to $\Gamma_{m+1}$ in $E_3$,
we can show that $w(\Gamma\cap{\rm Int}E_3)\ge1$.
Hence by (1)
we have 
$$\begin{array}{l}
w(\Gamma\cap{\rm Int}D_1)+w(\Gamma\cap{\rm Int}D_2)+w(\Gamma\cap{\rm Int}D_3)\vspace{2mm}\\
\ge w(w_4)+w(\Gamma\cap{\rm Int}E_2)+ w(\Gamma\cap{\rm Int}E_3)+w(e_2\cap e_3)\vspace{2mm}\\
\ge1+2+1+1=5.
\end{array}$$
This contradicts Condition~(c) of 
Section~\ref{s:Chart43}.
Therefore $w(\Gamma\cap{\rm Int}E_2)=0$.
Hence $w(\Gamma\cap{\rm Int}D)=w(\Gamma\cap{\rm Int}E_1)+w(\Gamma\cap{\rm Int}E_2)=0$.

By the similar way of the end of 
the proof of Lemma~\ref{D1OneD4Angled},
we can show that the 4-angled disk $D$ satisfies 
the condition of Lemma~\ref{4AngledDiskCrossNotMinimal}.
Thus by Lemma~\ref{4AngledDiskCrossNotMinimal},
the chart $\Gamma$ is not minimal.
This contradicts the fact that $\Gamma$ is minimal.
Therefore we complete the proof of 
Lemma~\ref{D1OneD4Angled2}.
\end{Proof}

%%%%%%%%%%%%%%%%%%%%%%%
%%%%%%%%%%%%%%%%%%%%%%%

\begin{lemma}
\label{D1OneD2OneD3TwoNotFig12}
There does not exist a minimal chart $\Gamma$ of type $(m;4,3)$
with $w(\Gamma\cap{\rm Int}D_1)=1$,
$w(\Gamma\cap{\rm Int}D_2)=1$ such that
a regular neighborhood of $D_3$
contains the pseudo chart
as shown in Fig.~\ref{fig13}.
\end{lemma}

\begin{figure}
\centerline{\includegraphics{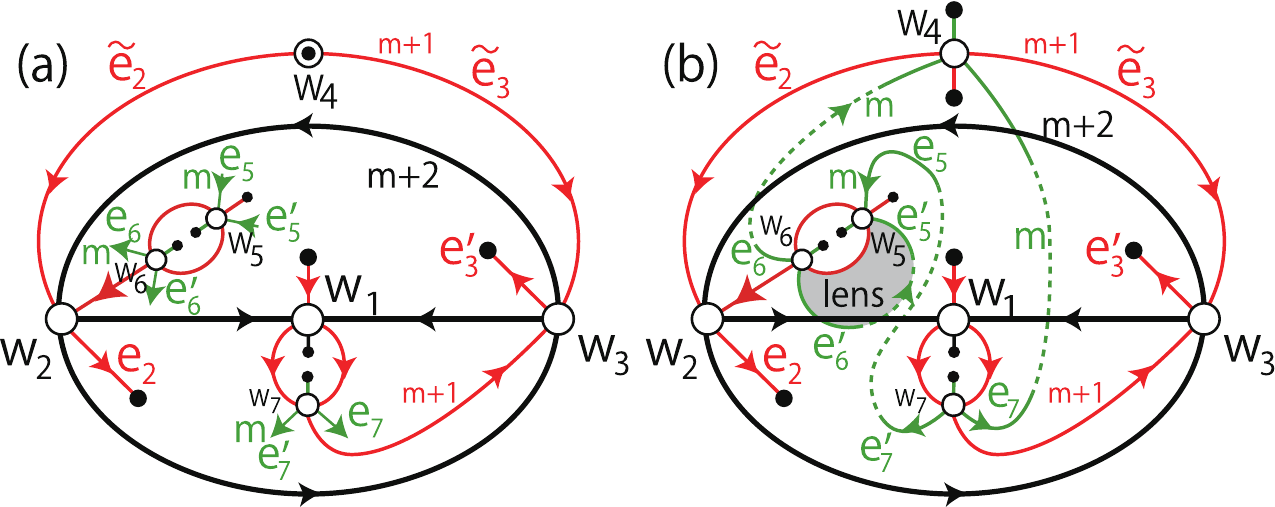}}
\caption{A regular neighborhood of $D_3$
contains the pseudo chart
as shown in Fig.~\ref{fig13}.
The gray region is a lens.
\label{fig34}
}
\end{figure}

\begin{Proof}
Suppose that there exists a minimal chart 
$\Gamma$ of type $(m;4,3)$
with $w(\Gamma\cap{\rm Int}D_1)=1$,
$w(\Gamma\cap{\rm Int}D_2)=1$ such that
a regular neighborhood of $D_3$
contains the pseudo chart
as shown in Fig.~\ref{fig13}.
Then
\begin{enumerate}
\item[(1)] $e_3'$ is a terminal edge.
\end{enumerate}

Since $w(\Gamma\cap{\rm Int}D_2)=1$ and
since $D_2$ has one feeler,
by Lemma~\ref{Theorem3AngledDisk}(b)
a regular neighborhood of $D_2$ contains one of the 
pseudo charts as shown in Fig.~\ref{fig09}(g),(h).

Suppose that a regular neighborhood of $D_2$ contains the 
pseudo charts as shown in Fig.~\ref{fig09}(h).
Then we have $e_2=e_3$.
By Corollary~\ref{D1OneDNotMinimal2},
the edge $e_3'$ is not a terminal edge.
This contradicts (1).
Thus a regular neighborhood of $D_2$ does not contain the 
pseudo charts as shown in Fig.~\ref{fig09}(h).

Suppose that a regular neighborhood of $D_2$ contains the 
pseudo charts as shown in Fig.~\ref{fig09}(g)
(see Fig.~\ref{fig34}(a)).
We use the notations as shown in Fig.~\ref{fig34}(a),
where $w_5,w_6$ are white vertices in ${\rm Int}D_3$, 
and
$w_7$ is the white vertex in ${\rm Int}D_2$.
Let $e_i,e_i'$ ($i=4,5,6,7$) be internal edges
(possibly terminal edges) of label $m$ at $w_i$
such that $e_4,e_4',e_5,e_5'$ are oriented inward at 
$w_4,w_4,w_5,w_5$, respectively, and
 $e_6,e_6',e_7,e_7'$ are oriented outward at 
$w_6,w_6,w_7,w_7$, respectively 
(see Fig.~\ref{fig34}(a)).
Then we can show that 
none of the eight edges are terminal edges.
Since $\Gamma_m$ contains exactly four white vertices
$w_4,w_5,w_6,w_7$,
we have $e_5=e_6$ or $e_5=e_7'$.

If $e_5=e_6$, then there exists a lens $E$ with
$\partial E\supset e_5$.
This contradicts Lemma~\ref{NoLens}.
If $e_5=e_7'$,
then $e_6\cap e_7\ni w_4$ and $e_5'=e_6'$
(see Fig.~\ref{fig34}(b)).
Thus there exists a lens $E$ with
$\partial E\supset e_5'$.
This contradicts Lemma~\ref{NoLens}.
Hence a regular neighborhood of $D_2$ does not contain the 
pseudo charts as shown in Fig.~\ref{fig09}(g).
\end{Proof}

%%%%%%%%%%%%%%%%%%
%%%%%%%%%%%%%%%%%% Figure
%%%%%%%%%%%%%%%%%%

%%%%%%%%%%%%%%%%%%%%%%%
%%%%%%%%%%%%%%%%%%%%%%%

\begin{lemma}
\label{D1OneD2OneD3TwoPseudoChart}
If there exists a minimal chart $\Gamma$ of type $(m;4,3)$
with $w(\Gamma\cap{\rm Int}D_1)=1$,
$w(\Gamma\cap{\rm Int}D_2)=1$ and $w(\Gamma\cap{\rm Int}D_3)=2$,
then a regular neighborhood of $D_2$ contains the 
pseudo chart as shown in Fig.~\ref{fig09}$($g$)$.
\end{lemma}

\begin{figure}
\centerline{\includegraphics{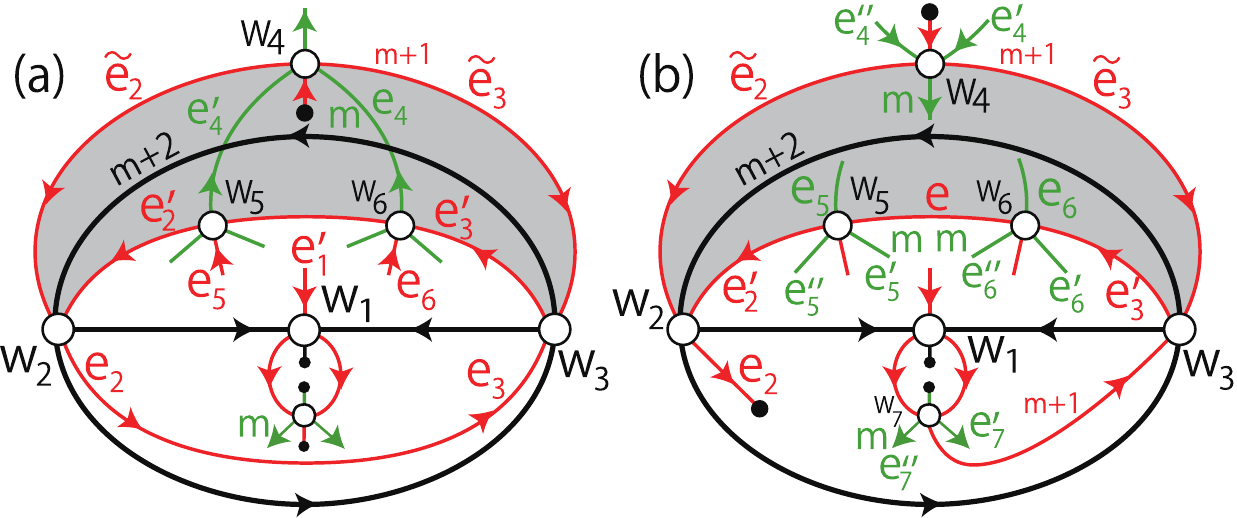}}
\caption{
The gray regions are the 5-angled disk $E$.
\label{fig35}
}
\end{figure}

\begin{Proof}
We use the notations as shown in Fig.~\ref{fig32}.
Since $w(\Gamma\cap{\rm Int}D_2)=1$
and since $D_2$ has one feeler,
by Lemma~\ref{Theorem3AngledDisk}(b)
a regular neighborhood of $D_2$ contains one of the 
pseudo charts as shown in Fig.~\ref{fig09}(g),(h).

Suppose that a regular neighborhood of $D_2$ contains the 
pseudo chart as shown in Fig.~\ref{fig09}(h).
Then $e_2=e_3$.
Thus by Corollary~\ref{D1OneDNotMinimal},
we have 
\begin{enumerate}
\item[(1)] $e_2'\not=e_3'$.
\end{enumerate}

Let $D$ be the 3-angled disk of $\Gamma_{m+1}$
with 
$\partial D=e_2\cup\widetilde{e_2}\cup \widetilde{e_3}$
and  $w(\Gamma\cap{\rm Int}D)=0$.
Then by Lemma~\ref{Theorem3AngledDisk}(a)
\begin{enumerate}
\item[(2)] the disk $D$ has no feelers.
\end{enumerate}

By Condition~(d) in Section~\ref{s:Chart43},
we have $w(\Gamma\cap{\rm Int}D_3)=w(\Gamma_{m+1}\cap{\rm Int}D_3)=2$.
Thus by Lemma~\ref{3angledDiskNoFeelerTwoWhiteVertex}(b)
and Lemma~\ref{D1OneD2OneD3TwoNotFig12},
there exists a simple arc $L$ in $\Gamma_{m+1}\cap D_3$
connecting $w_2$ and $w_3$.
Since $w(\Gamma\cap{\rm Int}D_3)=2$,
by (1) there are two cases:
(i) $w({\rm Int}L)=1$,
(ii) $w({\rm Int}L)=2$.

{\bf Case (i).}
Since $w({\rm Int}L)=1$,
the intersection $e_2'\cap e_3'$ is one white vertex.
However such a chart does not exist by Lemma~\ref{D1OneD4Angled}.
Thus Case (i) does not occur.

{\bf Case (ii).}
Let $E$ be the special 5-angled disk of $\Gamma_{m+1}$
with $\partial E\supset e_2'\cup e_3'\cup\widetilde{e_2}\cup \widetilde{e_3}$.
Then $w(\Gamma\cap{\rm Int}E)=0$.
By Lemma~\ref{NiceEdge5AngledDiskAtMostOneFeeler},
the disk $E$ has at most one feeler.
Thus by (2),
the disk $E$ has one feeler at $w_4$.
Hence by Lemma~\ref{NiceEdge5AngledDiskOneFeeler},
we can assume that
a regular neighborhood of $E$ contains the pseudo chart
as shown in 
Fig.~\ref{fig24}(c)
(see Fig.~\ref{fig35}(a)).

We use the notations as shown in Fig.~\ref{fig35}(a),
where $w_5,w_6$ are white vertices in $\partial E$
different from $w_2,w_3,w_4$.
Let $e_4,e_4'$ be internal edges of label $m$ at 
$w_4$ in $E$,
and $e_1',e_5,e_6$ internal edges (possibly terminal edges)
of label $m+1$ at $w_1,w_5,w_6$ in $D_3$,
respectively, with
$e_5\not\subset E$ and $e_6\not\subset E$.
Since both of $e_4$ and $e_4'$ are oriented inward at $w_4$,
the edges $e_5,e_6$ are oriented inward at $w_5,w_6$,
respectively.
Now $e_1'$ is oriented inward at $w_1$.

Let $E'=Cl(D_3-E)$. Then $w(\Gamma\cap{\rm Int}E')=0$.
Since the edges $e_1',e_5,e_6$ are oriented inward at $w_1,w_5,w_5$, all of  $e_1',e_5,e_6$ are terminal edges.

Since both of $e_3',e_6$ are oriented inward at $w_6$,
the internal edge of label $m+1$ connecting $w_5$ and $w_6$
is oriented from $w_6$ to $w_5$.
Since the edge $e_5$ is oriented inward at $w_5$,
 the terminal edge $e_5$ is not middle at $w_5$.
This contradicts Assumption~\ref{AssumeTerminal}.
Hence Case (ii) does not occur.

Therefore 
a regular neighborhood of $D_2$ does not contain the 
pseudo chart as shown in Fig.~\ref{fig09}(h).
Thus 
a regular neighborhood of $D_2$ contains the 
pseudo chart as shown in Fig.~\ref{fig09}(g).
\end{Proof}

%%%%%%%%%%%%%%%%%%
%%%%%%%%%%%%%%%%%% Figure
%%%%%%%%%%%%%%%%%%

%%%%%%%%%%%%%%%%%%%%%%%
%%%%%%%%%%%%%%%%%%%%%%%
%%%%%%%%%%%%%%%%%%%%%%%

\begin{proposition}
\label{D1OneD2OneD3Two}
There does not exist a minimal chart of type $(m;4,3)$
with $w(\Gamma\cap{\rm Int}D_1)=1$ and
$w(\Gamma\cap{\rm Int}D_2)=1$.
\end{proposition}

\begin{Proof}
Suppose that there exists a minimal chart of type $(m;4,3)$
with $w(\Gamma\cap{\rm Int}D_1)=1$ and
$w(\Gamma\cap{\rm Int}D_2)=1$.
We use the notations as shown in
Fig.~\ref{fig32} where
\begin{enumerate}
\item[(1)] $e_2'$ is oriented inward at $w_2$, and
$e_3'$ is oriented outward at $w_3$.
\end{enumerate}

By Table~1 in Section~\ref{s:Chart43},
we have $w(\Gamma\cap{\rm Int}D_3)=2$.
Thus by Lemma~\ref{D1OneD2OneD3TwoPseudoChart},
a regular neighborhood of $D_2$ contains
the pseudo chart as shown in Fig.~\ref{fig09}(g).
Hence
\begin{enumerate}
\item[(2)] $e_2$ is a terminal edge.
\end{enumerate}

By Condition~(d) in Section~\ref{s:Chart43},
we have $w(\Gamma\cap{\rm Int}D_3)=w(\Gamma_{m+1}\cap{\rm Int}D_3)=2$.
Thus by Lemma~\ref{3angledDiskNoFeelerTwoWhiteVertex}(b)
and Lemma~\ref{D1OneD2OneD3TwoNotFig12},
there exists a simple arc $L$ in $\Gamma_{m+1}\cap D_3$
connecting $w_2$ and $w_3$.
Since $w(\Gamma\cap{\rm Int}D_3)=2$,
 there are three cases:
(i) $w({\rm Int}L)=0$,
(ii) $w({\rm Int}L)=1$,
(iii) $w({\rm Int}L)=2$.

{\bf Case (i).}
Since $w({\rm Int}L)=0$,
we have $e_2'=e_3'$.
By Corollary~\ref{D1OneDNotMinimal3},
the edge $e_2$ is not a terminal edge.
This contradicts (2).
Hence Case (i) does not occur.

{\bf Case (ii).}
Since $w({\rm Int}L)=1$,
the intersection $e_2'\cap e_3'$ is one white vertex.
However such a chart does not exist 
by Lemma~\ref{D1OneD4Angled}.
Hence Case (ii) does not occur.

{\bf Case (iii).}
Let $E$ be the special 5-angled disk of $\Gamma_{m+1}$
with 
 $\partial E\supset e_2'\cup e_3'\cup\widetilde{e_2}\cup \widetilde{e_3}$.
Then $w(\Gamma\cap{\rm Int}E)=0$.
By Lemma~\ref{NiceEdge5AngledDiskAtMostOneFeeler},
the disk $E$ has at most one feeler.

If $E$ has one feeler,
then we have the same contradiction
by the similar way of the proof of Case (ii) in
Lemma~\ref{D1OneD2OneD3TwoPseudoChart}.
Hence $E$ has no feelers
(see Fig.~\ref{fig35}(b)).

We use the notations as shown in Fig.~\ref{fig35}(b),
where $w_5,w_6$ are white vertices in $\partial E$,
$e$ is the internal edge of label $m+1$ in $\partial E$
with $\partial e=\{w_5,w_6\}$, and
$e_5,e_6$ are internal edges (possibly terminal edges)
of label $m$ at $w_5,w_6$ in $E$, respectively.

We shall show that $e$ is oriented from $w_5$ to $w_6$.
If  $e$ is oriented from $w_6$ to $w_5$,
then neither $e_5$ nor $e_6$ is middle at $w_5$ or $w_6$.
Thus by Assumption~\ref{AssumeTerminal},
neither $e_5$ nor $e_6$ is a terminal edge.
Hence the condition $w(\Gamma\cap{\rm Int}E)=0$ implies
$e_5=e_6$.
Thus there exists a lens $E'$ with $\partial E'=e_5\cup e$.
This contradicts Lemma~\ref{NoLens}.
Hence
\begin{enumerate}
\item[(3)] $e$ is oriented from $w_5$ to $w_6$.
\end{enumerate}

Let $e_i',e_i''$ ($i=4,5,6,7$) be internal edges
(possibly terminal edges) of label $m$ at $w_i$
such that $e_4',e_4''$ are oriented inward at 
$w_4$, 
$e_7',e_7''$ are oriented outward at 
$w_7$, and 
none of $e_5',e_5'',e_6',e_6''$ are contained in $E$
(see Fig.~\ref{fig35}(b)).
Then we can show that 
\begin{enumerate}
\item[(4)] 
none of $e_4',e_4'',e_7',e_7''$ are terminal edges.
\end{enumerate}
By (1) and (3),
the edges $e_5',e_5''$ are oriented inward at 
$w_5$ and 
the edges $e_6',e_6''$ are oriented outward at 
$w_6$.
Thus we can show that  
\begin{enumerate}
\item[(5)]none of $e_5',e_5'',e_6',e_6''$ are terminal edges.
\end{enumerate}
Since $\Gamma_m$ contains exactly four white vertices
$w_4,w_5,w_6,w_7$,
the conditions (4) and (5) imply $e_5'=e_6''$.
Thus there exists a lens $E'$ with $\partial E'=e_5'\cup e$.
This contradicts Lemma~\ref{NoLens}.
Hence Case (iii) does not occur.

Therefore all the three cases do not occur.
We complete the proof of Proposition~\ref{D1OneD2OneD3Two}.
\end{Proof}

%\newpage
%%%%%%%%%%%%%%%%%%%%%%%%%%
%%%%%%%%%%%%%%%%%%%%%%%%%%
%%%%%%%%%%%%%%%%%%%%%%%%%%
\section{Case of a chart $\Gamma$ 
with $w(\Gamma\cap{\rm Int}D_1)=1$ and $w(\Gamma\cap{\rm Int}D_2)=2$}
\label{s:D1OneD2Two}

In this section
we shall show that there does not exist a minimal chart of type $(m;4,3)$
with $w(\Gamma\cap{\rm Int}D_1)=1$ and
$w(\Gamma\cap{\rm Int}D_2)=2$.
In this section,
 we use the notations as shown in Fig.~\ref{fig28}
and Fig.~\ref{fig32}.

%%%%%%%%%%%%%%%%%%%%%%%
%%%%%%%%%%%%%%%%%%%%%%%

\begin{lemma}
\label{D1OneD2TwoD3OnePseudoChart}
If there exists a minimal chart $\Gamma$ of type $(m;4,3)$
with $w(\Gamma\cap{\rm Int}D_1)=1$,
$w(\Gamma\cap{\rm Int}D_2)=2$ and $w(\Gamma\cap{\rm Int}D_3)=1$,
then a regular neighborhood of $D_3$ contains one of
the two pseudo charts as shown in
Fig.~\ref{fig09}$($d$)$ and $($e$)$.
\end{lemma}

%%%%%%%%%%%%%%%%%%%%%%
%%%%%%%%%%%%%%%%%%%%%%%%

\begin{Proof}
Since $w(\Gamma\cap{\rm Int}D_3)=1$ and since
$D_3$ has no feelers,
by Lemma~\ref{Theorem3AngledDisk}(b)
a regular neighborhood of $D_3$ contains one of
the three pseudo charts as shown in Fig.~\ref{fig09}(c),(d),(e).

If a regular neighborhood of $D_3$ contains 
the pseudo chart as shown in Fig.~\ref{fig09}(c),
then $e_2'\cap e_3'$ is one white vertex.
However such a chart does not exist by 
Lemma~\ref{D1OneD4Angled}.
Thus a regular neighborhood of $D_3$ contains one of
the two pseudo charts as shown in Fig.~\ref{fig09}(d),(e).
\end{Proof}

%%%%%%%%%%%%%%%%%%%%
%%%%%%%%%%%%%%%%%%%%
%%%%%%%%%%%%%%%%%%%%

\begin{lemma}
\label{D1OneD2TwoD3OneNotFig13}
There does not exist a minimal chart $\Gamma$ of type $(m;4,3)$
with $w(\Gamma\cap{\rm Int}D_1)=1$,
$w(\Gamma\cap{\rm Int}D_2)=2$, 
$w(\Gamma\cap{\rm Int}D_3)=1$ such that
a regular neighborhood of $D_2$
contains the pseudo chart
as shown in Fig.~\ref{fig14}$($a$)$.
\end{lemma}

\begin{figure}
\centerline{\includegraphics{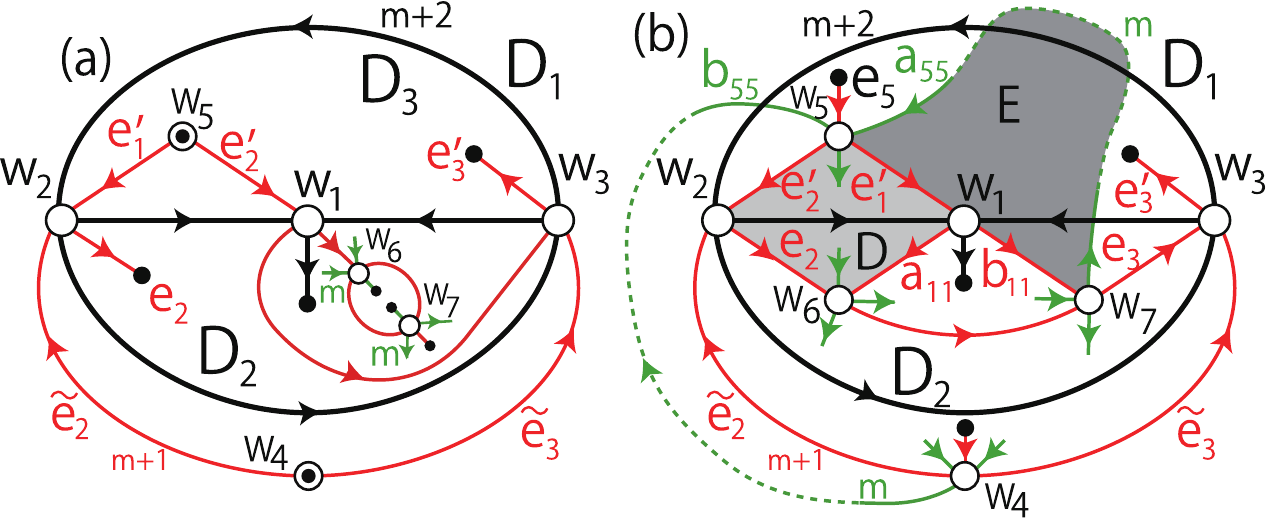}}
\caption{\label{fig36}
(a) A regular neighborhood of $D_2$
contains the pseudo chart
as shown in Fig.~\ref{fig14}$($a$)$.
(b) The light gray region is the 4-angled disk $D$, the dark gray region is the disk $E$.
}
\end{figure}

\begin{Proof}
Suppose that there exists a minimal chart $\Gamma$ of 
type $(m;4,3)$
with $w(\Gamma\cap{\rm Int}D_1)=1$,
$w(\Gamma\cap{\rm Int}D_2)=2$, 
$w(\Gamma\cap{\rm Int}D_3)=1$ such that
a regular neighborhood of $D_2$
contains the pseudo chart
as shown in 
Fig.~\ref{fig14}(a).
By Lemma~\ref{D1OneD2TwoD3OnePseudoChart},
a regular neighborhood of $D_3$ contains one of
the two pseudo charts as shown in Fig.~\ref{fig09}(d),(e).
Thus
the chart $\Gamma$ contains the pseudo chart 
as shown in Fig.~\ref{fig36}(a).
We use the notations as shown in Fig.~\ref{fig36}(a),
where $w_5$ is the white vertex in ${\rm Int}D_3$, and
$w_6,w_7$ are white vertices in ${\rm Int}D_2$.

Let $e_i,e_i'$ ($i=4,5,6,7$) be internal edges
(possibly terminal edges) of label $m$ at $w_i$
such that $e_4,e_4',e_5,e_5',e_6,e_6'$
are oriented inward at $w_4,w_4,w_5,w_5,w_6,w_6$,
respectively, and
$e_7,e_7'$ are oriented outward at $w_7$.
We can show that none of the eight edges are terminal edges.
Thus by IO-Calculation with respect to $\Gamma_m$ in $S^2$,
there exists a white vertex of $\Gamma_m$
different from $w_4,w_5,w_6,w_7$.
Hence $w(\Gamma_m)\ge5$.
This contradicts the fact that 
$\Gamma$ is of type $(m;4,3)$.
Therefore we complete the proof of 
Lemma~\ref{D1OneD2TwoD3OneNotFig13}. 
\end{Proof}

%%%%%%%%%%%%%%%%%%
%%%%%%%%%%%%%%%%%% Figure
%%%%%%%%%%%%%%%%%%

%%%%%%%%%%%%%%%%%%%%%%%
%%%%%%%%%%%%%%%%%%%%%%%
%%%%%%%%%%%%%%%%%%%%%%%

\begin{proposition}
\label{D1OneD2TwoD3One}
There does not exist a minimal chart of type $(m;4,3)$
with $w(\Gamma\cap{\rm Int}D_1)=1$ and
$w(\Gamma\cap{\rm Int}D_2)=2$.
\end{proposition}

\begin{Proof}
Suppose that there exists a minimal chart of type $(m;4,3)$
with $w(\Gamma\cap{\rm Int}D_1)=1$ and
$w(\Gamma\cap{\rm Int}D_2)=2$.
We use the notations as shown in Fig.~\ref{fig32}, where
\begin{enumerate}
\item[(1)] $e_3$ is oriented inward at $w_3$ but not middle at $w_3$,
\item[(2)] neither $a_{11}$ nor $b_{11}$ is middle at $w_1$.
\end{enumerate}

By Table~1 in Section~\ref{s:Chart43},
we have $w(\Gamma\cap{\rm Int}D_3)=1$.
Thus by Lemma~\ref{D1OneD2TwoD3OnePseudoChart},
 a regular neighborhood of $D_3$ contains one of
the two pseudo charts as shown in
Fig.~\ref{fig09}(d) and (e).
Thus 
\begin{enumerate}
\item[(3)] $e_3'$ is a terminal edge.
\end{enumerate}

By Condition~(d) in Section~\ref{s:Chart43},
we have $w(\Gamma\cap{\rm Int}D_2)=w(\Gamma_{m+1}\cap{\rm Int}D_2)=2$.
Thus by Lemma~\ref{3angledDiskOneFeelerTwoWhiteVertex}, 
either
\begin{enumerate}
\item[(4)] there exists a simple arc $L$ in 
$(\Gamma_{m+1}-w_1)\cap D_2$
connecting $w_2$ and $w_3$, or
\item[(5)] a regular neighborhood of $D_2$ contains one of the five pseudo charts as shown in Fig.~\ref{fig14}.
\end{enumerate}

If $\Gamma$ satisfies Condition (5),
then by (1)
a regular neighborhood of $D_2$ contains the pseudo chart as shown in Fig.~\ref{fig14}(a).
However such a chart does not exist by Lemma~\ref{D1OneD2TwoD3OneNotFig13}.
Thus $\Gamma$ satisfies Condition (4).

Since $w(\Gamma\cap{\rm Int}D_2)=2$,
 there are three cases:
(i) $w({\rm Int}L)=0$,
(ii) $w({\rm Int}L)=1$,
(iii) $w({\rm Int}L)=2$.

{\bf Case (i).}
Since $w({\rm Int}L)=0$,
we have $e_2=e_3$.
By Corollary~\ref{D1OneDNotMinimal2},
the edge $e_3'$ is not a terminal edge.
This contradicts (3).
Hence Case (i) does not occur.

{\bf Case (ii).}
Since $w({\rm Int}L)=1$,
the intersection $e_2\cap e_3$ is one white vertex.
However such a chart does not exist 
by Lemma~\ref{D1OneD4Angled2}.
Hence Case (ii) does not occur.

{\bf Case (iii).}
Let $w_6,w_7$ be the white vertices in ${\rm Int}D_2$
with $w_6\in e_2$ and $w_7\in e_3$.
Since neither $a_{11}$ nor $b_{11}$ is a terminal edge
by (2) and Assumption~\ref{AssumeTerminal},
we have $a_{11}\ni w_6$ and $b_{11}\ni w_7$.

Let $D$ be the 4-angled disk of $\Gamma_{m+1}$
with $\partial D=e_1'\cup e_2'\cup e_2\cup a_{11}$
and $w(\Gamma\cap{\rm Int}D)=0$.
Then by IO-Calculation with respect to $\Gamma_m$ in $D$,
\begin{enumerate}
\item[(6)] the terminal edge $e_5$ of label $m+1$
at $w_5$ is not contained in $D$.
\end{enumerate}

Let $a_{55},b_{55}$ be internal edges (possibly terminal edges)
of label $m+1$ oriented inward at $w_5$ such that
$a_{55},e_5,b_{55}$ lie anticlockwise around $w_5$.
Let $D'$ be the 6-angled disk of $\Gamma_{m+1}$ with
 $\partial D'=e_1'\cup e_2'\cup \widetilde{e_2}\cup \widetilde{e_3}\cup e_3\cup b_{11}$ and $w(\Gamma\cap{\rm Int}D')=0$.
Then by (6),
we have $a_{55}\cup b_{55}\subset D'$.
We can show that neither $a_{55}$ nor $b_{55}$ is a 
terminal edge. 
Thus by IO-Calculation with respect to $\Gamma_m$ in $D'$,
the terminal edge of label $m+1$ at $w_4$
is not contained in $D'$.
Moreover
the condition $w(\Gamma\cap{\rm Int}D')=0$
implies that $a_{55}\ni w_7$ and $b_{55}\ni w_4$.
Therefore the chart $\Gamma$ contains the pseudo chart
as shown in Fig.~\ref{fig36}(b).

Let $E$ be the disk in $D'$ with $\partial E=e_1'\cup b_{11}\cup a_{55}$ and $w(\Gamma\cap{\rm Int}E)=0$.
By using New Disk Lemma(Lemma~\ref{NewDiskLemma}),
we can assume that
\begin{enumerate}
\item[(7)] the chart $\Gamma$ is $(E,a_{55})$-arc free.
\end{enumerate}

We shall show that $\Gamma_{m+2}\cap a_{55}$ is one crossing.
Suppose that $\Gamma_{m+2}\cap a_{55}$ contains at least
two crossings.
Then $\Gamma_{m+2}\cap E$ consists of at least two proper arcs.
One of these arcs contains the white vertex $w_1$.
One of these arcs does not contain $w_1$,
say $\gamma$.
The arc $\gamma$ of label $m+2$
does not intersect the edges $e_1'$ and $b_{11}$ of 
label $m+1$.
Hence $\partial \gamma\subset a_{55}$.
Thus $\gamma$ is an $(E,a_{55})$-arc of label $m+2$.
This contradicts (7).
Hence $\Gamma_{m+2}\cap a_{55}$ is one crossing.

By C-II moves,
we move the black vertex in $e_5$ near
the black vertex in $e_3'$ along the edge $a_{55}$.
Then  by applying a C-I-M2 move between $e_5$ and $e_3'$,
we obtain a free edge of label $m+1$
(cf. Fig.~\ref{fig21}(d),(e)).
Hence the complexity of the chart decreases.
This contradicts the fact that $\Gamma$ is minimal.
Thus Case (iii) does not occur.

Therefore all the three cases do not occur.
Thus we complete the proof of Proposition~\ref{D1OneD2TwoD3One}.
\end{Proof}

%%%%%%%%%%%%%%%%%%%%%%%%%%
%%%%%%%%%%%%%%%%%%%%%%%%%%
%%%%%%%%%%%%%%%%%%%%%%%%%%

%\newpage
%%%%%%%%%%%%%%%%%%%%%%%%%%
%%%%%%%%%%%%%%%%%%%%%%%%%%
\section{Case of a chart $\Gamma$ 
with $w(\Gamma\cap{\rm Int}D_1)=1$ and $w(\Gamma\cap{\rm Int}D_2)=3$}
\label{s:D1OneD2Three}

In this section,
we shall show the main theorem,
i.e. we shall show that there does not exist a minimal chart
of type $(m;4,3)$ with $w(\Gamma\cap{\rm Int}D_1)=1$ and $w(\Gamma\cap{\rm Int}D_1)=3$.
In this section,
 we use the notations as shown in Fig.~\ref{fig28}
and Fig.~\ref{fig32}.

%%%%%%%%%%%%%%%%%%%%%%%
%%%%%%%%%%%%%%%%%%%%%%%

\begin{figure}
\centerline{\includegraphics{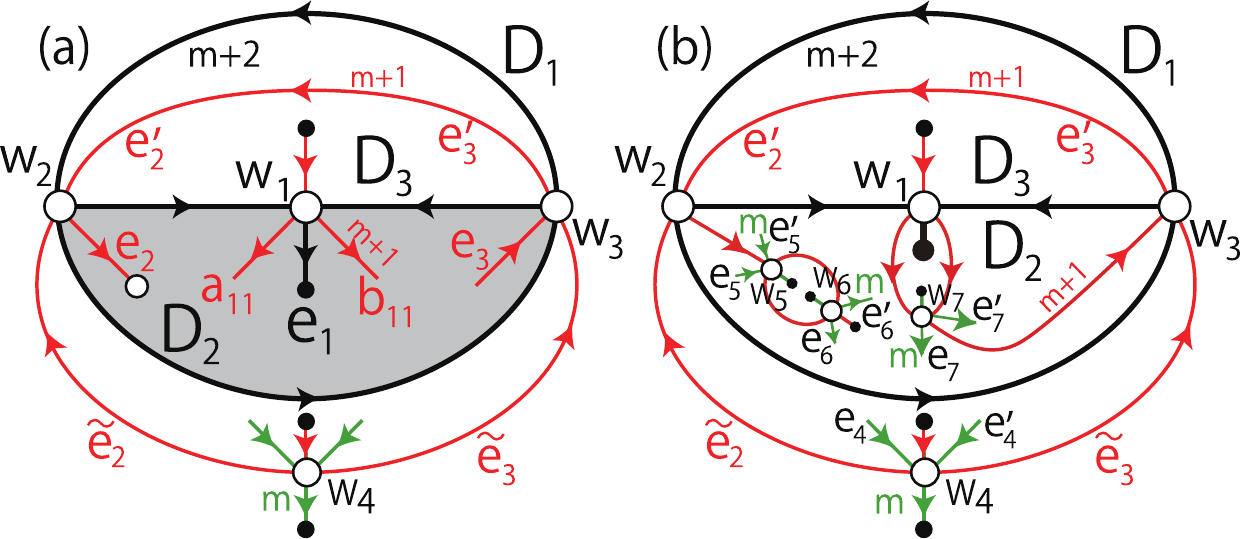}}
\caption{\label{fig37}
(a) The gray region is the disk $D_2$.
(b) A regular neighborhood of $D_2$
contains the pseudo chart
as shown in 
Fig.~\ref{fig15}(a).
}
\end{figure}

\begin{lemma}
\label{D1OneD2ThreeD3ZeroPseudoChart}
If there exists a minimal chart $\Gamma$ of type $(m;4,3)$
with $w(\Gamma\cap{\rm Int}D_1)=1$,
$w(\Gamma\cap{\rm Int}D_2)=3$ and $w(\Gamma\cap{\rm Int}D_3)=0$,
then the chart $\Gamma$ contains the pseudo chart
as shown in 
Fig.~\ref{fig37}$($a$)$.
Moreover
\begin{enumerate}
\item[{\rm (a)}] the edge $e_2$ contains a white vertex
 in ${\rm Int}D_2$,
\item[{\rm (b)}] neither $a_{11}$ nor $b_{11}$ is a terminal edge.
\end{enumerate}
\end{lemma}

%%%%%%%%%%%%%%%%%%
%%%%%%%%%%%%%%%%%% Figure
%%%%%%%%%%%%%%%%%%

%%%%%%%%%%%%%%%%%%%%%%
%%%%%%%%%%%%%%%%%%%%%%%%

\begin{Proof}
We use the notations as shown in Fig.~\ref{fig32},
where 
\begin{enumerate}
\item[(1)] 
$\widetilde{e_2}\cap \widetilde{e_3}$
is a white vertex $w_4$.
\end{enumerate}
Since $w(\Gamma\cap{\rm Int}D_3)=0$,
by Lemma~\ref{Theorem3AngledDisk}(a)
a regular neighborhood of $D_3$ contains the pseudo chart as shown in Fig.~\ref{fig09}(b).
Thus 
\begin{enumerate}
\item[(2)] 
$e_2'=e_3'$.
\end{enumerate}
Hence there exists the 3-angled disk $D$ 
of $\Gamma_{m+1}$ with
$\partial D=e_2'\cup\widetilde{e_2}\cup \widetilde{e_3}$ and
$w(\Gamma\cap {\rm Int}D)=0$.
Then by Lemma~\ref{Theorem3AngledDisk}(a)
a regular neighborhood of $D$ contains the pseudo chart as shown in Fig.~\ref{fig09}(a).
Hence the terminal edge of label $m+1$ at $w_4$
is not contained in $D$.

By (2) and Corollary~\ref{D1OneDNotMinimal3},
the edge $e_2$ is not a terminal edge.
Moreover by Corollary~\ref{D1OneDNotMinimal},
we have $e_2\not=e_3$.
Thus the edge $e_2$ contains a white vertex
in ${\rm Int}D_2$.
Hence the chart $\Gamma$ contains the pseudo chart
as shown in Fig.~\ref{fig37}(a).

Since neither $a_{11}$ nor $b_{11}$ is middle at $w_1$,
by Assumption~\ref{AssumeTerminal}
neither $a_{11}$ nor $b_{11}$ is a terminal edge.
\end{Proof}

%%%%%%%%%%%%%%%%%%%%
%%%%%%%%%%%%%%%%%%%%
%%%%%%%%%%%%%%%%%%%%

\begin{lemma}
\label{D1OneD2ThreeD3ZeroNotFig14A}
There does not exist a minimal chart $\Gamma$ of type $(m;4,3)$
with $w(\Gamma\cap{\rm Int}D_1)=1$,
$w(\Gamma\cap{\rm Int}D_2)=3$, 
$w(\Gamma\cap{\rm Int}D_3)=0$ such that
a regular neighborhood of $D_2$
contains the pseudo chart
as shown in Fig.~\ref{fig15}$($a$)$.
\end{lemma}

\begin{Proof}
Suppose that there exists a minimal chart $\Gamma$ of 
type $(m;4,3)$
with $w(\Gamma\cap{\rm Int}D_1)=1$,
$w(\Gamma\cap{\rm Int}D_2)=3$, 
$w(\Gamma\cap{\rm Int}D_3)=0$ such that
a regular neighborhood of $D_2$
contains the pseudo chart
as shown in 
Fig.~\ref{fig15}(a).
By Lemma~\ref{D1OneD2ThreeD3ZeroPseudoChart},
the chart $\Gamma$ contains the pseudo chart 
as shown in Fig.~\ref{fig37}(b).
We use the notations as shown in Fig.~\ref{fig37}(b).

Let $e_i,e_i'$ ($i=4,5,6,7$) be internal edges
(possibly terminal edges) of label $m$ at $w_i$
such that $e_4,e_4',e_5,e_5'$
are oriented inward at $w_4,w_4,w_5,w_5$,
respectively, and
$e_6,e_6',e_7,e_7'$ are oriented outward 
at $w_6,w_6,w_7,w_7$, respectively.
We can show that none of the eight edges are terminal edges.
Thus by looking at the edge $e_6$,
we have $e_6=e_5$.
Hence there exists a lens $E$ with $\partial E\supset e_6$.
This contradicts Lemma~\ref{NoLens}.
Therefore we complete the proof of 
Lemma~\ref{D1OneD2ThreeD3ZeroNotFig14A}. 
\end{Proof}

\begin{lemma}
\label{D1OneD2ThreeD3ZeroNotFig14B}
There does not exist a minimal chart $\Gamma$ of type $(m;4,3)$
with $w(\Gamma\cap{\rm Int}D_1)=1$,
$w(\Gamma\cap{\rm Int}D_2)=3$, 
$w(\Gamma\cap{\rm Int}D_3)=0$ such that
a regular neighborhood of $D_2$
contains the pseudo chart
as shown in Fig.~\ref{fig15}$($b$)$.
\end{lemma}

\begin{figure}
\centerline{\includegraphics{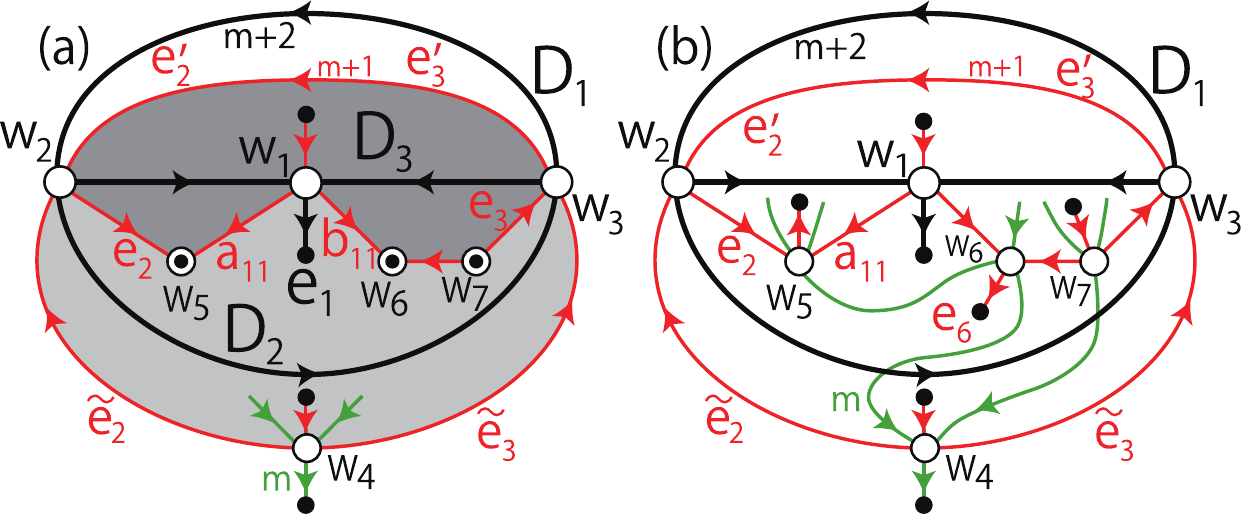}}
\caption{\label{fig38}
A regular neighborhood of $D_2$
contains the pseudo chart
as shown in Fig.~\ref{fig15}$($b$)$. 
The light gray region is the 7-angled disk $D$. the dark gray region is the 6-angled disk $D'$.}
\end{figure}

\begin{Proof}
Suppose that there exists a minimal chart $\Gamma$ of 
type $(m;4,3)$
with $w(\Gamma\cap{\rm Int}D_1)=1$,
$w(\Gamma\cap{\rm Int}D_2)=3$, 
$w(\Gamma\cap{\rm Int}D_3)=0$ such that
a regular neighborhood of $D_2$
contains the pseudo chart
as shown in 
Fig.~\ref{fig15}(b).
By Lemma~\ref{D1OneD2ThreeD3ZeroPseudoChart},
the chart $\Gamma$ contains the pseudo chart 
as shown in Fig.~\ref{fig38}(a).
We use the notations as shown in 
Fig.~\ref{fig38}(a).

Let $D$ be the special 7-angled disk of $\Gamma_{m+1}$
with at least one feeler, and
$D'$ the special 6-angled disk of $\Gamma_{m+1}$
 such that
$\partial D\supset e_2\cup e_3\cup a_{11}\cup b_{11}\cup\widetilde{e_2}\cup \widetilde{e_3}$ and 
$\partial D'\supset e_2'\cup e_2\cup e_3\cup a_{11}\cup b_{11}$.
Then 
\begin{enumerate}
\item[(1)] $w(\Gamma\cap{\rm Int}D)=0$, 
$w(\Gamma\cap{\rm Int}D')=0$.  
\end{enumerate}
By Lemma~\ref{NiceEdge7AngledDisk},
the disk $D$ has at most two feelers.

If $D$ has exactly one feeler at $w_4$,
then the disk $D'$ contains the three terminal edges of label $m+1$ containing $w_5,w_6,w_7$.
Let $a_{ii},b_{ii}$ $(i=5,6,7)$ be internal edges 
(possibly terminal edges) of label $m$ at $w_i$ in $D'$
such that
$a_{55},b_{55},a_{66},b_{66}$ are oriented outward at
$w_5,w_5,w_6,w_6$, respectively, and
$a_{77},b_{77}$ are oriented inward at $w_7$.
We can show that none of the six edges are terminal edges.
Thus by IO-Calculation with respect to $\Gamma_m$ in $D'$,
we have $w(\Gamma\cap {\rm Int}D')\ge1$.
This contradicts (1).
Hence the disk $D$ has exactly two feelers.

By Lemma~\ref{NiceEdge7AngledDisk},
a regular neighborhood of $D$ contains the pseudo chart as shown in Fig.~\ref{fig25}(b)
(see Fig.~\ref{fig38}(b)).
Let $e_6$ be the terminal edge of label $m+1$ at $w_6$.
By C-II moves, we can move the black vertex in $e_6$
near the white vertex $w_5$ along the internal edge of
label $m$ connecting $w_5$ and $w_6$.
By applying a C-I-M2 move between $e_6$ and
the internal edge $e_2$ of label $m+1$,
we obtain a new terminal edge of label $m+1$ at $w_2$.
Thus we obtain a minimal chart of type $(m;4,3)$
with $w(\Gamma\cap{\rm Int}D_1)=1$ and $e_2'=e_3'$,
but the edge $e_2$ is a terminal edge.
This contradicts Corollary~\ref{D1OneDNotMinimal3}.
Therefore we complete the proof of 
Lemma~\ref{D1OneD2ThreeD3ZeroNotFig14B}. 
\end{Proof}

%%%%%%%%%%%%%%%%%%
%%%%%%%%%%%%%%%%%% Figure
%%%%%%%%%%%%%%%%%%
%%%%%%%%%%%%%%%%%%%%
%%%%%%%%%%%%%%%%%%%%
%%%%%%%%%%%%%%%%%%%%

\begin{lemma}
\label{D1OneD2ThreeD3ZeroNotFig14C}
There does not exist a minimal chart $\Gamma$ of type $(m;4,3)$
with $w(\Gamma\cap{\rm Int}D_1)=1$,
$w(\Gamma\cap{\rm Int}D_2)=3$, 
$w(\Gamma\cap{\rm Int}D_3)=0$ such that
a regular neighborhood of $D_2$
contains the pseudo chart
as shown in Fig.~\ref{fig15}$($c$)$.
\end{lemma}

\begin{figure}
\centerline{\includegraphics{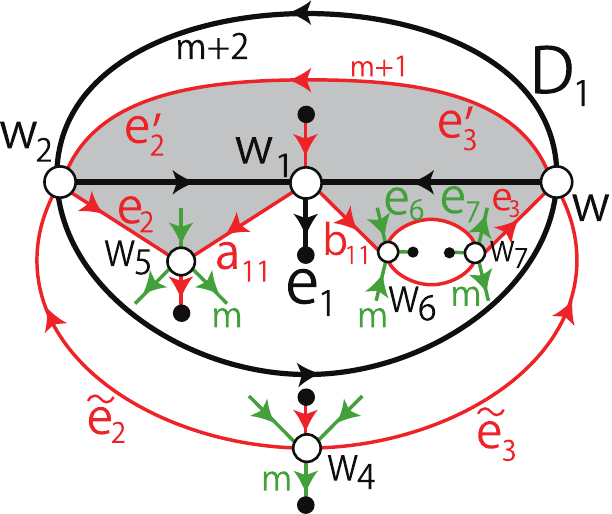}}
\caption{\label{fig39}
The gray region is the 6-angled disk $D$.
}
\end{figure}

\begin{Proof}
Suppose that there exists a minimal chart $\Gamma$ of 
type $(m;4,3)$
with $w(\Gamma\cap{\rm Int}D_1)=1$,
$w(\Gamma\cap{\rm Int}D_2)=3$, 
$w(\Gamma\cap{\rm Int}D_3)=0$ such that
a regular neighborhood of $D_2$
contains the pseudo chart
as shown in 
Fig.~\ref{fig15}(c).
By Lemma~\ref{D1OneD2ThreeD3ZeroPseudoChart},
the chart $\Gamma$ contains the pseudo chart 
as shown in Fig.~\ref{fig37}(a).

Let $D$ be the special 6-angled disk of $\Gamma_{m+1}$
with $\partial D\supset e_2'\cup e_2\cup e_3\cup a_{11}\cup b_{11}$.
Then $w(\Gamma\cap{\rm Int}D)=0$.
Let $w_5=e_2\cap a_{11}$ be the white vertex,
and $e_5$ the terminal edge of label $m+1$ at $w_5$.
Then $e_5\not\subset D$ by
IO-Calculation with respect to $\Gamma_m$ in $D$.
Thus the chart $\Gamma$ contains the pseudo chart 
as shown in Fig.~\ref{fig39}.
We use the notations as shown in Fig.~\ref{fig39}.

Let $e_6,e_7$ be internal edges (possibly terminal edges)
of label $m$ at $w_6,w_7$ in $D$, respectively.
Then  neither $e_6$ nor $e_7$ is middle at 
$w_6$ or $w_7$.
Thus by Assumption~\ref{AssumeTerminal},
neither $e_6$ nor $e_7$ is a terminal edge.
Hence we have $e_6=e_7$.
Hence there exists a lens $E$ with $\partial E\supset e_6$.
This contradicts Lemma~\ref{NoLens}.
Therefore we complete the proof of 
Lemma~\ref{D1OneD2ThreeD3ZeroNotFig14C}. 
\end{Proof}

%%%%%%%%%%%%%%%%%%
%%%%%%%%%%%%%%%%%% Figure
%%%%%%%%%%%%%%%%%%

%%%%%%%%%%%%%%%%%%%%%%%

{\it Proof of the main theorem$($Theorem~\ref{MainTheorem}$)$.}
Suppose that there exists a minimal chart of type $(m;4,3)$.
Then $\Gamma$ contains the pseudo chart as shown in
Fig.~\ref{fig28} (see Section~\ref{s:Chart43}).
By Table~1 in Section~\ref{s:Chart43},
either $w(\Gamma\cap{\rm Int}D_1)=1$ or $3$.
Thus by Proposition~\ref{D1ThreeD2OneD3Zero},
we have $w(\Gamma\cap{\rm Int}D_1)=1$.
Hence by Table~1,
$w(\Gamma\cap{\rm Int}D_2)=1,2$ or $3$.
Thus by Proposition~\ref{D1OneD2OneD3Two} and
Proposition~\ref{D1OneD2TwoD3One},
we have $w(\Gamma\cap{\rm Int}D_2)=3$.
Therefore by Table~1, we have
\begin{enumerate}
\item[(1)] $w(\Gamma\cap{\rm Int}D_1)=1$, $w(\Gamma\cap{\rm Int}D_2)=3$, $w(\Gamma\cap{\rm Int}D_3)=0$.
\end{enumerate}

By Lemma~\ref{D1OneD2ThreeD3ZeroPseudoChart},
the chart $\Gamma$ contains the pseudo chart as shown in
Fig.~\ref{fig37}(a).
Moreover by Lemma~\ref{D1OneD2ThreeD3ZeroPseudoChart}(a),
the edge $e_2$ is middle at $w_2$ but not a terminal edge.
Hence by Condition~(d) in Section~\ref{s:Chart43} and
by four lemmata 
(Lemma~\ref{3angledDiskOneFeelerThreeWhiteVertex},
Lemma~\ref{D1OneD2ThreeD3ZeroNotFig14A},
Lemma~\ref{D1OneD2ThreeD3ZeroNotFig14B}, 
Lemma~\ref{D1OneD2ThreeD3ZeroNotFig14C})
there exists a simple arc $L$ in $(\Gamma_{m+1}-w_1)\cap D_2$ connecting $w_2$ and $w_3$.
Furthermore by Lemma~\ref{D1OneD2ThreeD3ZeroPseudoChart}(a),
 the edge $e_2$ contains a white vertex
in ${\rm Int}D_2$.
Thus by the condition $w(\Gamma\cap{\rm Int}D_2)=3$,
 there are three cases:
(i) $w({\rm Int}L)=1$,
(ii) $w({\rm Int}L)=2$,
(iii) $w({\rm Int}L)=3$.

{\bf Case (i).}
Since $w({\rm Int}L)=1$,
the intersection $e_2\cap e_3$ is one white vertex.
However such a chart does not exist 
by Lemma~\ref{D1OneD4Angled2}.
Hence Case (i) does not occur.

{\bf Case (ii).}
Let $w_5,w_6,w_7$ be the white vertices in ${\rm Int}D_2$
with $w_5\in e_2$ and $w_6\in e_3$.
The arc $L$ divides $D_2$ into two disks.
One of the two disks contains the white vertex $w_1$, 
say $E_1$.
Let $E_2$ be the other disk.
Then (ii-1) $w_7\in E_1$ or (ii-2) $w_7\in E_2$.

\begin{figure}
\centerline{\includegraphics{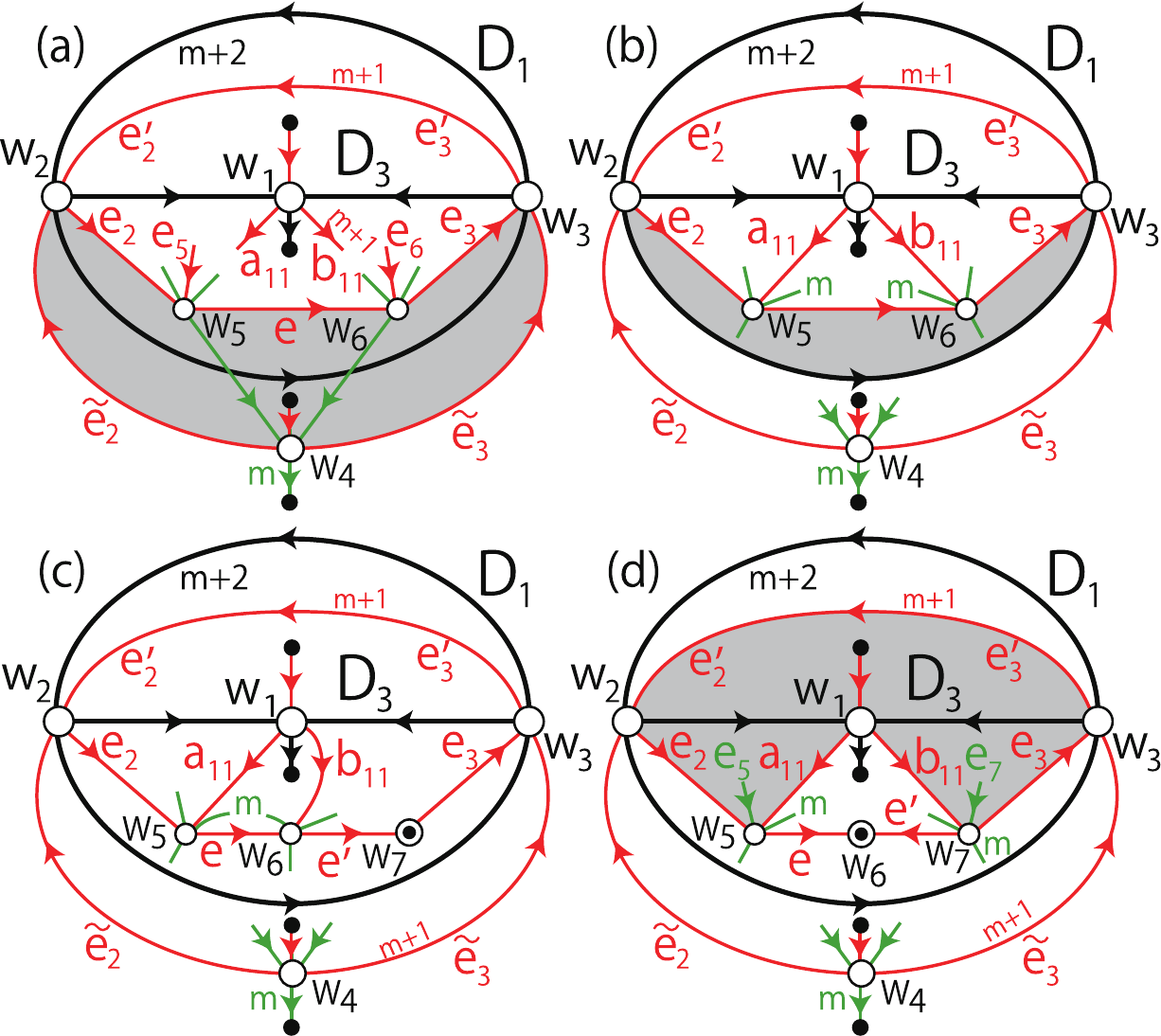}}
\caption{\label{fig40}
(a) The gray region is the 5-angled disk $D$.
(b) The gray region is the disk $E_2$.
(c) $w_7$ is the BW-vertex with respect to $\Gamma_{m+1}$.
(d) The gray region is the 5-angled disk $D$.
}
\end{figure}

{\bf Case (ii-1).}
Let $D$ be the special 5-angled disk of $\Gamma_{m+1}$
with at least one feeler at $w_4$ such that
$\partial D\supset e_2\cup e_3\cup\widetilde{e_2}\cup \widetilde{e_3}$.
Since $w_7\in E_1$, we have $w(\Gamma\cap {\rm Int}D)=0$.
Thus by Lemma~\ref{NiceEdge5AngledDiskAtMostOneFeeler} and Lemma~\ref{NiceEdge5AngledDiskOneFeeler},
a regular neighborhood of $D$ contains the pseudo chart as
shown in 
Fig.~\ref{fig24}(c)
(see Fig.~\ref{fig40}(a)).

Let $e_5,e_6$ be internal edges (possibly terminal edges)
of label $m+1$ at $w_5,w_6$ in $E_1$,
respectively, and
$e$ the internal edge of label $m+1$
connecting $w_5$ and $w_6$.
Then the edges $e_5,e_6$ are oriented inward at $w_5,w_6$,
respectively.
Since both of $e_2$ and $e_5$ are oriented inward at $w_5$,
\begin{enumerate}
\item[(1)] 
the edge $e_5$ is not middle at $w_5$, and
\end{enumerate}
the edge $e$ is oriented from $w_5$ to $w_6$.
Since both of $e$ and $e_6$ are oriented inward at $w_6$,
\begin{enumerate}
\item[(2)] 
the edge $e_6$ is not middle at $w_6$.
\end{enumerate}
Thus by (1),(2) and Assumption~\ref{AssumeTerminal},
neither $e_5$ nor $e_6$ is a terminal edge.
Moreover
by Lemma~\ref{D1OneD2ThreeD3ZeroPseudoChart}(b),
neither $a_{11}$ nor $b_{11}$ is a terminal edge.
Thus the condition $w(\Gamma\cap {\rm Int}E_1)\ge1$
implies $w(\Gamma\cap {\rm Int}E_1)\ge2$ 
by IO-Calculation with respect to $\Gamma_{m+1}$ in $E_1$.
Hence 
$$3=w(\Gamma\cap{\rm Int}D_2)\ge w(\Gamma\cap {\rm Int}E_1)+w(\{w_5,w_6\})\ge2+2=4.$$
This is a contradiction.
Thus Case (ii-1) does not occur.

{\bf Case (ii-2).}
Since $w_7\in E_2$,
we have $w(\Gamma\cap{\rm Int}E_1)=0$.
Thus by Lemma~\ref{D1OneD2ThreeD3ZeroPseudoChart}(b),
we have $a_{11}\ni w_5$ and $b_{11}\ni w_6$
(see Fig.~\ref{fig40}(b)).
Hence there exists a connected component of $\Gamma_{m+1}$
with exactly one white vertex $w_7$ in $E_2$.
This contradicts Lemma~\ref{LemmaWithTerminal3}.
Hence Case (ii-2) does not occur.
Therefore Case (ii) does not occur.

{\bf Case (iii).}
Let $w_5,w_6,w_7$ be the white vertices in ${\rm Int}D_2$
with $w_5\in e_2$ and $w_7\in e_3$
(see Fig.~\ref{fig40}(c),(d)).
Let $e,e'$ be internal edges of label $m+1$ with
$\partial e=\{w_5,w_6\}$ and
$\partial e'=\{w_6,w_7\}$.
By Lemma~\ref{D1OneD2ThreeD3ZeroPseudoChart}(b),
there are three cases:
(iii-1) $a_{11}\ni w_5$, $b_{11}\ni w_6$
(see Fig.~\ref{fig40}(c)),
(iii-2) $a_{11}\ni w_5$, $b_{11}\ni w_7$
(see Fig.~\ref{fig40}(d)),
(iii-3) $a_{11}\ni w_6$, $b_{11}\ni w_7$.

{\bf Case (iii-1) and Case (iii-3).}
We shall show only that Case (iii-1) does not occur.
Suppose that 
$a_{11}\ni w_5$, $b_{11}\ni w_6$.
Then there exists a terminal edge of label $m+1$ at $w_7$.

Since $e_2,a_{11}$ are oriented inward at $w_5$,
the edge $e$ is oriented outward at $w_5$.
Thus the edge $e$ is 
oriented inward at $w_6$.
Since the edge $b_{11}$ is oriented inward at $w_6$,
the edge $e'$ is oriented from $w_6$ to $w_7$
(see Fig.~\ref{fig40}(c)).
Moreover since $e_3$ is oriented outward at $w_7$,
the terminal edge of label $m+1$ at $w_7$
is not middle at $w_7$.
This contradicts Assumption~\ref{AssumeTerminal}.
Hence Case (iii-1) does not occur.

Similarly we can show that Case (iii-3) does not occur.

{\bf Case (iii-2).}
Since $a_{11}\ni w_5$ and $b_{11}\ni w_7$,
 the white vertex $w_6$ is contained in a terminal edge
of label $m+1$.
Let $D$ be the 5-angled disk of $\Gamma_{m+1}$
with $\partial D=e_2'\cup e_2\cup e_3\cup a_{11}\cup b_{11}$ and $w(\Gamma\cap{\rm Int}D)=0$.
Let $e_5,e_7$ be internal edges (possibly terminal edges)
of label $m$ at $w_5,w_7$ in $D$, respectively.
Since both of $e_2$ and $a_{11}$ are oriented 
inward at $w_5$,
\begin{enumerate}
\item[(3)] $e_5$ is oriented inward at $w_5$, and 
\end{enumerate}
the edge $e$ is oriented from $w_5$ to $w_6$.
Thus by Lemma~\ref{OriBWvertex}
the edge $e'$ is oriented from $w_7$ to $w_6$.
Hence 
\begin{enumerate}
\item[(4)] $e_7$ is oriented inward at $w_7$ 
(see Fig.~\ref{fig40}(d)). 
\end{enumerate}
Since $e_3$ is oriented outward at $w_7$,
the edge $e_7$ is not middle at $w_7$.
Thus by Assumption~\ref{AssumeTerminal},
the edge $e_7$ is not a terminal edge.
By (3),(4) and by IO-Calculation with respect to $\Gamma_m$
in $D$, we have $w(\Gamma\cap{\rm Int}D)\ge1$.
This contradicts $w(\Gamma\cap{\rm Int}D)=0$.
Thus Case (iii-2) does not occur.

Therefore all the three cases do not occur.
Hence there does not exist a minimal chart
of type $(m;4,3)$.
Thus we complete the proof of Theorem~\ref{MainTheorem}.
{\hfill {$\square$}\vspace{1.5em}}

%%%%%%%%%%%%%%%%%%
%%%%%%%%%%%%%%%%%% Figure
%%%%%%%%%%%%%%%%%%

%%%%%%%%%%%%%%%%%%%%%
%%%%%%%%%%%%%%%%%%%%%

%%%%%%%%%%%%%%%%%%%%%%%%%%
%%%%%%%%%%%%%%%%%%%%%%%%%%
%%%%%%%%%%%%%%%%%%%%%%%%%%

%%\newpage

%%%%%%%%%%%%%%%%%%%%%%%%%%
%%%   Reference
%%%%%%%%%%%%%%%%%%%%%%%%%%

%%%%%%%%%%%%%%%%%%%%%%
%%%%%%%%%%%%%%%%%%%%%%%
%%%%%%%%%%%%%%%%%%%%%%%

\vspace{5mm}

\begin{minipage}{65mm}
{Teruo NAGASE
\\
{\small Tokai University \\
4-1-1 Kitakaname, Hiratuka \\
Kanagawa, 259-1292 Japan\\
\\
nagase@keyaki.cc.u-tokai.ac.jp
}}
\end{minipage}
\begin{minipage}{65mm}
{Akiko SHIMA 
\\
{\small Department of Mathematics, 
\\
Tokai University
\\
4-1-1 Kitakaname, Hiratuka \\
Kanagawa, 259-1292 Japan\\
shima@keyaki.cc.u-tokai.ac.jp
}}
\end{minipage}

%%\newpage

\vspace{0.7cm}

{\bf List of terminologies}\vspace{2mm}\\
{\small $
\begin{array}{ll||}
\text{$k$-angled disk} & p6 \\
\text{admissible disk} & p5\\
\text{associated disk of a loop} & p5\\
\text{BW-vertex} & p7 \\
\text{C-move equivalent} & p3 \\
\text{chart} & p3 \\
\text{complexity $(w(\Gamma),-f(\Gamma))$} & p3 \\
\text{$(D,\alpha)$-arc} & p25\\
\text{$(D,\alpha)$-arc free} & p25\\
\text{feeler} & p6 \\
\text{free edge} & p3 \\
\text{hoop} & p4 \\
\text{internal edge} & p6 \\
\text{inward} & p3 \\
\text{inward arc} & p23 \\
\text{IO-Calculation} & p24 \\
\text{keeping $X$ fixed} & p25 \\
\text{lens} & p15 \\
\end{array}
~~
\begin{array}{ll}
\text{loop} & p5 \\
\text{M4-psuedo chart} & p30\\
\text{middle arc} & p3 \\
\text{middle at $v$} & p3 \\
\text{minimal chart} & p3 \\
\text{nice edge} & p28\\
\text{outward} & p3 \\
\text{outward arc} & p23 \\
\text{point at infinity $\infty$} & p4 \\
\text{proper arc} & p7\\
\text{pseudo chart} & p6 \\
\text{ring} & p4 \\
\text{RO-family} & p9 \\
\text{simple hoop} & p4 \\
\text{special $k$-angled disk} & p9 \\
\text{terminal edge} & p4 \\
\text{type $(m;n_1,n_2,\cdots,n_k)$ for a chart} & p2 \\
& \\
\end{array}
$}

\vspace{0.5cm}

{\bf List of notations}\vspace{2mm}\\
{\small $
\begin{array}{ll||}
\text{$\Gamma_m$} & p2 \\
\text{$w(\Gamma)$} & p3 \\
\text{$f(\Gamma)$} & p3 \\
\text{${\rm Int}X$} & p4 \\
\text{$\partial X$} & p4 \\
\end{array}
$~~
$\begin{array}{ll}
\text{$Cl(X)$} & p4 \\
\text{$\partial \alpha$} & p4 \\
\text{${\rm Int}\alpha$} & p4\\
\text{$w(X)$} & p5 \\
\text{$a_{ij},b_{ij}$} & p9 \\
\end{array}
$
}

\end{document}